\documentclass{article}

\usepackage{arxiv}

\usepackage[utf8]{inputenc}
\usepackage[T1]{fontenc}
\usepackage{url}
\usepackage{microtype}
\usepackage[square,sort,comma,numbers]{natbib}

\usepackage{placeins}
\usepackage{amsmath,amssymb}
\usepackage{graphicx}
\usepackage{subcaption}
\usepackage{booktabs}
\usepackage{multirow}
\usepackage{float}
\usepackage{xcolor}
\usepackage{tikz}
\usepackage{hyperref}
\usepackage[capitalise,noabbrev]{cleveref}

\graphicspath{{figures/}}

\usepackage{amsthm}
\newtheorem{proposition}{Proposition}

\theoremstyle{remark}
\newtheorem{remark}{Remark}
\theoremstyle{plain}

\newcommand{\SB}{SB}
\newcommand{\RODE}{ROD-E}
\newcommand{\RODL}{\texorpdfstring{ROD-$L^2$}{ROD-L2}}
\newcommand{\RODESP}{ROD-E-SP}
\newcommand{\RODLSP}{\texorpdfstring{ROD-$L^2$-SP}{ROD-L2-SP}}
\newcommand{\RODEMP}{ROD-E-MP}
\newcommand{\RODLMP}{\texorpdfstring{ROD-$L^2$-MP}{ROD-L2-MP}}
\newcommand{\PP}[1]{\mathbb{P}^{#1}}

\newcommand{\bphi}{\boldsymbol{\varphi}}
\newcommand{\bpsi}{\boldsymbol{\psi}}
\newcommand{\bu}{\mathbf{u}}
\newcommand{\bv}{\mathbf{v}}
\newcommand{\bx}{\mathbf{x}}
\newcommand{\bxb}{\bar{\mathbf{x}}}
\newcommand{\bxt}{\tilde{\mathbf{x}}}
\newcommand{\ba}{\mathbf{a}}
\newcommand{\bn}{\mathbf{n}}
\newcommand{\btaub}{\bar{\boldsymbol\tau}}
\newcommand{\bd}{\mathbf{d}}

\newcommand{\balpha}{\boldsymbol{\alpha}}
\newcommand{\blambda}{\boldsymbol{\lambda}}
\newcommand{\uD}{\mathbf{u}_D}
\newcommand{\MM}{\mathbf{M}}
\newcommand{\KK}{\mathbf{L}}
\newcommand{\Phib}{\Phi(\bxb)}
\newcommand{\phit}{\tilde{\bphi}}

\newcommand{\CFLmax}{\mathrm{CFL}^p_{\max}}
\newcommand{\aSPE}{\alpha^{\text{SP-E}}}
\newcommand{\aSPL}{\alpha^{\text{SP-}L^2}}
\newcommand{\aSPW}{\alpha^{\text{SP-W}}}
\newcommand{\aMPE}{\balpha^{\text{MP-E}}}
\newcommand{\aMPL}{\balpha^{\text{MP-}L^2}}
\newcommand{\aMPW}{\balpha^{\text{MP-W}}}
\newcommand{\spanrows}[1]{\multicolumn{2}{c}{\multirow{4}{*}{#1}}}
\newcommand{\spanskip}{\multicolumn{2}{c}{}}
\newcommand{\spancell}[1]{\multicolumn{2}{c}{#1}}

\title{High-order polynomial corrections for embedded curved boundaries in discontinuous Galerkin methods for hyperbolic systems:
from staircase instabilities to a gap method}
\author{Mirco Ciallella \\
	Laboratoire Jacques-Louis Lions\\
	Universit\'e Paris Cit\'e\\
	75013, Paris \\
	\texttt{mirco.ciallella@u-pariscite.fr}
}

\renewcommand{\headeright}{Ciallella (2026)}
\renewcommand{\shorttitle}{An \textit{arXiv} Preprint}

\hypersetup{
pdftitle={High-order polynomial corrections for embedded curved boundaries in discontinuous Galerkin methods for hyperbolic systems: from staircase instabilities to a gap method},
pdfauthor={Mirco Ciallella},
}

\begin{document}
\maketitle

\begin{abstract}
High order polynomial corrections, which impose the boundary conditions on an embedded boundary by modifying the polynomial of the surrogate boundary cells, give accurate and robust discretizations of elliptic problems on unfitted meshes. We show that, used out of the box in discontinuous Galerkin discretizations of hyperbolic systems on Cartesian meshes, the shifted boundary and minimization-based (ROD) corrections become unstable, with eigenvalues of positive real part that no restriction of the time step can remove. Linear advection, whose operator is block triangular, is the least affected case: a tilted straight boundary makes every correction unstable with the geometry-normal map, while the face-orthogonal map, which likely gives a one-dimensional character to the boundary condition with tensor-product bases, keeps the 1D stability up to one polynomial degree less. The slip wall of the linear acoustic system is much more critical: all the corrections are unstable already on a tilted straight wall, and an energy balance shows that the cause is the misalignment, in direction and in position, between the faces of the staircase and the wall. Since the origin of the instability is geometric, we act on the surrogate boundary itself, and develop the gap paradigm for discontinuous Galerkin discretizations of hyperbolic systems: each boundary cell is extended, by extrapolating its polynomial without new degrees of freedom, up to the segments joining the projections of its vertices on the wall, where the boundary condition is imposed, and neighbouring extended cells exchange information through the numerical fluxes used on the interior interfaces. Numerical experiments confirm the energy stability of the gap method on a straight wall and show high order convergence on curved walls, in configurations where the classical corrections applied on the staircase diverge.

\end{abstract}



\section{Introduction}
\label{sec:intro}

Thanks to their compact stencil and to their high order of accuracy, discontinuous Galerkin (DG) methods~\cite{reed1973triangular,cockburn2001runge} are a natural choice to solve hyperbolic problems accurately at an affordable computational cost~\cite{wang2013high}. When the geometry of the domain is complex, however, the accuracy of the numerical solution is intrinsically linked to the accuracy of its representation~\cite{bassi1997high}. The isoparametric approach, which describes curved boundaries with high order polynomial mappings, recovers the optimal convergence rates, but the generation of valid high order curvilinear meshes, with nonlinear mappings, complex quadrature rules and mesh quality control, remains a difficult task for realistic applications~\cite{persson2009curved,moxey2015isoparametric,coppeans2025aerodynamic}.

Embedded, or immersed, boundary methods follow a different paradigm: the domain is discretized with a simple fixed mesh, possibly Cartesian, over which the boundaries are free to move, so that the difficulty is shifted from the generation of the mesh to the formulation of accurate and stable boundary conditions. Since the pioneering work of Peskin on the flow inside the heart~\cite{peskin1977numerical,peskin2002immersed}, a large number of embedded methods has been developed, among which we only recall the immersed boundary and ghost-point methods in finite differences~\cite{mittal2005immersed,fadlun2000combined,fedkiw1999non}, the cut-cell finite volume methods~\cite{berger2012progress,berger2017cut,monasse2012conservative}, and the cut finite element methods~\cite{hansbo2002unfitted,burman2015cutfem,burman2025cut}, also in the DG context~\cite{qin2013discontinuous,giuliani2022two,adjerid2019immersed}. These approaches modify the integration domain of the cells cut by the boundary, and must deal with arbitrarily small cut cells, which restrict the time step or spoil the conditioning, and are handled by cell merging~\cite{BERGER2015180,may2017explicit}, by stabilization terms such as the ghost penalty~\cite{burman2010ghost}, and by careful quadrature on the cut cells.

A different family of methods avoids the cut cells altogether. The cells crossed by the boundary are discarded, and the boundary conditions are imposed on the surrogate boundary formed by the faces of the remaining cells, after a modification that compensates the distance, of the order of the mesh size, between the surrogate and the true boundary. In the shifted boundary method (SBM)~\cite{main2018shifted}, the condition is transferred with a truncated Taylor expansion, and the method has been extended to several problems~\cite{song2018shifted,ZENG2025113512} and to high order~\cite{atallah2022high,collins2023penalty}. A Taylor expansion at the boundary is also the basis of the inverse Lax--Wendroff (ILW) procedure~\cite{tan2010inverse,tan2012efficient}, developed for finite differences on Cartesian meshes, in which the values at the ghost points are extrapolated from the boundary and the normal derivatives are obtained from the boundary condition and from the time and tangential derivatives through the equation itself; its stability with explicit time integration was analysed in~\cite{li2016stability}, and it has recently been extended to DG methods~\cite{yang2025inverse}. In the reconstruction for off-site data (ROD)~\cite{costa2018very,costa2019very}, which extends to unfitted boundaries the constrained least-squares reconstructions developed for finite volume methods on body-fitted meshes~\cite{ollivier2002high,michalak2009unstructured}, a modified polynomial of the boundary cell is obtained by a constrained minimization that enforces the condition at points of the true boundary. The approach has also been developed for DG methods on polygonal meshes~\cite{santos2024very,ciallella2024very}. The shifted boundary polynomial correction~\cite{ciallella2023shifted} was proposed for polygonal meshes as a simple high order version of the high order SBM~\cite{atallah2022high}: all the derivatives of the Taylor expansion are replaced by the evaluation of the polynomial of the boundary cell at the true boundary, so that the modified condition becomes a correction of the polynomial that requires no derivative. It was then applied to high order arbitrary-Lagrangian-Eulerian methods for hyperbolic equations with moving curved boundaries~\cite{BOSCHERI2025114215} and to spectral elements for the Poisson problem~\cite{visbech2025spectral}. The stability of DG methods with this correction, in the fully unfitted configurations in which the distance is of the order of the mesh size, was studied in~\cite{ciallella2027stability} for 1D linear hyperbolic equations, which showed that the correction restricts the stability region of explicit schemes more and more as the degree increases. In~\cite{ciallella2026minimization} it was shown, from a 1D analysis, that the ROD minimization can also be written as a polynomial correction, which differs from the SBM one only by a scalar coefficient, and that this coefficient makes the ROD corrections, in their Euclidean (ROD-E) and $L^2$ (ROD-$L^2$) versions, stable under the standard CFL condition up to degree three and four, respectively. These minimization-based polynomial corrections have been validated to reach high order on both fixed and moving domains discretized with polygonal meshes~\cite{ciallellaboscheri2026minimization}. They have then been generalized to Neumann and Robin conditions for the Poisson problem on unfitted spectral elements~\cite{ciallella2026generalized}, where they give accurate and robust discretizations with improved conditioning with respect to the original correction.

The goal of this work is to continue the stability study of high order embedded DG methods for hyperbolic equations started in~\cite{ciallella2027stability,ciallella2026minimization}, by moving from 1D to 2D Cartesian meshes. Motivated by failing numerical experiments, we show that, unfortunately, the 1D stability results extend to 2D only when the distance between the surrogate and the true boundary is constant along the boundary, and do not extend to general boundaries. For linear advection, the operator of the mesh is block triangular, so that the stability is decided cell by cell. We show that already a single straight boundary tilted with respect to the mesh is enough to make every correction unstable when the boundary data are taken at the closest point of the true boundary, while, if each point of a surrogate face is mapped along the normal of the face, the corrections keep the stability of the 1D problem up to one polynomial degree less. This is probably related to the tensor-product bases used in this work, which give this map a one-dimensional character. The impact is even larger for the slip wall boundary condition of the linear acoustic system, used as a preliminary model for the nonlinear Euler equations for gas dynamics. Unlike linear advection, the operator of the acoustic system is not block triangular, every cell exchanges information with all its neighbours, and the stability of a wall boundary treatment depends on the whole mesh. All the corrections are unstable already on a tilted straight wall. A simple energy balance of the surrogate wall flux shows that two terms can produce energy: a normal correction term, which appears as soon as the normal of the true wall is used on a face that is not aligned with it, as in the boundary correction developed for body-fitted meshes in~\cite{krivodonova2006high}, and a location correction term, introduced by the polynomial correction that moves the condition from the face to the wall. Both terms are $\mathcal{O}(1)$, i.e.\ they do not vanish as the mesh is refined, on the surrogate boundaries obtained by removing cells from a mesh that is generated independently of the geometry, whose faces are not aligned with the wall and lie at a distance of the order of the mesh size from it. The development of stable high order polynomial corrections for hyperbolic systems that use the faces of a Cartesian mesh as surrogate boundary therefore remains an open question.

The energy balance shows that the terms that can produce energy depend on the geometry of the surrogate boundary, through its angle with the wall and its distance from it relative to the cell size, more than on the correction itself. Therefore, instead of designing a new correction, we propose to stabilize the problem by changing the geometry of the surrogate boundary, and we develop the gap paradigm in the DG framework. The Gap-SBM was introduced for continuous Galerkin discretizations of elliptic operators, to recover the optimal accuracy with Neumann boundary conditions~\cite{collins2026gap}, and was then extended to isogeometric discretizations~\cite{gorgi2026isogeometric,antonelli2026isogeometric} and to solid mechanics~\cite{li2026gap}. Each element with a face on the surrogate boundary is extended by the quadrilateral between this face and the segment joining the closest-point projections of its vertices on the true boundary, so that the boundary is approximated by a polygon with vertices on it, and the polynomial of the element is used on the extension without new degrees of freedom. The extensions of two neighbouring elements do not match on the edge that they share in the gap, and in the continuous formulation these edges must be treated with jump and average terms with a penalty, as in DG methods~\cite{li2026gap}. In this work we develop, for the first time, the gap paradigm for high order DG discretizations of hyperbolic systems. The DG framework allows a more natural introduction of the gap cells: an extended cell is still a cell of a broken polynomial space, and the gap cells of neighbouring boundary cells exchange information through the same numerical fluxes that are used on the interior interfaces, with no additional terms or penalty parameters. We show that this formulation gives a much more stable method, which is energy stable for the slip wall of the linear acoustic system on straight walls. The price to pay is the construction of the gap cells and the integration on them, of the extended mass matrix, of the volume terms and of the fluxes on the chords and on the lateral edges, which makes the method more expensive than the standard corrections on the surrogate Cartesian faces. These, however, are of no use when they are unstable. On curved walls the energy production terms do not vanish, but they are reduced from $\mathcal{O}(1)$ to terms proportional to the local curvature, with an angle between the chords and the wall of $\mathcal{O}(\kappa h)$ and a distance of $\mathcal{O}(\kappa h^2)$. They still act where the curvature is poorly resolved by the mesh, which can be mitigated by a local refinement. 

The rest of the paper is organized as follows. In \cref{sec:methods}, the DG discretization on embedded Cartesian meshes and the polynomial corrections are presented. \Cref{sec:advection} studies the stability of the corrections for the linear advection equation, with a single-cell model, a tilted straight boundary, and the numerical experiments on a disk. In \cref{sec:acoustics}, the slip wall of the linear acoustic system is studied, together with its energy balance, which identifies the cause of the instabilities. The gap method is presented in \cref{sec:gap}, with its energy balance and the numerical experiments on a pressure pulse reflected by a straight wall, on the potential flow around a cylinder and on the start-up of the flow around an inclined ellipse. Some conclusions are drawn in \cref{sec:conclusions}.

\section{Discretization, embedded boundaries and polynomial corrections}
\label{sec:methods}

\subsection{Setting and notation}

We consider a hyperbolic system in conservation form
\begin{equation}
  \partial_t\mathbf U+\partial_x\mathbf F(\mathbf U)+\partial_y\mathbf G(\mathbf U)=\mathbf S(\bx)\quad\text{in }\bar\Omega,
  \label{eq:system}
\end{equation}
with $\mathbf U(\bx,t)\in\mathbb R^m$, completed by boundary conditions on the waves that enter the domain: at a point of the boundary with outward normal $\bar\bn$, as many conditions are prescribed as there are negative eigenvalues of the normal Jacobian $\bar n_x\,\partial\mathbf F/\partial\mathbf U+\bar n_y\,\partial\mathbf G/\partial\mathbf U$. Here $\bar\Omega\subset\mathbb{R}^2$ is a curved physical domain with boundary
$\bar\Gamma=\partial\bar\Omega$, the \emph{real} (or physical) boundary, of normal $\bar\bn$. The
domain is not meshed: it is embedded in a Cartesian background mesh of square cells $K_{ij}$ of
size $\Delta x=\Delta y$. The \emph{active} cells are the cells kept in the computation (in this work,
the cells fully contained in $\bar\Omega$), and their union $\Omega=\bigcup_{ij}K_{ij}\subset\bar\Omega$
is the \emph{computational domain}. Its boundary $\tilde\Gamma=\partial\Omega$, the
\emph{computational boundary}, is made of the faces of active cells whose neighbour is not active
(\emph{computational faces}). The boundary condition is not imposed on $\bar\Gamma$, which does not
coincide with any face, but on $\tilde\Gamma$, through a modified condition.

To this purpose, each point $\bxt$ of a computational face is associated with a point $\bxb$ of the real
boundary by the mapping $\bxb=\bxt+\bd$, where $\bd=\bd(\bxt)$ is the distance vector.
Two mappings are considered in this work (\cref{fig:mapping}). In the \emph{geometry-normal} map, which is the closest-point projection and the standard choice in the literature, $\bd$ is parallel to the normal $\bar\bn$ of $\bar\Gamma$. In the \emph{face-orthogonal} map, $\bd$ is parallel to the normal $\tilde\bn$ of the computational face, so that the point moves along the Cartesian direction orthogonal to the face and keeps its coordinate along it. On a straight boundary parallel to the faces both maps coincide and $|\bd|$ is constant.

\subsection{Discontinuous Galerkin discretization of hyperbolic systems}
\label{sec:dg}

The system \eqref{eq:system} is discretized with the DG method on the active cells. In each cell, the numerical solution $\mathbf U_h$ is a polynomial of degree $p$, whose components belong to the space $V_{ij}=\PP{p}(K_{ij})$ of the polynomials of degree at most $p$ in each variable (often denoted $\mathbb{Q}^p$), spanned by a tensor-product basis $\bphi=[\varphi_0,\dots,\varphi_{D-1}]^T$ with $D=(p+1)^2$ degrees of freedom, and $M=\int_{K_{ij}}\bphi\bphi^T\,d\bx$ is the mass matrix of the cell. Multiplying \eqref{eq:system} by a test function $\varphi_m\in V_{ij}$, integrating over $K_{ij}$ and integrating by parts the flux divergence, we obtain the weak formulation
\begin{equation}
  \begin{aligned}
  \int_{K_{ij}}\varphi_m\,\partial_t\mathbf U_h\,d\bx
  &-\int_{K_{ij}}\bigl(\partial_x\varphi_m\,\mathbf F(\mathbf U_h)+\partial_y\varphi_m\,\mathbf G(\mathbf U_h)\bigr)\,d\bx\\
  &+\sum_{f\subset\partial K_{ij}}\int_f\varphi_m\,\hat{\mathbf F}\bigl(\mathbf U_h^-,\mathbf U_h^+,\tilde\bn_f\bigr)\,ds
  =\int_{K_{ij}}\varphi_m\,\mathbf S\,d\bx,
  \end{aligned}
  \label{eq:dg-weak}
\end{equation}
for $m=0,\dots,D-1$, where $\tilde\bn_f$ is the outward normal of the face $f$, $\mathbf U_h^-$ is the trace of the cell on $f$ and $\mathbf U_h^+$ the exterior state. On an interior face, $\mathbf U_h^+$ is the trace of the neighbouring active cell. On a computational face, which lies on $\tilde\Gamma$, the neighbouring cell is not active and the exterior state is a ghost state, $\mathbf U_h^+=\mathbf U^g$, through which the boundary condition is imposed weakly. In an embedded method the boundary integrals of \eqref{eq:dg-weak} are therefore computed on the computational boundary $\tilde\Gamma$, and not on the real boundary $\bar\Gamma$, where the condition is prescribed. The integrals are computed with the tensor-product Gauss rule with $(p+1)^2$ points in the cell and with $N_q=p+1$ Gauss points $\bxt_k$ on each face.

The numerical flux is the local Lax--Friedrichs, or Rusanov, flux
\begin{equation}
  \hat{\mathbf F}\bigl(\mathbf U^-,\mathbf U^+,\bn\bigr)=\tfrac12\bigl(\mathbf F_{\bn}(\mathbf U^-)+\mathbf F_{\bn}(\mathbf U^+)\bigr)-\tfrac12\,s_{\max}\bigl(\mathbf U^+-\mathbf U^-\bigr),\qquad
  \mathbf F_{\bn}=n_x\mathbf F+n_y\mathbf G,
  \label{eq:llf}
\end{equation}
where $s_{\max}$ is the largest absolute value of the eigenvalues of the normal Jacobian $\partial\mathbf F_{\bn}/\partial\mathbf U$. For the linear advection equation of \cref{sec:advection}, $s_{\max}=|\ba\cdot\bn|$ and \eqref{eq:llf} is the upwind flux, while for the linear acoustic system of \cref{sec:acoustics}, with unit sound speed, $s_{\max}=1$. The semi-discrete system obtained from \eqref{eq:dg-weak} is advanced in time with an explicit scheme of order $p+1$: the classical Runge--Kutta methods of order $2$, $3$ and $4$ for $p=1,2,3$, and the deferred correction (DeC) method of order $5$~\cite{dutt2000spectral,micalizzi2024new} for $p=4$.

On a body-fitted mesh, the ghost state is built from the data at the face: $u^g=u_D$ for a Dirichlet condition, or the mirror state for a slip wall. On $\tilde\Gamma$ the data are known only on $\bar\Gamma$, at a distance of the order of the mesh size from the face, and this is the role of the corrections. All the methods considered here are \emph{polynomial corrections}: the cell polynomial of a quantity prescribed at the boundary is replaced by a modified polynomial that satisfies the boundary condition on $\bar\Gamma$, and the ghost state at each Gauss point $\bxt_k$ of the computational face is built from its value there. For clarity, in the rest of this section we present the corrections for one scalar quantity $u_h=\bphi^T\bu$, with $\bu$ its coefficients, subject to the Dirichlet condition $u=u_D$: the cell polynomial $u_h$ is replaced by $v_h=\bphi^T\bv\in\PP{p}(K_{ij})$, and the exterior trace is $u_h^+=v_h(\bxt_k)$. This is the case of linear advection, and of each prescribed component of a system; the slip wall, in which the corrected quantity is the normal velocity, is described in \cref{sec:wall-condition}.
In the following, $\bphi(\bxt_k)$ is the basis at the computational boundary point and
$\bphi(\bxb_k)$ the basis at the mapped point, both evaluated with the polynomial of the
cell (extrapolated when $\bxb_k\notin K_{ij}$).

\subsection{The shifted boundary polynomial correction (\SB)}

The \SB\ correction~\cite{ciallella2023shifted} starts from the shifted boundary method~\cite{main2018shifted,atallah2022high}, in which the Dirichlet value at the mapped point is transferred to the computational boundary point with a Taylor expansion along $\bd$,
\begin{equation*}
  v_h(\bxt_k) = u_D(\bxb_k) - \Bigl(\bd\cdot\nabla u_h(\bxt_k) + \tfrac12\,\bd^T\nabla^2u_h(\bxt_k)\,\bd + \dots\Bigr).
\end{equation*}
Since $u_h$ is a polynomial, its Taylor expansion is finite, and when all its terms are kept the sum in brackets is exactly $u_h(\bxb_k)-u_h(\bxt_k)$. The derivatives are then no longer needed, and the correction reads
\begin{equation}
  v_h(\bxt_k) = u_D(\bxb_k) - \big(u_h(\bxb_k)-u_h(\bxt_k)\big)
  = \big[\bphi(\bxt_k)-\bphi(\bxb_k)\big]^T\bu + u_D(\bxb_k).
\end{equation}

\subsection{Single-point (SP) ROD}

The minimization-based corrections (ROD) look for the polynomial closest to $u_h$ that satisfies the
boundary condition at the mapped point. Given a symmetric positive definite weight matrix $W$, this
is the stationary point of the functional
\begin{equation}
  \mathcal{L}(\bv,\lambda) = \tfrac12(\bv-\bu)^TW(\bv-\bu) + \lambda\big(v_h(\bxb_k)-u_D(\bxb_k)\big),
\end{equation}
with $W=I$ for the Euclidean distance (\RODE) and $W=M$, the mass matrix, for the
$L^2$ distance (\RODL). As shown in~\cite{ciallella2026minimization}, solving the
optimality conditions and evaluating $v_h$ at $\bxt_k$ gives the polynomial correction
\begin{equation}
  v_h(\bxt_k) = \big[\bphi(\bxt_k)-\aSPW(\bxt_k)\,\bphi(\bxb_k)\big]^T\bu + \aSPW(\bxt_k)\,u_D(\bxb_k),
  \label{eq:sp}
\end{equation}
\begin{equation}
  \aSPW(\bxt_k)=\frac{\bphi(\bxt_k)^TW^{-1}\bphi(\bxb_k)}{\bphi(\bxb_k)^TW^{-1}\bphi(\bxb_k)}.
  \label{eq:sp-alpha}
\end{equation}
We write $\aSPE$ for $W=I$ and $\aSPL$ for $W=M$.
The \SB\ correction is retrieved if $\aSPW(\bxt_k)=1$. A different modified polynomial is used for
each quadrature point; this is the multi-dimensional correction of~\cite{ciallella2026minimization}.

\subsection{Multi-point (MP) ROD}

The multi-point ROD~\cite{santos2024very,ciallella2024very} imposes the condition at the $N_q$ mapped points
$\{\bxb_k\}_{k=1}^{N_q}$ of a set of computational boundary points with one modified polynomial. With
$[\Phib]_{jk}=\varphi_j(\bxb_k)$ and $\uD=[u_D(\bxb_1),\dots,u_D(\bxb_{N_q})]^T$, the functional
\begin{equation}
  \mathcal{L}(\bv,\blambda) = \tfrac12(\bv-\bu)^TW(\bv-\bu)
  + \sum_{k=1}^{N_q}\lambda_k\big(v_h(\bxb_k)-u_D(\bxb_k)\big)
\end{equation}
has the optimality conditions
\begin{equation}
  \bv = \bu - W^{-1}\Phib\blambda,\qquad
  \blambda = \big(\Phi^T(\bxb)W^{-1}\Phib\big)^{-1}\big(\Phi^T(\bxb)\bu-\uD\big),
\end{equation}
and, evaluating $v_h$ at a computational boundary point $\bxt$, one obtains the vectorial polynomial
correction~\cite{ciallella2026minimization}
\begin{equation}
  v_h(\bxt)=\big[\bphi(\bxt)-\Phib\aMPW(\bxt)\big]^T\bu+[\aMPW(\bxt)]^T\uD,
  \label{eq:mp}
\end{equation}
\begin{equation}
  [\aMPW(\bxt)]^T=\bphi^T(\bxt)W^{-1}\Phib\big(\Phi^T(\bxb)W^{-1}\Phib\big)^{-1}.
  \label{eq:mp-alpha}
\end{equation}
We write $\aMPE$ for $W=I$ and $\aMPL$ for $W=M$.
For $W=M$ this is the ROD2D-$L^2$ correction of~\cite{ciallella2026minimization}. For $N_q=1$ it is the SP ROD. The
formula requires $\Phi^T(\bxb)W^{-1}\Phib$ to be invertible, i.e.\
the $N_q$ conditions to be linearly independent on $\PP{p}$.

\begin{remark}[MP ROD face by face]\label{rem:mp-face}
On linear triangulations~\cite{santos2024very,ciallella2024very} each boundary element has one face on the boundary, and the $N_q=p+1$ Gauss points of that face are unisolvent for the trace. On a Cartesian staircase, instead, a cell can have two or three computational faces. Collecting all their points in one modified polynomial might either produce a matrix $\Phi^T(\bxb)W^{-1}\Phib$ that is singular, or produce an over-constrained system. We therefore use the MP ROD face by face: each computational face has its own modified polynomial with its own $N_q=p+1$ points, and \eqref{eq:mp} is used as it is.
\end{remark}

In the following, SP denotes the single-point ROD and MP the multi-point ROD face by face; e.g.\
\RODESP, \RODLSP, \RODEMP\ and \RODLMP.

\subsection{Equivalence of SP and MP ROD for tensor-product bases}

In the numerical experiments, MP and SP gave the same results to machine precision in some configurations. This is explained by two exact equivalences, which hold for tensor-product bases. Both are written for the left face
$x=0$ of a cell, with $\bphi(x,y)=\bpsi(x)\otimes\bpsi(y)$, $\bpsi$ the 1D basis, and computational
boundary points $\bxt_k=(0,\eta_k)$, $\eta_k$ being the $p+1$ Gauss points of the face. The mapped points
keep the coordinate along the face, $\bxb_k=(\bar x_k,\eta_k)$. We set
$V=[\bpsi(\eta_1)\ \cdots\ \bpsi(\eta_{p+1})]$, which is invertible. Both weights are separable,
$W=W_1\otimes W_1$, with $W_1=I$ or the 1D mass matrix $M_1$.

Both proofs rely on one identity. Since $W^{-1}=W_1^{-1}\otimes W_1^{-1}$, the product of two tensor basis
vectors separates into an $x$-part times a $y$-part:
\begin{equation}
  \big(\bpsi(x)\otimes\bpsi(y)\big)^TW^{-1}\big(\bpsi(x')\otimes\bpsi(y')\big)
  =\big(\bpsi(x)^TW_1^{-1}\bpsi(x')\big)\big(\bpsi(y)^TW_1^{-1}\bpsi(y')\big).
  \label{eq:mp-separation}
\end{equation}

We first apply it to the SP coefficient. The numerator and the denominator of \eqref{eq:sp-alpha} have the
same $y$-part, $\bpsi(\eta_k)^TW_1^{-1}\bpsi(\eta_k)$, which cancels. What remains is the 1D ROD coefficient
for the distance $\bar x_k$:
\begin{equation}
  \aSPW(\bxt_k)=\frac{\bpsi(\bar x_k)^TW_1^{-1}\bpsi(0)}{\bpsi(\bar x_k)^TW_1^{-1}\bpsi(\bar x_k)} .
  \label{eq:sp-1d}
\end{equation}

We then apply it to the MP coefficient. By \eqref{eq:mp-alpha}, $\aMPW(\bxt_k)$ is the solution
of the linear system with matrix $\Phi^T(\bxb)W^{-1}\Phib$ and right-hand side $\Phi^T(\bxb)W^{-1}\bphi(\bxt_k)$.
Since the columns of $\Phib$ are the vectors $\bphi(\bxb_m)$, each entry is a product of the form
\eqref{eq:mp-separation}, so that
\begin{align}
  \big[\Phi^T(\bxb)W^{-1}\Phib\big]_{lm}
  &=\bphi(\bxb_l)^TW^{-1}\bphi(\bxb_m)
   =\big(\bpsi(\bar x_l)^TW_1^{-1}\bpsi(\bar x_m)\big)\big(\bpsi(\eta_l)^TW_1^{-1}\bpsi(\eta_m)\big),
  \label{eq:mp-sys-matrix}\\
  \big[\Phi^T(\bxb)W^{-1}\bphi(\bxt_k)\big]_l
  &=\bphi(\bxb_l)^TW^{-1}\bphi(\bxt_k)
   =\big(\bpsi(\bar x_l)^TW_1^{-1}\bpsi(0)\big)\big(\bpsi(\eta_l)^TW_1^{-1}\bpsi(\eta_k)\big).
  \label{eq:mp-sys-rhs}
\end{align}

Finally, MP and SP coincide as soon as this system has the solution
\begin{equation}
  \aMPW(\bxt_k)=\aSPW(\bxt_k)\,\mathbf{e}_k .
  \label{eq:mp-equals-sp}
\end{equation}
Indeed, since $\Phib\mathbf{e}_k=\bphi(\bxb_k)$, in this case
$\Phib\aMPW(\bxt_k)=\aSPW(\bxt_k)\,\bphi(\bxb_k)$ and only $u_D(\bxb_k)$ remains in the last term of \eqref{eq:mp}, with the SP
coefficient, so that \eqref{eq:mp} reduces to \eqref{eq:sp}. Each of the two proofs below only has to show \eqref{eq:mp-equals-sp}.

\begin{proposition}[straight boundary at constant distance]\label{prop:mp-constant}
Consider a straight boundary, without curvature, parallel to a computational face. Each computational boundary
point is then mapped orthogonally to the face, and the distance $d$ from the face is constant across the points:
$(0,\eta_k)$ is mapped to $(d,\eta_k)$. In this configuration, the MP correction of the face gives the same result as the SP
correction applied at each point, for both \RODE\ and \RODL.
\end{proposition}

\begin{proof}
With $\bar x_l=d$ for all $l$, the $x$-parts in \eqref{eq:mp-sys-matrix} and \eqref{eq:mp-sys-rhs} are the same
for all the entries:
\begin{align*}
  \big[\Phi^T(\bxb)W^{-1}\Phib\big]_{lm}
  &=\big(\bpsi(d)^TW_1^{-1}\bpsi(d)\big)\big(\bpsi(\eta_l)^TW_1^{-1}\bpsi(\eta_m)\big),\\
  \big[\Phi^T(\bxb)W^{-1}\bphi(\bxt_k)\big]_l
  &=\big(\bpsi(d)^TW_1^{-1}\bpsi(0)\big)\big(\bpsi(\eta_l)^TW_1^{-1}\bpsi(\eta_k)\big).
\end{align*}
Since the $y$-parts satisfy
\[
  \big[\bpsi(\eta_l)^TW_1^{-1}\bpsi(\eta_k)\big]_l=\big[\bpsi(\eta_l)^TW_1^{-1}\bpsi(\eta_m)\big]_{lm}\mathbf{e}_k
\]
and the matrix is invertible, even if it is not diagonal, the solution of the system is
\[
  \aMPW(\bxt_k)=\frac{\bpsi(d)^TW_1^{-1}\bpsi(0)}{\bpsi(d)^TW_1^{-1}\bpsi(d)}\,\mathbf{e}_k .
\]
The fraction is the SP coefficient \eqref{eq:sp-1d} with $\bar x_k=d$, so that the condition \eqref{eq:mp-equals-sp} is satisfied and MP coincides with SP. In other words, the conditions mean that the trace of
$v_h$ on the line $x=d$ equals $u_D$, and the minimizer is the 1D ROD in $x$ applied on each line, which is
what the SP ROD does point by point.
\end{proof}

\begin{proposition}[face-orthogonal map]\label{prop:mp-faceorth}
Consider a boundary of arbitrary shape, possibly curved, and the face-orthogonal map: each computational boundary
point moves orthogonally to the face and keeps its coordinate along it, while the distance from the face can change
from point to point: $(0,\eta_k)$ is mapped to $(\bar x_k,\eta_k)$. In this configuration, the MP correction of the face gives
the same result as the SP correction applied at each point, for \RODL. This is not true in general for \RODE.
\end{proposition}

\begin{proof}
Now the $x$-parts change from entry to entry, so the previous argument is not possible. However, with
$W_1=M_1$ the $y$-parts vanish for $l\ne m$: the $p+1$ Gauss points of the face, of weights $w_k$, integrate $M_1$ exactly, so that
$M_1=V\operatorname{diag}(w_k)V^T$ and
\[
  \bpsi(\eta_l)^TM_1^{-1}\bpsi(\eta_m)=\frac{\delta_{lm}}{w_l} .
\]
Therefore, \eqref{eq:mp-sys-matrix} and \eqref{eq:mp-sys-rhs} become
\begin{align*}
  \big[\Phi^T(\bxb)M^{-1}\Phib\big]_{lm}
  &=\big(\bpsi(\bar x_l)^TM_1^{-1}\bpsi(\bar x_m)\big)\,\frac{\delta_{lm}}{w_l},\\
  \big[\Phi^T(\bxb)M^{-1}\bphi(\bxt_k)\big]_l
  &=\big(\bpsi(\bar x_l)^TM_1^{-1}\bpsi(0)\big)\,\frac{\delta_{lk}}{w_l}.
\end{align*}
The matrix is diagonal, so each equation contains only one unknown and can be solved directly. For $l\ne k$ the right-hand side is zero and the
number multiplying the unknown is positive, so that $\big[\aMPL(\bxt_k)\big]_l=0$. For $l=k$ the weight $w_k$ appears on both
sides and cancels, so that
\[
  \big[\aMPL(\bxt_k)\big]_k=\frac{\bpsi(\bar x_k)^TM_1^{-1}\bpsi(0)}{\bpsi(\bar x_k)^TM_1^{-1}\bpsi(\bar x_k)} .
\]
The fraction is the SP coefficient \eqref{eq:sp-1d} with $W_1=M_1$, so that the condition \eqref{eq:mp-equals-sp} is satisfied and MP coincides with SP. In other words, each line $y=\eta_k$ is an
independent 1D ROD problem with its own distance $\bar x_k$.
\end{proof}

\section{Linear advection equation}
\label{sec:advection}

In this section \eqref{eq:system} is the linear advection equation with a source term,
\begin{equation}
  \partial_t u + \nabla\cdot(\ba u) = s(\bx) \quad \text{in } \bar\Omega,\qquad
  u = u_D \quad \text{on } \bar\Gamma^- = \{\bxb\in\bar\Gamma:\ \ba\cdot\bar\bn<0\},
  \label{eq:advection}
\end{equation}
with a constant speed $\ba$. The Dirichlet condition $u=u_D$ should be imposed on the inflow part $\bar\Gamma^-=\{\bx\in\bar\Gamma:\ \ba\cdot\bar\bn<0\}$ of the real boundary. In the embedded discretization, instead, the condition is imposed on the surrogate inflow boundary $\tilde\Gamma^-=\{\bxt\in\tilde\Gamma:\ \ba\cdot\tilde\bn<0\}$, made of the inflow faces of the computational boundary, where the data $u_D$ on $\bar\Gamma$ are not available: on these faces the upwind flux $\hat u(u_h^-,u_h^+,\tilde\bn)$ uses the ghost state $u_h^+=v_h(\bxt_k)$, given by one of the corrections of \cref{sec:methods}, which transfers the condition from $\bar\Gamma$ to $\tilde\Gamma$. The faces of the surrogate outflow boundary $\tilde\Gamma^+=\{\bxt\in\tilde\Gamma:\ \ba\cdot\tilde\bn>0\}$ need no boundary condition.

\subsection{Discrete operator of linear advection in a single cell}
\label{sec:single-cell}

We consider one Cartesian cell $K=[0,1]^2$ and the linear advection equation \eqref{eq:advection} with a constant unit speed
\begin{equation}
  \ba=(\cos\beta,\sin\beta)^T .
  \label{eq:beta}
\end{equation}
All the lengths of this section are in units of the mesh size, i.e.\ $\Delta x=1$. In the same way as for $\tilde\Gamma^\pm$, we denote by $\partial K^-$ and $\partial K^+$ the inflow ($\ba\cdot\tilde\bn<0$) and outflow ($\ba\cdot\tilde\bn>0$) faces of $K$, with $\tilde\bn$ the outward normal; the faces with $\ba\cdot\tilde\bn=0$ belong to neither set, since the upwind flux vanishes on them. An inflow face whose neighbour is cut by the real boundary lies on $\tilde\Gamma^-$ and carries one of the corrections of \cref{sec:methods}, with homogeneous data $u_D=0$. All the results use equispaced Lagrange bases and $p+1$ Gauss points per face.

Let $\bu\in\mathbb R^D$ be the coefficients of $u_h=\bphi^T\bu$ in the cell. We introduce the stiffness matrices $Q_x=\int_K\partial_x\bphi\,\bphi^T\,d\bx$ and $Q_y=\int_K\partial_y\bphi\,\bphi^T\,d\bx$ and, for each face $f$ with Gauss points $\bxt_k$ and weights $w_k$, $k=1,\dots,N_q=p+1$, the $D\times N_q$ matrices $\tilde\Phi_f=[\bphi(\bxt_1),\dots,\bphi(\bxt_{N_q})]$ and $\bar\Phi_f=[\bphi(\bxb_1),\dots,\bphi(\bxb_{N_q})]$ of the basis at the face points and at the mapped points, $\mathcal W_f={\rm diag}(w_k)$, and the face mass matrix $E_f=\tilde\Phi_f\mathcal W_f\tilde\Phi_f^T$. With the upwind flux, the weak form \eqref{eq:dg-weak} of the cell reads
\begin{equation}
  M\dot\bu=\bigl(a_xQ_x+a_yQ_y\bigr)\bu-\sum_{f\subset\partial K^+}(\ba\cdot\tilde\bn_f)\,E_f\bu-\sum_{f\subset\partial K^-}(\ba\cdot\tilde\bn_f)\,\tilde\Phi_f\mathcal W_f\,\mathbf v_f,
  \label{eq:adv-cell}
\end{equation}
where $\mathbf v_f=(v_h(\bxt_1),\dots,v_h(\bxt_{N_q}))^T$ collects the exterior values on an inflow face. On a face $f\subset\tilde\Gamma^-$, with $u_D=0$, the corrections \eqref{eq:sp} and \eqref{eq:mp} give $\mathbf v_f=\bigl(\tilde\Phi_f^T-\mathcal A_f\bar\Phi_f^T\bigr)\bu$, where $\mathcal A_f={\rm diag}(\alpha_k)$ for SP, with $\alpha_k=\aSPW(\bxt_k)$ and $\alpha_k=1$ for \SB, while for MP the $k$-th row of $\mathcal A_f$ is $[\aMPW(\bxt_k)]^T$. On an inflow face with an active neighbour, $\mathbf v_f$ is the trace of the neighbour, which does not depend on $\bu$. The semi-discrete scheme of the cell is then $\MM\dot\bu=\KK\bu$, with $\MM=M$ and
\begin{equation}
  \KK=a_xQ_x+a_yQ_y-\sum_{f\subset\partial K^+}(\ba\cdot\tilde\bn_f)\,E_f-\sum_{f\subset\partial K^-\cap\tilde\Gamma^-}(\ba\cdot\tilde\bn_f)\,\tilde\Phi_f\mathcal W_f\bigl(\tilde\Phi_f^T-\mathcal A_f\bar\Phi_f^T\bigr),
  \label{eq:adv-K}
\end{equation}
and we study the spectrum of the operator $\MM^{-1}\KK$. The correction enters only through the last term, a matrix of rank at most $N_q$ for each face on $\tilde\Gamma^-$.

\begin{remark}[One cell is enough]
As also discussed in~\cite{ciallella2027stability} for the 1D case, with a constant $\ba$ and the upwind flux, a cell receives information only from the neighbours behind its inflow faces, and the correction uses only the cell itself. Ordering the cells along the flow, the trace of a neighbour in \eqref{eq:adv-cell} only couples a cell with the cells that precede it, and the operator of the whole mesh is block lower triangular, with the matrices \eqref{eq:adv-K} of the single cells on the diagonal. Its spectrum is therefore the union of the spectra of the single cells, for any degree $p$, and the single-cell test decides the stability of the whole mesh, provided that it covers the geometry of the actual boundary cells.
\end{remark}

The fully discrete scheme uses the explicit time integration of order $p+1$ of \cref{sec:dg}. Applied to the linear system $\MM\dot\bu=\KK\bu$, all these methods have the amplification factor
\begin{equation}
  \bu^{n+1}=z(\mu)\,\bu^n,\qquad z(\mu)=\sum_{k=0}^{p+1}\frac{\mu^k}{k!},\qquad \mu=\lambda\Delta t,
  \label{eq:amplification}
\end{equation}
with $\lambda$ an eigenvalue of $\MM^{-1}\KK$, and we use the time step of the full solver,
\begin{equation}
  \Delta t=\frac{\CFLmax\,\text{CFL}}{|a_x|+|a_y|},\qquad \CFLmax=0.33,\,0.19,\,0.12,\,0.11\quad(p=1,\dots,4).
  \label{eq:dt}
\end{equation}
A configuration is stable at a given CFL when $|z(\mu)|\le1$ for all $\mu=\lambda\Delta t$, $\lambda\in\sigma(\MM^{-1}\KK)$. Two situations must be distinguished. If all the eigenvalues have $\operatorname{Re}\lambda\le0$, the semi-discrete scheme is stable, and the only question is the largest stable time step. If an eigenvalue has $\operatorname{Re}\lambda>0$, the scheme is unstable whatever the time step. Since the lengths are in units of the mesh size, on a mesh of size $\Delta x$ the eigenvalues are $\lambda/\Delta x$, and the growth rate increases when the mesh is refined.

\subsection{Constant distance: 2D = 1D}
\label{sec:constant}

\begin{figure}
\centering
\includegraphics[width=0.82\linewidth]{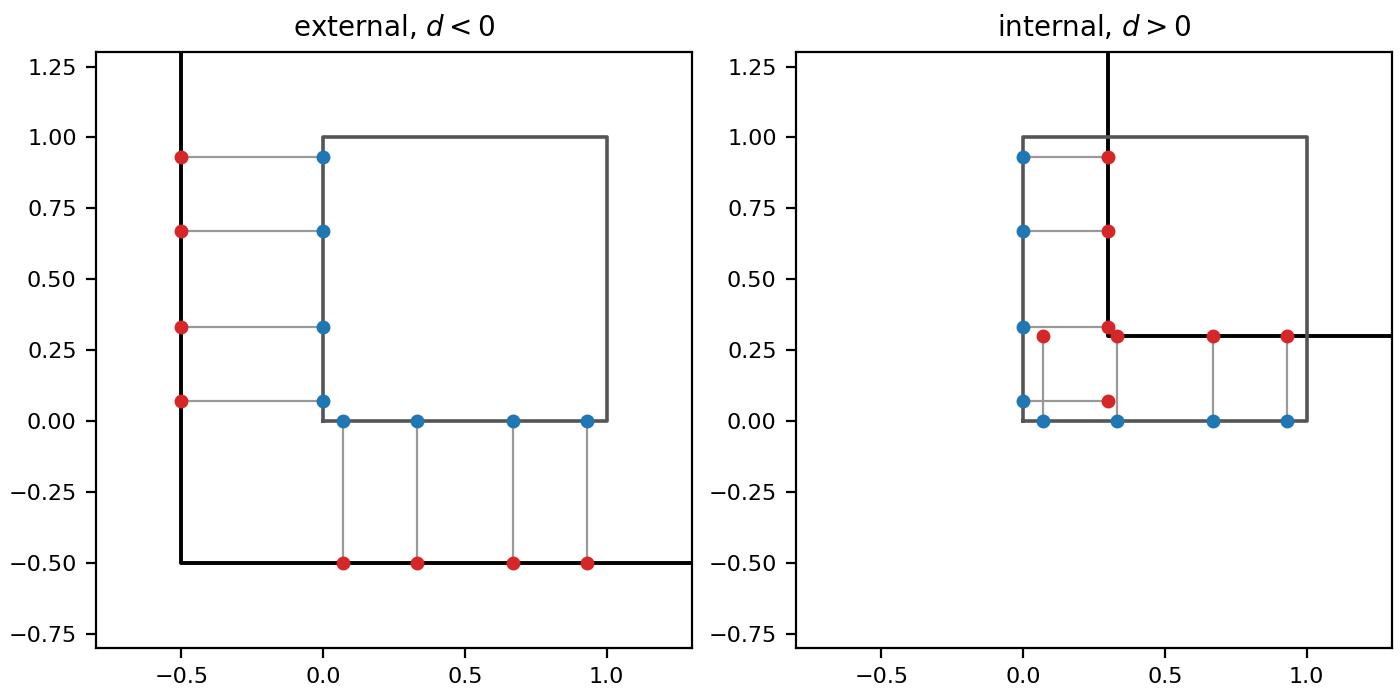}
\caption{Straight boundary at constant distance, $\PP{3}$: the cell, the real boundary (black), computational boundary points $\bxt_k$ (blue) and their images $\bxb_k$ (red), for an external (left) and an internal (right) boundary.}
\label{fig:const-config}
\end{figure}

We first recall the simplest configuration (\cref{fig:const-config}): a straight boundary parallel to the faces, for which each computational boundary point $(0,\eta_k)$ is
mapped to $(d,\eta_k)$, and $(\xi_k,0)$ to $(\xi_k,d)$, with $d\in[-1,1]$; $d<0$ is an external boundary (the real boundary lies outside the cell) and $d>0$ an internal one.

In this case the 2D stability problem reduces to the 1D one studied in~\cite{ciallella2027stability,ciallella2026minimization}.

\begin{proposition}[2D reduces to 1D at constant distance]\label{prop:constant-1d}
Consider one Cartesian cell with advection speed $\ba=[1,1]^T$, upwind outflow on the right and top faces, and a
boundary made of two straight segments, parallel to the left and bottom faces, which are computational. All the points of these faces are mapped at the same
distance $d$ from the face: $(0,\eta_k)$ to $(d,\eta_k)$ and $(\xi_k,0)$ to $(\xi_k,d)$. Then, for \SB, \RODE\ and \RODL, with either SP or MP, the 2D
single-cell operator is the Kronecker sum of the 1D single-cell operator $A_1$, i.e.\ the operator of the cell $[0,1]$ with
speed $1$ and boundary at the same distance $d$:
\begin{equation}
  \MM^{-1}\KK = A_1\otimes I + I\otimes A_1,\qquad
  \sigma(\MM^{-1}\KK)=\{\lambda_i+\lambda_l:\ \lambda_i,\lambda_l\in\sigma(A_1)\} .
  \label{eq:kronecker-sum}
\end{equation}
The factor $\tfrac12$ in the time step $\Delta t=\tfrac12\,\CFLmax\,\text{CFL}$ averages the two eigenvalues: the 2D values are
$\mu=\tfrac12(\mu_i+\mu_l)$, where $\mu_i=\text{CFL}\,\CFLmax\lambda_i$ are the 1D values.
\end{proposition}

\begin{proof}
By \cref{prop:mp-constant}, MP gives the same correction as SP, so we consider SP only. By \eqref{eq:sp-1d} with $\bar x_k=d$, the coefficient
$\alpha$ does not depend on the point and equals the 1D coefficient, and $\alpha=1$ for \SB. The modified polynomial at $\bxt_k=(0,\eta_k)$ then separates:
\[
  \bphi(\bxt_k)-\alpha\,\bphi(\bxb_k)=\big(\bpsi(0)-\alpha\bpsi(d)\big)\otimes\bpsi(\eta_k).
\]
Since $M_1=\sum_kw_k\bpsi(\eta_k)\bpsi(\eta_k)^T$, the sum over the Gauss points of the left face gives
\[
  \sum_kw_k\,\bphi(\bxt_k)\big(\bphi(\bxt_k)-\alpha\,\bphi(\bxb_k)\big)^T
  =\bpsi(0)\big(\bpsi(0)-\alpha\bpsi(d)\big)^T\otimes M_1 .
\]
In the same way, the outflow term of the right face is $-\bpsi(1)\bpsi(1)^T\otimes M_1$. In the volume term with the derivative in $x$, the derivative acts
only on the first factor, $\partial_x\bphi=\bpsi'(x)\otimes\bpsi(y)$, so that the integral splits as
\[
  \int_0^1\!\!\int_0^1\partial_x\bphi\,\bphi^T\,dx\,dy
  =\Big(\int_0^1\bpsi'\bpsi^T\,dx\Big)\otimes\Big(\int_0^1\bpsi\bpsi^T\,dy\Big)
  =S_1\otimes M_1,\qquad S_1=\int_0^1\bpsi'\bpsi^T\,dx .
\]
All the terms related to the direction $x$ thus sum to $B_1\otimes M_1$, where
\[
  B_1=S_1-\bpsi(1)\bpsi(1)^T+\bpsi(0)\big(\bpsi(0)-\alpha\bpsi(d)\big)^T ,
\]
and, by symmetry, the terms related to the direction $y$ sum to $M_1\otimes B_1$. Since $\MM=M_1\otimes M_1$, we get \eqref{eq:kronecker-sum} with $A_1=M_1^{-1}B_1$, which is the 1D single-cell operator. The values of $\mu$ follow by multiplying the eigenvalues by $\Delta t$.
\end{proof}

\begin{figure}
\centering
\includegraphics[width=0.9\linewidth]{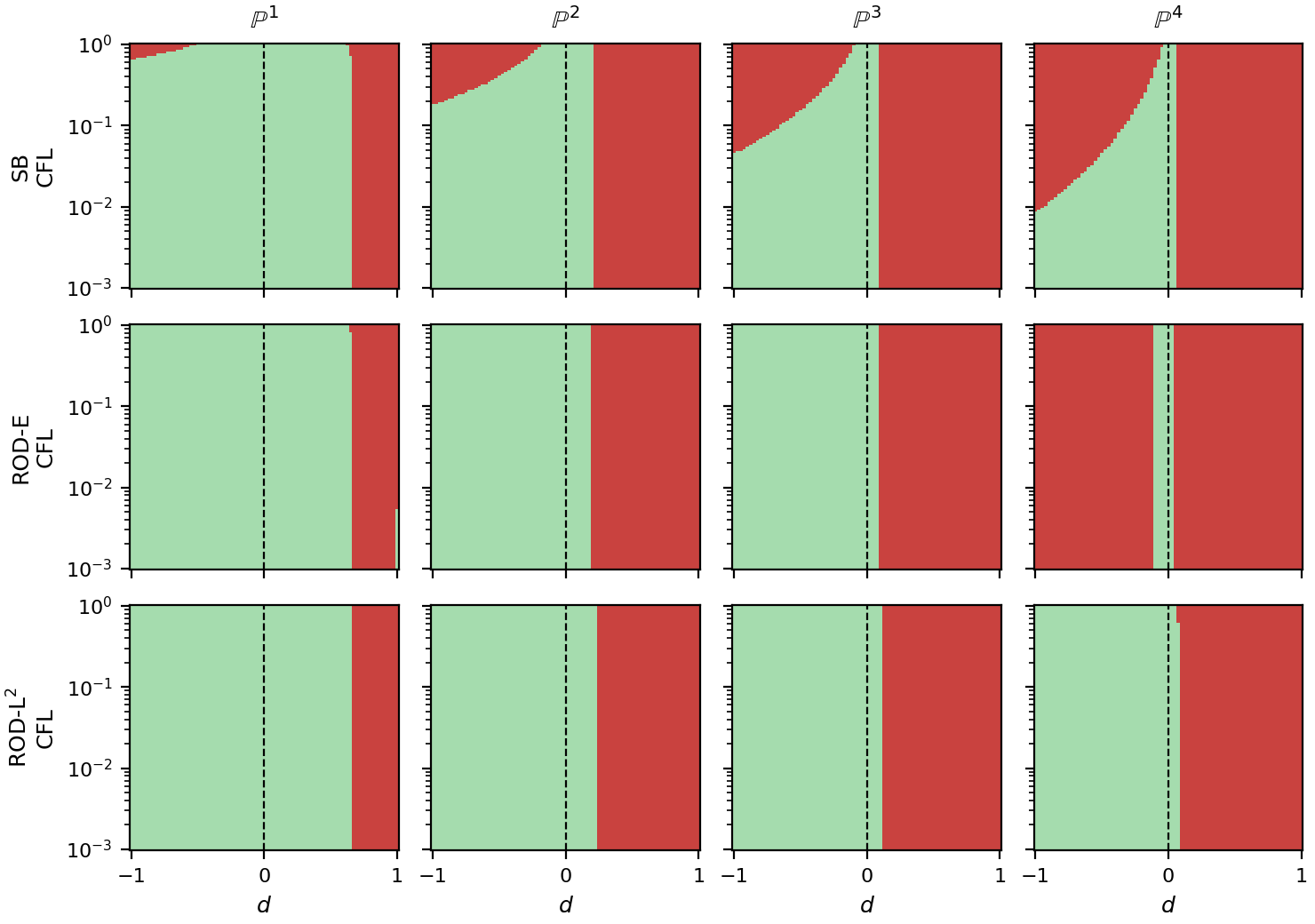}
\caption{Constant distance: stability regions in the plane $(d,\text{CFL})$, green stable, red unstable. Dashed line: $d=0$, with external boundaries on its left.}
\label{fig:const-maps}
\end{figure}

The stability regions of the corrections in the plane $(d,\text{CFL})$ are shown in \cref{fig:const-maps}. For external boundaries, \RODL\ is stable up to $\PP{4}$ and \RODE\ up to $\PP{3}$ at the nominal time step, and their stability is independent of the CFL, i.e.\ it is decided by the sign of the real part of the eigenvalues. The \SB\ correction is stable for every external distance up to $\PP{4}$, but it is limited by the time step: the largest stable CFL for all the external distances is $0.63$, $0.18$, $0.045$, $0.008$ for $\PP{1},\dots,\PP{4}$.

\begin{remark}[Internal configurations]\label{rem:internal}
For $d>0$ the real boundary lies inside the cell and the computational face is upstream of it. The value at $\bar\Gamma$ is then already determined by what entered through the face at an earlier time, and satisfying the condition would require transporting it against the characteristic, which a scheme that advances in time cannot do. Consistently, \cref{fig:const-maps} shows that internal configurations are stable only very close to the face, within a distance that decreases with the degree. Since they never occur when the cut cells are discarded, internal configurations are not considered in the rest of the paper.
\end{remark}

\subsection{A tilted straight boundary}
\label{sec:tilted}

\begin{figure}
\centering
\includegraphics[width=0.36\linewidth]{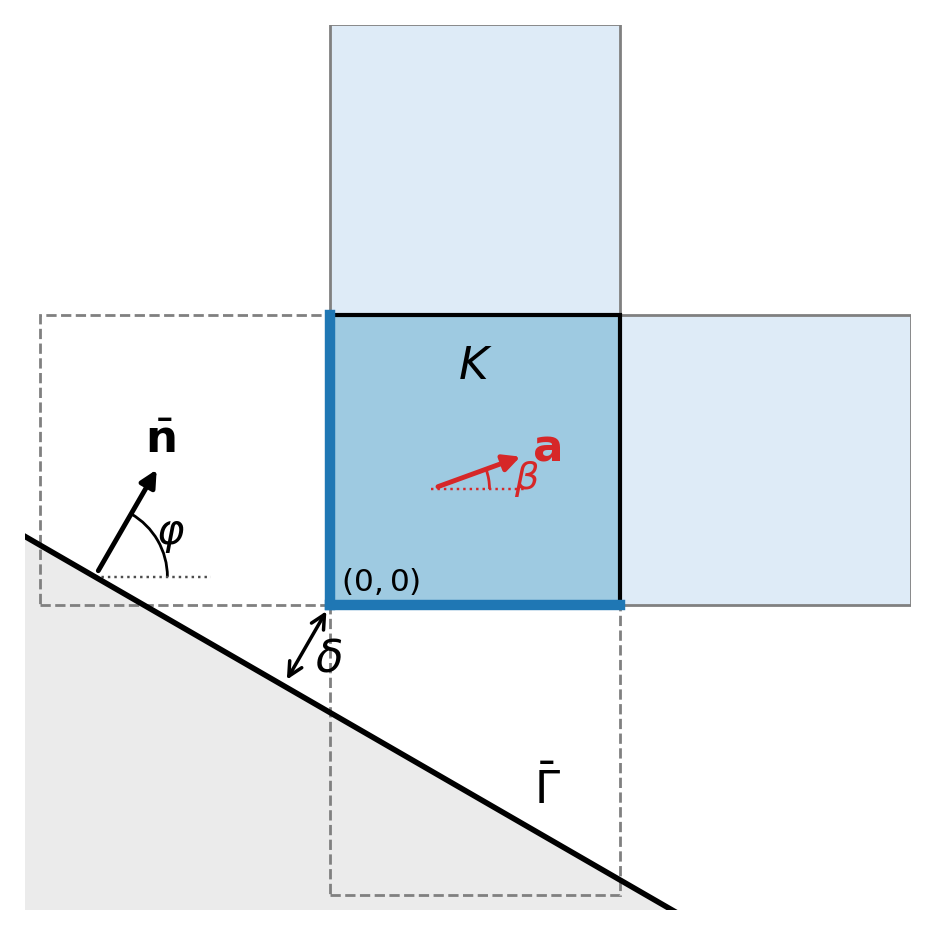}
\caption{Tilted straight boundary: the cell $K$, its active right and top neighbours (light) and its left and bottom neighbours cut by the real boundary $\bar\Gamma$ (dashed), the computational faces (thick blue), the angle $\varphi$ of the normal $\bar\bn$, the distance $\delta$ of the corner $(0,0)$ from $\bar\Gamma$, and the direction $\beta$ of the flow $\ba$.}
\label{fig:tilted-setup}
\end{figure}

At constant distance the stability of the 2D cell is that of the 1D cell (\cref{prop:constant-1d}). We now show that this is lost as soon as the straight boundary is tilted with respect to the mesh. The real boundary is the line
\begin{equation}
  \bar\Gamma=\{\bx:\ \bar\bn\cdot\bx=-\delta\},\qquad \bar\bn=(\cos\varphi,\sin\varphi)^T,
  \label{eq:line}
\end{equation}
where $\varphi$ is the angle of the normal $\bar\bn$, pointing into the fluid, and $\delta>0$ the distance of the corner $(0,0)$ of the cell from the line (\cref{fig:tilted-setup}). The values $\varphi=0^\circ$ and $\varphi=90^\circ$ are the constant-distance configurations of \cref{sec:constant}. We take the configuration

\begin{equation}
  \varphi=85^\circ,\qquad \delta=0.025,\qquad \beta=10^\circ:
  \label{eq:tilted-config}
\end{equation}
the line makes an angle of $5^\circ$ with the bottom face and passes $0.025$ below the corner, both neighbours are cut, and both the left and the bottom faces are inflow faces, with $\ba\cdot\tilde\bn=-\cos10^\circ$ and $-\sin10^\circ$.

\begin{figure}
\centering
\includegraphics[width=0.92\linewidth]{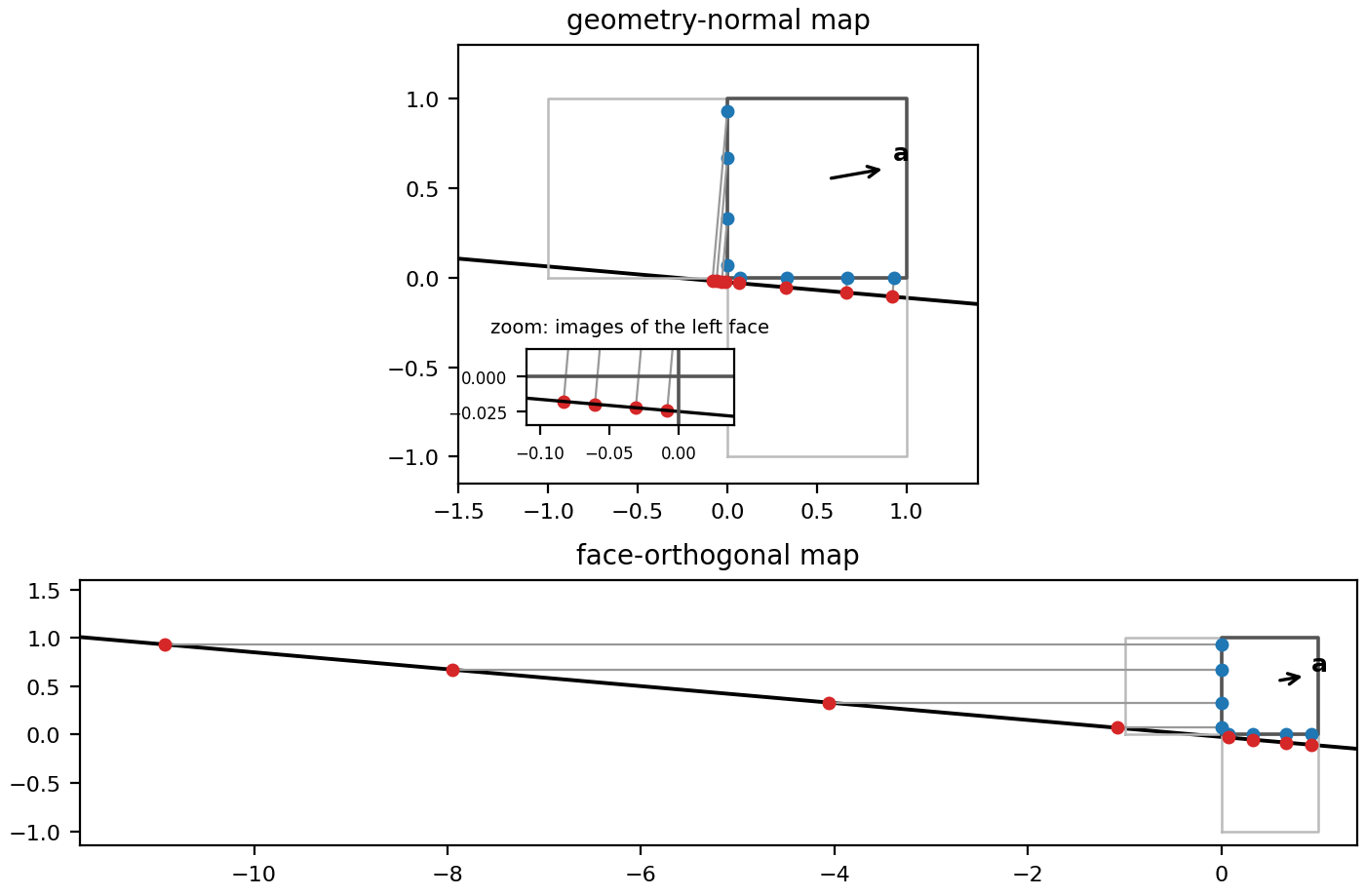}
\caption{The tilted straight boundary \eqref{eq:tilted-config}, $\PP{3}$: the cell (dark), its left and bottom neighbours (light), the real boundary (black), the computational boundary points $\bxt_k$ (blue), their images $\bxb_k$ (red) and the direction $\ba$ of the flow. Top: geometry-normal map, with a zoom on the images of the left face. Bottom: face-orthogonal map; the scale of the horizontal axis is different.}
\label{fig:mapping}
\end{figure}

\Cref{fig:mapping} shows the mapped points for $\PP{3}$. On the bottom face, which is almost parallel to $\bar\Gamma$, the two maps are close to each other and to the constant-distance case. On the left face, which is almost orthogonal to $\bar\Gamma$, they behave in opposite ways. The geometry-normal map projects the four points of the face close to each other, so that the correction of the whole face is built from the values of $u_h$ in this small region. The face-orthogonal map, instead, moves each point along the horizontal line through it and keeps its height $\eta_k$; since this line is almost parallel to $\bar\Gamma$, the distances are much larger, and go up to $\sim11$ cells.

\begin{table}
\centering
\caption{Tilted straight boundary \eqref{eq:tilted-config}: largest real part of the eigenvalues of the single-cell operator. Bold: positive, i.e.\ unstable for any time step.}
\label{tab:tilted}
\small
\begin{tabular}{lcccc@{\qquad}cccc}
\toprule
 & \multicolumn{4}{c}{geometry-normal map} & \multicolumn{4}{c}{face-orthogonal map}\\
\cmidrule(lr){2-5}\cmidrule(lr){6-9}
 & $\PP{1}$ & $\PP{2}$ & $\PP{3}$ & $\PP{4}$ & $\PP{1}$ & $\PP{2}$ & $\PP{3}$ & $\PP{4}$\\
\midrule
\SB     & $\mathbf{0.19}$ & $\mathbf{1.29}$ & $\mathbf{3.11}$ & $\mathbf{5.61}$ & $-0.59$ & $-0.45$ & $-0.30$ & $\mathbf{0.042}$\\
\RODESP & $\mathbf{0.025}$ & $\mathbf{0.52}$ & $\mathbf{1.70}$ & $\mathbf{3.31}$ & $-0.57$ & $-0.55$ & $-0.34$ & $\mathbf{0.38}$\\
\RODLSP & $-0.049$ & $\mathbf{0.55}$ & $\mathbf{1.85}$ & $\mathbf{3.78}$ & $-0.56$ & $-0.55$ & $-0.33$ & $-0.041$\\
\RODEMP & $-0.23$ & $\mathbf{0.65}$ & $\mathbf{1.01}$ & $\mathbf{1.86}$ & $-0.57$ & $-0.53$ & $\mathbf{0.45}$ & $\mathbf{4.97}$\\
\RODLMP & $-0.30$ & $\mathbf{0.88}$ & $\mathbf{0.91}$ & $\mathbf{2.56}$ & $-0.56$ & $-0.55$ & $-0.33$ & $-0.041$\\
\bottomrule
\end{tabular}
\end{table}

With the geometry-normal map (\cref{tab:tilted}, left) all the corrections have eigenvalues with positive real part from $\PP{2}$ on, and \SB\ and \RODESP{} already with $\PP{1}$. A tilt of $5^\circ$ is therefore enough to make every correction unstable, and this cannot be fixed by the time step or by the time integrator. The direction of the flow matters: with $\ba=(1,1)^T/\sqrt2$, i.e.\ $\beta=45^\circ$, on the same line, \SB\ is unstable from $\PP{2}$, \RODESP{} and \RODLSP{} from $\PP{3}$, and \RODEMP{} and \RODLMP{} only with $\PP{4}$. The diagonal direction is favourable, and a stability test made only with $\ba=(1,1)^T$ would miss most of these instabilities.

\begin{remark}[MP with the geometry-normal map: loss of rank]\label{rem:mp-degenerate}
When the mapped points of a face collapse, as on the left face of \cref{fig:mapping}, the matrix $\Phib$ of the MP correction is close to rank deficient. This is not a limitation in practice: on such faces the MP correction can simply be replaced by the SP correction \eqref{eq:sp}, which only needs $\bphi(\bxb_k)$ at each point and is always defined. The loss of rank is thus a question of implementation, and not the cause of the instabilities of \cref{tab:tilted}, which also affect the SP corrections.
\end{remark}

In this context, the face-orthogonal map is much more robust. On the same configuration (\cref{tab:tilted}, right) it is stable for all the corrections up to $\PP{2}$, and up to $\PP{3}$ for all but \RODEMP{}; with $\PP{4}$ only the \RODL\ corrections are stable. This happens although the left face is mapped up to $11$ cells away, i.e.\ although the polynomial of the cell is extrapolated much farther than with the geometry-normal map. The extrapolation distance is thus not what decides the stability. A likely explanation is that, with a tensor-product basis, the face-orthogonal map keeps the coordinate of each point along the face, so that the correction at a point only involves the polynomial along the line of the cell through it, and the boundary condition keeps a one-dimensional character, close to the 1D problems whose stability is known. The geometry-normal map, instead, also moves the points along the face, and the correction of one line uses the polynomial of other lines, as on the left face of \cref{fig:mapping}, where all the points are mapped close to the corner.

\subsection{Numerical experiments: manufactured solution on a disk}

The single-cell model says whether a correction can be stable, but not how accurate it is, and it contains a single cell. All the corrections have therefore been tested in the 2D DG solver on a curved boundary, with both maps, to see where they are stable and where they converge. We solve the linear advection equation with $\ba=(1,1)^T$ and the stationary manufactured solution $u=1+\sin(x+y)$, with the corresponding source term, in the unit disk embedded in $[-1.2,1.2]^2$. The cells cut by the circle are discarded, so that all the configurations are external, and we use the nested meshes with $8\times8$ to $64\times64$ cells, the explicit time integration of order $p+1$ and no stabilization. As in the 1D study, no final time is prescribed: the computation is advanced until the time residual $\|\bu^{n+1}-\bu^n\|_2$ reaches $10^{-13}$, i.e.\ until the discrete stationary state is reached, with a cap at $t=500$, and a computation that diverges or never reaches the stationary state is considered unstable. The error is measured on the converged stationary state. 

\begin{figure}
\centering
\includegraphics[width=0.95\linewidth]{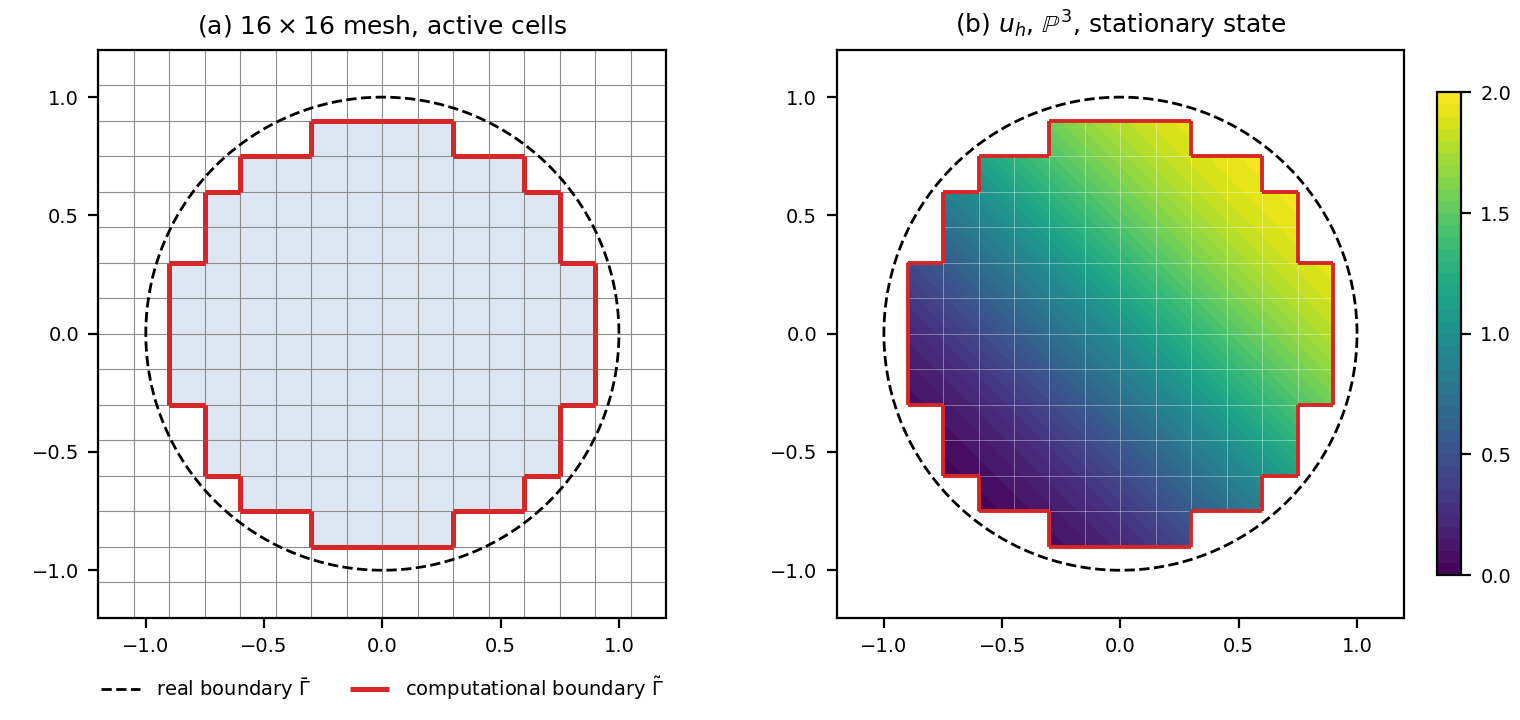}
\caption{Linear advection in the unit disk, $\ba=(1,1)^T$, $16\times16$ cells. (a) Active cells (shaded), computational boundary $\tilde\Gamma$ (faces between active and discarded cells) and real boundary $\bar\Gamma$. (b) Stationary numerical solution $u_h$ with \RODLSP{} and the face-orthogonal map, $\PP{3}$; $L^2$ error $3.98\cdot10^{-5}$ (\cref{tab:acc-faceorth}).}
\label{fig:disk}
\end{figure}

\Cref{fig:disk} shows the $16\times16$ mesh and a stationary solution. Since the cut cells are discarded, the computational boundary $\tilde\Gamma$ is a staircase inscribed in the circle, and its faces lie at distances from the circle that range from a small fraction of the cell to more than one cell, with all the orientations of the boundary relative to the faces. Each boundary cell is thus in a different configuration of the single-cell model, and a correction converges on a given mesh only if all of them are stable. With $\ba=(1,1)^T$, the condition is imposed on the inflow faces of the staircase, i.e.\ the faces with outward normal $(-1,0)$, on the left half of the circle, and $(0,-1)$, on the lower half, and the correction transfers the data from the circle to these faces. The solution computed with \RODLSP{}, face-orthogonal map and $\PP{3}$ is shown in \cref{fig:disk}(b).

\begin{table}
\centering
\caption{Linear advection equation in 2D, unit disk, geometry-normal map: $L^2$ errors $\|u-u_{ex}\|$
and estimated order of accuracy (EOA) with the \SB\ and ROD polynomial corrections and explicit time
integration for polynomials $\PP{p}$.}
\label{tab:acc-normal}
\small
\setlength{\tabcolsep}{5pt}
\begin{tabular}{lcccccc}
\toprule
\SB & \spancell{$\PP{1}$, CFL=0.2} & \spancell{$\PP{2}$, CFL=0.03} & \spancell{$\PP{3}$, CFL=0.005}\\
$N_e$ & $L^2$ error & EOA & $L^2$ error & EOA & $L^2$ error & EOA\\
\midrule
$8\times8$   & 1.07E-01 & --   & \spanrows{unstable} & \spanrows{unstable}\\
$16\times16$ & 3.51E-02 & 1.61 & \spanskip & \spanskip\\
$32\times32$ & 1.18E-02 & 1.57 & \spanskip & \spanskip\\
$64\times64$ & 2.38E-03 & 2.31 & \spanskip & \spanskip\\
\midrule
\RODESP & \spancell{$\PP{1}$, CFL=1.0} & \spancell{$\PP{2}$, CFL=1.0} & \spancell{$\PP{3}$, CFL=1.0}\\
$N_e$ & $L^2$ error & EOA & $L^2$ error & EOA & $L^2$ error & EOA\\
\midrule
$8\times8$   & 3.16E-02 & --   & \spanrows{unstable} & \spanrows{unstable}\\
$16\times16$ & 9.23E-03 & 1.78 & \spanskip & \spanskip\\
$32\times32$ & 1.89E-03 & 2.28 & \spanskip & \spanskip\\
$64\times64$ & 5.64E-04 & 1.75 & \spanskip & \spanskip\\
\midrule
\RODLSP & \spancell{$\PP{1}$, CFL=1.0} & \spancell{$\PP{2}$, CFL=1.0} & \spancell{$\PP{3}$, CFL=1.0}\\
$N_e$ & $L^2$ error & EOA & $L^2$ error & EOA & $L^2$ error & EOA\\
\midrule
$8\times8$   & 3.17E-02 & --   & \spanrows{unstable} & \spanrows{unstable}\\
$16\times16$ & 9.24E-03 & 1.78 & \spanskip & \spanskip\\
$32\times32$ & \spancell{unstable} & \spanskip & \spanskip\\
$64\times64$ & \spancell{unstable} & \spanskip & \spanskip\\
\midrule
\RODEMP & \spancell{$\PP{1}$, CFL=1.0} & \spancell{$\PP{2}$, CFL=1.0} & \spancell{$\PP{3}$, CFL=1.0}\\
$N_e$ & $L^2$ error & EOA & $L^2$ error & EOA & $L^2$ error & EOA\\
\midrule
$8\times8$   & 5.37E-02 & --   & 4.53E-03 & --   & \spanrows{unstable}\\
$16\times16$ & 1.52E-02 & 1.82 & 3.58E-04 & 3.66 & \spanskip\\
$32\times32$ & 5.06E-03 & 1.58 & 8.75E-05 & 2.03 & \spanskip\\
$64\times64$ & 1.11E-03 & 2.20 & 6.21E-06 & 3.82 & \spanskip\\
\midrule
\RODLMP & \spancell{$\PP{1}$, CFL=1.0} & \spancell{$\PP{2}$, CFL=1.0} & \spancell{$\PP{3}$, CFL=1.0}\\
$N_e$ & $L^2$ error & EOA & $L^2$ error & EOA & $L^2$ error & EOA\\
\midrule
$8\times8$   & 4.55E-02 & --   & 3.16E-03 & --   & 1.91E-04 & --\\
$16\times16$ & 1.40E-02 & 1.70 & 5.14E-04 & 2.62 & 2.67E-05 & 2.84\\
$32\times32$ & 4.45E-03 & 1.65 & 8.83E-05 & 2.54 & 2.47E-06 & 3.44\\
$64\times64$ & 9.61E-04 & 2.21 & 8.30E-06 & 3.41 & 1.25E-07 & 4.31\\
\bottomrule
\end{tabular}

\smallskip
{\footnotesize Unstable: the computation diverges.}
\end{table}

\begin{table}
\centering
\caption{Linear advection equation in 2D, unit disk, face-orthogonal map: $L^2$ errors $\|u-u_{ex}\|$
and estimated order of accuracy (EOA) with the \SB\ and ROD polynomial corrections and explicit time
integration for polynomials $\PP{p}$. \RODLSP{} and \RODLMP{} give identical results
(\cref{prop:mp-faceorth}).}
\label{tab:acc-faceorth}
\small
\setlength{\tabcolsep}{5pt}
\begin{tabular}{lcccccc}
\toprule
\SB & \spancell{$\PP{1}$, CFL=0.3} & \spancell{$\PP{2}$, CFL=0.03} & \spancell{$\PP{3}$, CFL=0.0025}\\
$N_e$ & $L^2$ error & EOA & $L^2$ error & EOA & $L^2$ error & EOA\\
\midrule
$8\times8$   & 9.25E-02 & --   & 4.49E-03 & --   & 5.02E-04 & --\\
$16\times16$ & 3.19E-02 & 1.54 & 7.96E-04 & 2.49 & 5.84E-05 & 3.10\\
$32\times32$ & 1.08E-02 & 1.56 & 1.13E-04 & 2.81 & \spancell{not run$^\ast$}\\
$64\times64$ & 2.30E-03 & 2.23 & \spancell{not run$^\ast$} & \spancell{not run$^\ast$}\\
\midrule
\RODESP & \spancell{$\PP{1}$, CFL=1.0} & \spancell{$\PP{2}$, CFL=1.0} & \spancell{$\PP{3}$, CFL=1.0}\\
$N_e$ & $L^2$ error & EOA & $L^2$ error & EOA & $L^2$ error & EOA\\
\midrule
$8\times8$   & 6.84E-02 & --   & 2.01E-03 & --   & \spanrows{unstable}\\
$16\times16$ & 2.25E-02 & 1.61 & 4.01E-04 & 2.32 & \spanskip\\
$32\times32$ & 7.46E-03 & 1.59 & 5.87E-05 & 2.77 & \spanskip\\
$64\times64$ & 1.56E-03 & 2.26 & 6.37E-06 & 3.20 & \spanskip\\
\midrule
\RODEMP & \spancell{$\PP{1}$, CFL=1.0} & \spancell{$\PP{2}$, CFL=1.0} & \spancell{$\PP{3}$, CFL=1.0}\\
$N_e$ & $L^2$ error & EOA & $L^2$ error & EOA & $L^2$ error & EOA\\
\midrule
$8\times8$   & 6.85E-02 & --   & 2.01E-03 & --   & \spanrows{unstable}\\
$16\times16$ & 2.25E-02 & 1.61 & 4.01E-04 & 2.32 & \spanskip\\
$32\times32$ & 7.46E-03 & 1.59 & 5.86E-05 & 2.77 & \spanskip\\
$64\times64$ & 1.56E-03 & 2.26 & 6.36E-06 & 3.20 & \spanskip\\
\midrule
\RODLSP{} = MP & \spancell{$\PP{1}$, CFL=1.0} & \spancell{$\PP{2}$, CFL=1.0} & \spancell{$\PP{3}$, CFL=1.0}\\
$N_e$ & $L^2$ error & EOA & $L^2$ error & EOA & $L^2$ error & EOA\\
\midrule
$8\times8$   & 6.36E-02 & --   & 2.80E-03 & --   & 3.39E-04 & --\\
$16\times16$ & 2.07E-02 & 1.62 & 5.43E-04 & 2.37 & 3.98E-05 & 3.09\\
$32\times32$ & 6.90E-03 & 1.58 & 7.88E-05 & 2.79 & 3.31E-06 & 3.59\\
$64\times64$ & 1.43E-03 & 2.28 & 8.55E-06 & 3.20 & 1.64E-07 & 4.34\\
\bottomrule
\end{tabular}

\smallskip
{\footnotesize Unstable: the computation diverges. $^\ast$ Too expensive at this CFL.}
\end{table}

The results are given in \cref{tab:acc-normal,tab:acc-faceorth}. \RODLSP{} $\PP{4}$ with the geometry-normal map, not reported in the tables, diverges. All the configurations that are stable converge, with an order on the finest meshes close to $p+1$. On the coarse meshes the order is irregular, as expected when the number of boundary cells and their configurations change with the mesh.

With the face-orthogonal map, \RODE\ is stable up to $\PP{2}$ and \RODL\ up to $\PP{3}$, and MP gives the same errors as SP. In 1D, and in 2D at constant distance, \RODE\ is stable up to $\PP{3}$ and \RODL\ up to $\PP{4}$ for all the external distances (\cref{sec:constant}). The face-orthogonal map therefore keeps the stability of the 1D problems up to one degree less. As discussed in \cref{sec:tilted}, this map likely gives a one-dimensional character to the boundary condition, and the degree that is lost can be attributed to the coupling, inside the cell, of lines with different distances. \SB\ with this map has no eigenvalue with positive real part up to $\PP{3}$, but it is limited by a time step that must be lowered: the tables use $0.3$, $0.03$ and $0.0025$. On the faces almost orthogonal to the boundary the face-orthogonal map gives large distances (\cref{fig:mapping}), over which \SB\ extrapolates the polynomial, hence reducing the time step.

With the geometry-normal map only MP is stable beyond $\PP{1}$, and \RODLMP{} is the only correction stable up to $\PP{3}$. \RODLSP{} $\PP{1}$ converges on the two coarse meshes and diverges on the two fine ones: the boundary cells of the coarse meshes do not reach an unstable configuration, while refining the mesh adds boundary cells in new positions, which eventually do. These results, however, hold only for $\ba=(1,1)^T$, which is the favourable direction for this map (\cref{sec:tilted}). To check this in the full solver, the combinations of the geometry-normal map that converge with $\ba=(1,1)^T$ have been run again with $\beta=10^\circ$ in \eqref{eq:beta}, on the same disk, with the same manufactured solution and the same CFL (\cref{tab:acc-beta}). All the MP corrections that converge with $\ba=(1,1)^T$ diverge with $\beta=10^\circ$, on both meshes. With $\PP{1}$, \SB\ and \RODESP{} still reach the stationary state, but the error of \RODESP{} is much larger and converges with order close to $1$ (\cref{tab:acc-beta}).

\begin{table}
\centering
\caption{Unit disk, geometry-normal map: the combinations that converge with $\ba=(1,1)^T$ (\cref{tab:acc-normal}), run with $\beta=10^\circ$. $L^2$ errors of the stationary state, or a label that the computation has diverged.}
\label{tab:acc-beta}
\small
\begin{tabular}{llcccc}
\toprule
 & & \multicolumn{2}{c}{$\ba=(1,1)^T$} & \multicolumn{2}{c}{$\beta=10^\circ$}\\
\cmidrule(lr){3-4}\cmidrule(lr){5-6}
correction & degree, CFL & $16\times16$ & $32\times32$ & $16\times16$ & $32\times32$\\
\midrule
\SB     & $\PP{1}$, 0.2 & 3.51E-02 & 1.18E-02 & 3.74E-02 & 1.17E-02 \\
\RODESP & $\PP{1}$, 1.0 & 9.23E-03 & 1.89E-03 & 2.57E-01 & 1.17E-01 \\
\RODEMP & $\PP{2}$, 1.0 & 3.58E-04 & 8.75E-05 & diverges & diverges \\
\RODLMP & $\PP{2}$, 1.0 & 5.14E-04 & 8.83E-05 & diverges & diverges \\
\RODLMP & $\PP{3}$, 1.0 & 2.67E-05 & 2.47E-06 & diverges & diverges \\
\bottomrule
\end{tabular}
\end{table}

In summary, the extension of the polynomial corrections from 1D to 2D is not trivial, even for linear advection on Cartesian meshes with the cut cells discarded. With the geometry-normal map no correction is stable for all the configurations beyond $\PP{1}$, since the stability depends both on the position of the boundary and on the direction of the flow, and \RODLMP{} reaches $\PP{3}$ only for the favourable diagonal flow. With the face-orthogonal map, which likely gives a one-dimensional character to the boundary condition, \RODE\ is stable up to $\PP{2}$ and \RODL\ up to $\PP{3}$ on the disk, one degree less than in 1D, while \SB\ is stable only with a small enough time step that adapts to the mesh size.

\section{Linear acoustics and the slip wall}
\label{sec:acoustics}

We now move to a hyperbolic system, the 2D linear acoustic system with sound speed $1$,
\begin{equation}
  \partial_tp+\nabla\cdot\bu=0,\qquad \partial_t\bu+\nabla p=0,
  \label{eq:acoustics}
\end{equation}
which is \eqref{eq:system} with $\mathbf U=(p,u^x,u^y)^T$, $\mathbf F(\mathbf U)=(u^x,p,0)^T$, $\mathbf G(\mathbf U)=(u^y,0,p)^T$ and $\mathbf S=\mathbf 0$, discretized as in \cref{sec:dg}. With a face of outward normal $\tilde\bn$, the normal flux is $A_{\tilde\bn}\mathbf U$, $s_{\max}=1$, and the Rusanov flux \eqref{eq:llf} reads
\begin{equation}
  \hat{\mathbf F}(\mathbf U^-,\mathbf U^+,\tilde\bn)=\tfrac12A_{\tilde\bn}\bigl(\mathbf U^-+\mathbf U^+\bigr)+\tfrac12\bigl(\mathbf U^--\mathbf U^+\bigr),\qquad
  A_{\tilde\bn}=\begin{bmatrix}0&\tilde n_x&\tilde n_y\\\tilde n_x&0&0\\\tilde n_y&0&0\end{bmatrix}.
  \label{eq:rusanov}
\end{equation}
Along $\tilde\bn$ the system has an outgoing wave $p+\bu\cdot\tilde\bn$ with speed $+1$, an incoming wave $p-\bu\cdot\tilde\bn$ with speed $-1$, and the tangential velocity, which is not transported. At every boundary one wave enters the domain and exactly one condition is imposed. This is the main difference with linear advection, where the outflow faces carry no condition at all, and it is the reason why the operator of a mesh is no longer block triangular: every cell exchanges information with all its neighbours, and the stability of a boundary correction can only be decided on a set of cells.

On a straight boundary at constant distance each wave then behaves as linear advection, and one cell reproduces the stability regions of the 1D inlet~\cite{ciallella2026minimization}. We do not repeat this analysis, and we focus on the slip wall, which has no counterpart for a scalar equation.

\subsection{The slip wall condition}
\label{sec:wall-condition}

At a wall the condition is $\bu\cdot\bar\bn=0$, with $\bar\bn$ the normal of the wall. It is carried by the incoming wave, while the pressure and the tangential velocity are not prescribed, so that the polynomial correction acts on the normal velocity only. Let $\bar\bn_k=(\bar n_{x,k},\bar n_{y,k})^T$ be the normal of the wall at the mapped point $\bxb_k$ of the face point $\bxt_k$, and $\btaub_k=(\bar\tau_{x,k},\bar\tau_{y,k})^T=(-\bar n_{y,k},\bar n_{x,k})^T$ the tangent. The correction is built from the normal velocity of the cell polynomial evaluated on the real boundary.
In the following we consider first the single-point (SP) correction, which uses only the mapped point of the face point, and then we extend it to the multi-point (MP) correction, which uses all the mapped points of the face.

For SP, the corrected velocity at the face point $\bxt_k$ is obtained from the vector polynomial $\bv_h=(\bphi^T\bv^x,\bphi^T\bv^y)$ closest to $\bu_h=(\bphi^T\bu^x,\bphi^T\bu^y)$ in the norm of the weight $W$ applied to both components, under the single condition $\bar\bn_k\cdot\bv_h(\bxb_k)=0$. This is the stationary point of
\begin{equation}
  \mathcal L(\bv^x,\bv^y,\lambda)=\tfrac12(\bv^x-\bu^x)^TW(\bv^x-\bu^x)+\tfrac12(\bv^y-\bu^y)^TW(\bv^y-\bu^y)+\lambda\,\bar\bn_k\cdot\bv_h(\bxb_k),
\end{equation}
the vector counterpart of the scalar functional of \cref{sec:methods}. With $\phit_k=\bphi(\bxt_k)$ and $\bar\bphi_k=\bphi(\bxb_k)$, and since $\bar\bn_k\cdot\bv_h(\bxb_k)=\bar n_{x,k}\,\bar\bphi_k^T\bv^x+\bar n_{y,k}\,\bar\bphi_k^T\bv^y$, the derivatives with respect to $\bv^x$ and $\bv^y$ give
\begin{equation*}
  \bv^x=\bu^x-\lambda\,\bar n_{x,k}\,W^{-1}\bar\bphi_k,\qquad
  \bv^y=\bu^y-\lambda\,\bar n_{y,k}\,W^{-1}\bar\bphi_k,
\end{equation*}
i.e.\ both components are modified by the same function $\bphi^TW^{-1}\bar\bphi_k$, and the modification $\bu_h-\bv_h$ is everywhere parallel to $\bar\bn_k$. Inserting this in the constraint, with $|\bar\bn_k|=1$, gives
\begin{equation*}
  \bar\bn_k\cdot\bu_h(\bxb_k)-\lambda\,\bar\bphi_k^TW^{-1}\bar\bphi_k=0
  \quad\Longrightarrow\quad
  \lambda=\frac{\bar\bn_k\cdot\bu_h(\bxb_k)}{\bar\bphi_k^TW^{-1}\bar\bphi_k},
\end{equation*}
and evaluating $\bv_h$ at the face point,
\begin{equation*}
  \bv_h(\bxt_k)=\bu_h(\bxt_k)-\lambda\,\phit_k^TW^{-1}\bar\bphi_k\,\bar\bn_k
  =\bu_h(\bxt_k)-\alpha_k\,\bigl(\bar\bn_k\cdot\bu_h(\bxb_k)\bigr)\,\bar\bn_k,
\end{equation*}
with $\alpha_k=\aSPW(\bxt_k)$ the coefficient of the scalar case \eqref{eq:sp-alpha}. The SP correction therefore applies the scalar correction to the normal velocity $\bar\bn_k\cdot\bu_h$, and leaves unchanged the velocity along the tangent of the wall at $\bxb_k$. The \SB\ correction is again retrieved with $\alpha_k=1$.

For MP, the corrected velocity is the vector polynomial $\bv_h$ closest to $\bu_h$ in the same norm, under the $2N_q$ conditions
\begin{equation}
  \bar\bn_k\cdot\bv_h(\bxb_k)=0,\qquad \btaub_k\cdot\bv_h(\bxt_k)=\btaub_k\cdot\bu_h(\bxt_k),\qquad k=1,\dots,N_q,
  \label{eq:wall-mp-cond}
\end{equation}
where $\btaub_k$ is the tangent of the wall at $\bxb_k$, so that the normal velocity vanishes at the mapped points and the tangential velocity of the cell is kept at the face points. In SP the second condition is not needed, since the modification is parallel to $\bar\bn_k$ and the tangential velocity at $\bxt_k$ is kept automatically. In MP, instead, each mapped point modifies the velocity along its own normal $\bar\bn_j$, and on a curved wall these normals differ, so that without the second condition the tangential velocity at the face points would also change. On a straight wall all the $\bar\bn_j$ coincide, and the second condition is satisfied automatically, as in SP. The corrected velocity is the stationary point of
\begin{equation}
  \begin{split}
  \mathcal L(\bv^x,\bv^y,\boldsymbol\lambda,\boldsymbol\mu)={}&\tfrac12(\bv^x-\bu^x)^TW(\bv^x-\bu^x)+\tfrac12(\bv^y-\bu^y)^TW(\bv^y-\bu^y)\\
  &+\sum_{j=1}^{N_q}\lambda_j\,\bar\bn_j\cdot\bv_h(\bxb_j)
  +\sum_{j=1}^{N_q}\mu_j\,\btaub_j\cdot\bigl(\bv_h(\bxt_j)-\bu_h(\bxt_j)\bigr).
  \end{split}
\end{equation}
The derivatives with respect to $\bv^x$ and $\bv^y$ give, at any point $\bx$,
\begin{equation}
  \bv_h(\bx)=\bu_h(\bx)-\sum_{j}\Bigl(\lambda_j\,\bphi(\bx)^TW^{-1}\bar\bphi_j\,\bar\bn_j+\mu_j\,\bphi(\bx)^TW^{-1}\phit_j\,\btaub_j\Bigr).
  \label{eq:wall-mp-v}
\end{equation}
Inserting \eqref{eq:wall-mp-v} in the conditions \eqref{eq:wall-mp-cond}, evaluated at $\bx=\bxb_k$ and $\bx=\bxt_k$ respectively, gives the linear system
\begin{equation*}
  \begin{bmatrix}G&H\\H^T&S\end{bmatrix}\begin{bmatrix}\boldsymbol\lambda\\\boldsymbol\mu\end{bmatrix}=\begin{bmatrix}\bar{\mathbf u}_n\\\mathbf 0\end{bmatrix},
  \qquad [\bar{\mathbf u}_n]_k=\bar\bn_k\cdot\bu_h(\bxb_k),
\end{equation*}
with the $N_q\times N_q$ matrices
\begin{equation*}
  G_{kj}=(\bar\bn_k\cdot\bar\bn_j)\,\bar\bphi_k^TW^{-1}\bar\bphi_j,\quad
  H_{kj}=(\bar\bn_k\cdot\btaub_j)\,\bar\bphi_k^TW^{-1}\phit_j,\quad
  S_{kj}=(\btaub_k\cdot\btaub_j)\,\phit_k^TW^{-1}\phit_j.
\end{equation*}
The second block row gives $\boldsymbol\mu=-S^{-1}H^T\boldsymbol\lambda$, and the first one $\boldsymbol\lambda=\bigl(G-HS^{-1}H^T\bigr)^{-1}\bar{\mathbf u}_n$. At the face point $\bxt_k$ the tangential component of $\bu_h-\bv_h$ vanishes by \eqref{eq:wall-mp-cond}, so that $\bu_h(\bxt_k)-\bv_h(\bxt_k)$ is parallel to $\bar\bn_k$, and its component along $\bar\bn_k$ is obtained by projecting \eqref{eq:wall-mp-v} at $\bx=\bxt_k$ on $\bar\bn_k$,
\begin{equation*}
  \begin{split}
  \bar\bn_k\cdot\bigl(\bu_h(\bxt_k)-\bv_h(\bxt_k)\bigr)&=\sum_j\bigl(A_{kj}\lambda_j+B_{kj}\mu_j\bigr)
  =\bigl[\bigl(A-BS^{-1}H^T\bigr)\boldsymbol\lambda\bigr]_k\\
  &=\bigl[\bigl(A-BS^{-1}H^T\bigr)\bigl(G-HS^{-1}H^T\bigr)^{-1}\bar{\mathbf u}_n\bigr]_k,
  \end{split}
\end{equation*}
with $A_{kj}=(\bar\bn_k\cdot\bar\bn_j)\,\phit_k^TW^{-1}\bar\bphi_j$ and $B_{kj}=(\bar\bn_k\cdot\btaub_j)\,\phit_k^TW^{-1}\phit_j$. Hence, writing the $k$-th component of the matrix-vector product,
\begin{equation*}
  \bv_h(\bxt_k)=\bu_h(\bxt_k)-\Bigl(\sum_j\alpha_{kj}\,\bar\bn_j\cdot\bu_h(\bxb_j)\Bigr)\bar\bn_k,
\end{equation*}
with
\begin{equation}
  [\alpha_{kj}]=\bigl(A-BS^{-1}H^T\bigr)\bigl(G-HS^{-1}H^T\bigr)^{-1}.
  \label{eq:wall-mp}
\end{equation}
For $N_q=1$, $H=B=0$ since $\bar\bn_1\cdot\btaub_1=0$, and $\alpha_{11}=A_{11}/G_{11}$ is the SP coefficient.

Both corrections thus modify only the normal velocity at the face point, and we write the corrected velocity as $\bu_h(\bxt_k)-c_k\bar\bn_k$, with
\begin{equation}
  c_k=\alpha_k\,\bar\bn_k\cdot\bu_h(\bxb_k)\quad\text{(SP, with $\alpha_k=1$ for \SB)},\qquad c_k=\sum_j\alpha_{kj}\,\bar\bn_j\cdot\bu_h(\bxb_j)\quad\text{(MP)}.
  \label{eq:wall-c}
\end{equation}
The state at the face point is then $(p_s,\bu_s)=\bigl(p_h(\bxt_k),\ \bu_h(\bxt_k)-c_k\bar\bn_k\bigr)$, with the pressure and the tangential velocity of the cell.
The ghost state is the mirror of the interior state $\bigl(p_h(\bxt_k),\bu_h(\bxt_k)\bigr)$ with respect to the face state, so that the average of the two is the face state,
\begin{equation*}
  p_g=2p_s-p_h(\bxt_k)=p_h(\bxt_k),\qquad
  \bu_g=2\bu_s-\bu_h(\bxt_k)=\bu_h(\bxt_k)-2c_k\bar\bn_k,
\end{equation*}
and the Rusanov flux \eqref{eq:rusanov} gives
\begin{equation}
  \hat F^p=\bigl(\bu_h(\bxt_k)-c_k\bar\bn_k\bigr)\cdot\tilde\bn,\qquad
  \hat{\mathbf F}^{\bu}=p_h(\bxt_k)\,\tilde\bn+c_k\,\bar\bn_k .
  \label{eq:wall-flux}
\end{equation}
On a wall parallel to the face, $\bar\bn_k=\tilde\bn$, the correction only changes the normal velocity in the normal direction, and the wall is the 1D wall applied on each line. On a tilted wall the term $c_k\bar\bn_k$ of the momentum flux has a component along the face, and this is where the difficulties of the following sections originate.

\subsection{Discrete operator of linear acoustics with the slip wall}
\label{sec:wall-operator}

In each cell the pressure and the velocity are $p_h=\bphi^T\mathbf p$, $u^x_h=\bphi^T\bu^x$ and $u^y_h=\bphi^T\bu^y$, and the coefficients of the discrete state of the cell are collected in $\mathbf U_{ij}=(\mathbf p,\bu^x,\bu^y)^T\in\mathbb R^{3D}$. With the matrices $Q_x$, $Q_y$, $\tilde\Phi_f$, $\bar\Phi_f$, $\mathcal W_f$ and $E_f$ of \cref{sec:single-cell}, the weak form \eqref{eq:dg-weak} of \eqref{eq:acoustics} reads
\begin{equation}
  (I_3\otimes M)\,\dot{\mathbf U}_{ij}=\KK^{\rm vol}\,\mathbf U_{ij}-\sum_f\bigl(I_3\otimes\tilde\Phi_f\mathcal W_f\bigr)\hat{\mathbf F}_f,\qquad
  \KK^{\rm vol}=\begin{bmatrix}0&Q_x&Q_y\\Q_x&0&0\\Q_y&0&0\end{bmatrix},
  \label{eq:acoustic-cell}
\end{equation}
where $\hat{\mathbf F}_f\in\mathbb R^{3N_q}$ collects the three components of the numerical flux at the $N_q$ points of the face $f$, and the sum runs over the four faces. The Rusanov flux \eqref{eq:rusanov} is linear in the two traces, so that each face adds to \eqref{eq:acoustic-cell} a term linear in the coefficients, of one of three kinds.
On an interior face, shared with the active neighbour $b$, the flux gives the self term $-\tfrac12\bigl(A_{\tilde\bn}+I_3\bigr)\otimes E_f$, which multiplies $\mathbf U_{ij}$, and the coupling term $-\tfrac12\bigl(A_{\tilde\bn}-I_3\bigr)\otimes\tilde\Phi_f\mathcal W_f(\tilde\Phi_f^{b})^T$, which multiplies the coefficients $\mathbf U_b$ of the neighbour, with $\tilde\Phi_f^{b}$ the basis of the neighbour at the same points. On a face with exact data, the exterior state does not depend on the cell: the self term is the same, $-\tfrac12\bigl(A_{\tilde\bn}+I_3\bigr)\otimes E_f$, and the data only add a forcing, which does not affect the stability. On a wall face, the flux \eqref{eq:wall-flux} depends on the cell also through the coefficients $c_k$ of \eqref{eq:wall-c}, which are linear in the velocity at the mapped points. Computing the integral of the flux with the Gauss points of the face, the self term of a wall face is
\begin{equation}
  \KK^{w}_f\,\mathbf U_{ij}=-\begin{bmatrix}0&\tilde n_xE_f&\tilde n_yE_f\\\tilde n_xE_f&0&0\\\tilde n_yE_f&0&0\end{bmatrix}\mathbf U_{ij}
  +\sum_{k=1}^{N_q}w_k\,c_k\begin{bmatrix}(\bar\bn_k\cdot\tilde\bn)\,\phit_k\\-\bar n_{x,k}\,\phit_k\\-\bar n_{y,k}\,\phit_k\end{bmatrix},
  \label{eq:K-wall}
\end{equation}
where $c_k=\sum_{j}\alpha_{kj}\,\bar\bn_j\cdot\bu_h(\bxb_j)$, with $\alpha_{kj}=\alpha_k\delta_{kj}$ for \SB\ and SP.
The first term is the flux of the uncorrected cell state, and the second one is the correction: the normal velocity of the cell at the mapped points gives the $c_k$, which act on the pressure through the alignment $\bar\bn_k\cdot\tilde\bn$ and on the velocity along $\bar\bn_k$. Since the $c_k$ are linear in the coefficients, $\KK^w_f$ is a matrix, and its correction part has rank at most $N_q$.
To assemble the operator of a set of $N_e$ active cells, numbered $i=1,\dots,N_e$, we stack their coefficients in the global vector $\mathbf U=(\mathbf U_1,\dots,\mathbf U_{N_e})^T\in\mathbb R^{3DN_e}$. The semi-discrete scheme is $\MM\dot{\mathbf U}=\KK\mathbf U$, with the block diagonal mass matrix $\MM=I_{N_e}\otimes I_3\otimes M$ and the block matrix $\KK=[\KK_{ib}]_{i,b=1}^{N_e}$, whose blocks are
\begin{equation}
  \begin{aligned}
  \KK_{ii}&=\KK^{\rm vol}-\tfrac12\sum_{f\subset\partial K_i\setminus\tilde\Gamma_w}\bigl(A_{\tilde\bn}+I_3\bigr)\otimes E_f+\sum_{f\subset\partial K_i\cap\tilde\Gamma_w}\KK^{w}_f,\\
  \KK_{ib}&=-\tfrac12\bigl(A_{\tilde\bn}-I_3\bigr)\otimes\tilde\Phi_f\mathcal W_f(\tilde\Phi_f^{b})^T\qquad\text{if $i$ and $b$ share the face $f$},
  \end{aligned}
  \label{eq:acoustic-global}
\end{equation}
and $\KK_{ib}=0$ otherwise, where $\tilde\Gamma_w$ denotes the computational faces on which the wall condition is imposed. Unlike for advection, both $\KK_{ib}$ and $\KK_{bi}$ are nonzero, since one wave crosses every interior face in each direction: $\KK$ is not block triangular, and the spectrum of $\MM^{-1}\KK$ must be computed on the whole set of cells. 

\subsection{A wall at constant distance}
\label{sec:wall-const}

We consider one cell $K=[0,1]^2$ with the wall parallel to its left face, at the distance $d$, so that $\bxt_k=(0,\eta_k)$ is mapped to $\bxb_k=(d,\eta_k)$ and $\bar\bn_k=\tilde\bn=(-1,0)$, and with the exact condition on the other three faces. In 1D, on the cell $[0,1]$ with the wall at the distance $d$ from its left end, $\tilde n=\bar n=-1$, the flux \eqref{eq:wall-flux} reads $\hat F^p=-\bigl(u_h(0)-\alpha\,u_h(d)\bigr)$ and $\hat F^u=-p_h(0)+\alpha\,u_h(d)$, and with $\bpsi_0=\bpsi(0)$, $\bpsi_d=\bpsi(d)$ the wall block on the coefficients $(\mathbf p,\bu)$ is
\begin{equation}
  L^{w}_{1D}=\begin{bmatrix}0&\bpsi_0\bigl(\bpsi_0-\alpha\,\bpsi_d\bigr)^T\\[2pt]\bpsi_0\bpsi_0^T&-\alpha\,\bpsi_0\bpsi_d^T\end{bmatrix}.
  \label{eq:wall-1d}
\end{equation}

\begin{proposition}[slip wall at constant distance]\label{prop:wall-const}
With the tensor-product basis $\bphi=\bpsi(x)\otimes\bpsi(y)$ and $W=W_1\otimes W_1$, the block \eqref{eq:K-wall} of the left face is
\begin{equation}
  \KK^{w}_f=\begin{bmatrix}L^{w}_{1D}\otimes M_1&\mathbf 0\\\mathbf 0&0\end{bmatrix}\qquad\text{on }\bigl((\mathbf p,\bu^x),\ \bu^y\bigr),
  \label{eq:wall-const}
\end{equation}
with $\alpha=1$ for \SB\ and $\alpha$ the 1D coefficient \eqref{eq:sp-1d} with $\bar x_k=d$ for ROD (MP reduces to SP). The correction is therefore the 1D wall applied on each horizontal line, and it does not act on $u^y$.
\end{proposition}

\begin{proof}
The SP coefficients are given by \eqref{eq:sp-1d} and do not depend on $k$, and MP reduces to SP as in \cref{prop:mp-constant}. On the left face, $\tilde\bn=\bar\bn_k=(-1,0)$, so that $\bar\bn_k\cdot\tilde\bn=1$, $c_k=-\alpha\,u^x_h(d,\eta_k)$, and \eqref{eq:K-wall} becomes
\begin{equation*}
  \KK^{w}_f\,\mathbf U_{ij}=\begin{bmatrix}E_f\bu^x\\E_f\mathbf p\\\mathbf 0\end{bmatrix}
  -\alpha\sum_kw_k\begin{bmatrix}\phit_k\\\phit_k\\\mathbf 0\end{bmatrix}u^x_h(d,\eta_k).
\end{equation*}
Since the $p+1$ Gauss points integrate exactly the products of two polynomials of degree $p$ in $y$,
\begin{equation*}
  \begin{aligned}
  E_f&=\sum_kw_k\,\bphi(0,\eta_k)\bphi(0,\eta_k)^T=\bpsi_0\bpsi_0^T\otimes M_1,\\
  \sum_kw_k\,\phit_k\,u^x_h(d,\eta_k)&=\sum_kw_k\,\bphi(0,\eta_k)\bphi(d,\eta_k)^T\bu^x=\bigl(\bpsi_0\bpsi_d^T\otimes M_1\bigr)\bu^x,
  \end{aligned}
\end{equation*}
which gives \eqref{eq:wall-const}.
\end{proof}

Unlike for advection (\cref{prop:constant-1d}), this does not make the operator of the cell the 1D one. The wall itself, however, is the 1D wall, with the same coefficients.

The spectra confirm the 1D behaviour on the wall outside the cell, $d\in[-1,0]$. All the ROD corrections are stable with $\text{CFL}=1$ up to $\PP{3}$, and \RODL\ also with $\PP{4}$, while \RODE\ with $\PP{4}$ is unstable. The \SB\ correction is limited by the time step, down to $\text{CFL}=0.97$, $0.33$, $0.1$ and $0.017$ for $\PP{1},\dots,\PP{4}$ at $d=-1$. With the wall inside the cell, all the corrections are stable only close to the face from $\PP{2}$ on. At constant distance the slip wall is therefore as stable as the 1D inlet, and this is the reference for the next section, where the wall is tilted with respect to the mesh.

\subsection{A tilted straight wall}
\label{sec:wall-tilted}

\begin{figure}
\centering
\includegraphics[width=0.8\linewidth]{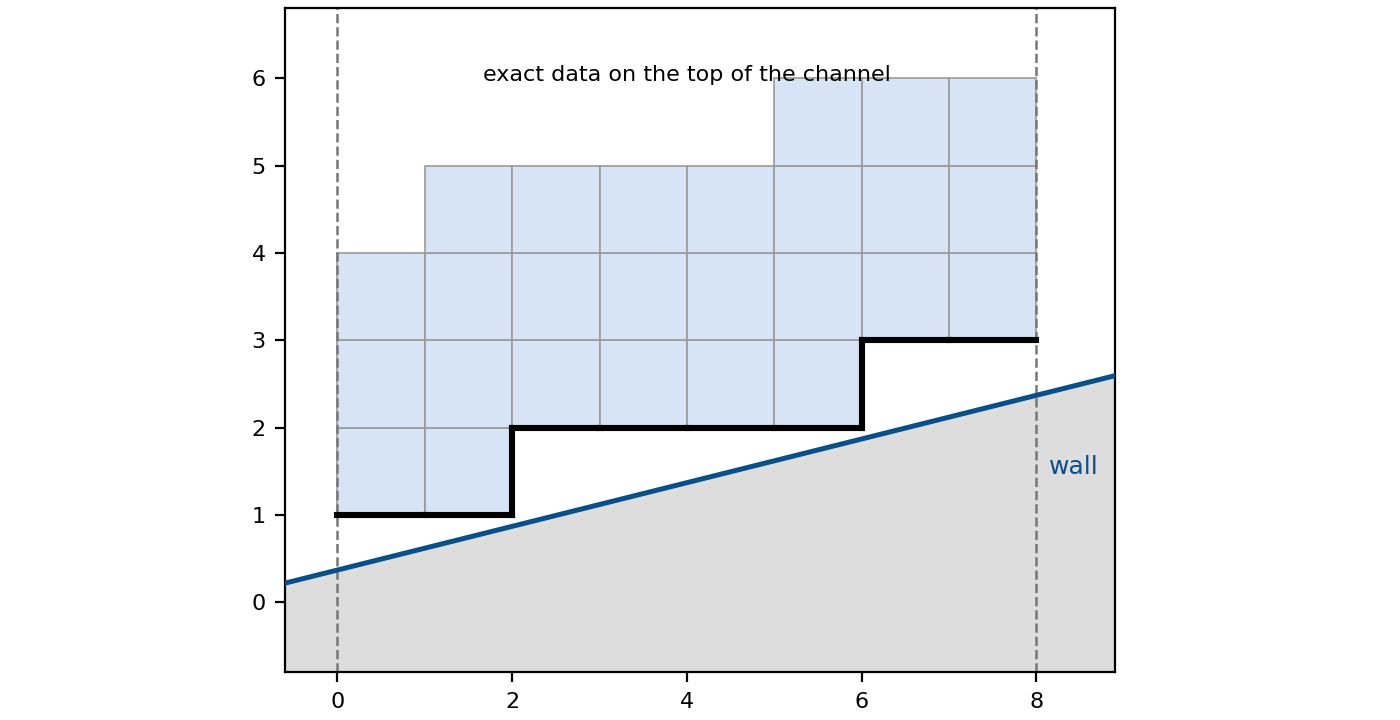}
\caption{The periodic channel of \cref{sec:wall-tilted}, $N=8$, $\theta=14^\circ$: one period (between the dashed lines) of the active cells (blue), the wall (blue line) and the computational wall (black), with the exact condition on the top of the channel.}
\label{fig:channel}
\end{figure}

As for linear advection, the constant-distance picture is lost as soon as the wall is tilted with respect to the mesh, and a single straight wall is enough to show it. Since the acoustic operator couples all the cells, the configuration must contain a set of cells and a portion of the wall. To do this, we use a periodic channel (\cref{fig:channel}). Lengths are measured in units of the cell, $h=1$, so that the vertices of the mesh have integer coordinates. The wall is the line $y=y_0+x\tan\theta$, with the fluid above it and $y_0=0.37$, so that the wall does not pass through the vertices of the mesh, and its slope is rational, $\tan\theta=k/N$ with integers $k$ and $N$: over $N$ cells along $x$ the wall rises by exactly $k$ cells, and the channel is periodic along the wall, the point $(x+N,y+k)$ being identified with $(x,y)$, so that the geometry repeats exactly and the wall never meets the sides of the domain. The computational domain is a strip along the wall. Its active cells are the cells that lie entirely in the fluid, above the wall, and entirely below a second line, parallel to the wall at the distance $N/2$ above it. The faces of the active cells towards the cells cut by the wall form the computational wall, a staircase, and carry the slip-wall boundary treatment, while on the faces towards the cells cut by the upper line the exact condition is imposed. We take $N=8$ and $k=2$, i.e.\ $\theta=14^\circ$, so that the computational wall is a staircase with two steps per period. Besides the corrections of \cref{sec:wall-condition} with the geometry-normal map, we consider two boundary treatments without correction of the location, $\bxb_k=\bxt_k$: the mirror about the normal of the face, which is the fitted wall on the staircase, and the mirror about the normal of the real wall, in the spirit of the boundary correction of~\cite{krivodonova2006high} developed to impose slip-wall conditions on body-fitted linear meshes.

\begin{table}
\centering
\caption{Periodic channel with $\theta=14^\circ$ (\cref{fig:channel}): largest real part of the eigenvalues of the whole operator, computed with lengths in units of the cell ($h=1$), i.e.\ $h\max\operatorname{Re}\lambda$ for a mesh of size $h$, for one period of $N=8$ and $N=16$ cells. Bold: positive; values below $10^{-10}$ in absolute value are reported as $0$.}
\label{tab:channel}
\small
\begin{tabular}{lccccc}
\toprule
 & \multicolumn{2}{c}{$\PP{1}$} & \multicolumn{2}{c}{$\PP{2}$} & $\PP{3}$\\
\cmidrule(lr){2-3}\cmidrule(lr){4-5}\cmidrule(lr){6-6}
boundary treatment & $N=8$ & $N=16$ & $N=8$ & $N=16$ & $N=8$\\
\midrule
fitted staircase       & $-6.4\cdot10^{-3}$ & $-1.8\cdot10^{-4}$ & $0$ & $0$ & $0$\\
normal of the wall & $-5.8\cdot10^{-4}$ & $\mathbf{4.1\cdot10^{-3}}$ & $\mathbf{0.60}$ & $\mathbf{0.60}$ & $\mathbf{2.48}$\\
\SB     & $-7.3\cdot10^{-4}$ & $\mathbf{9.6\cdot10^{-3}}$ & $\mathbf{1.02}$ & $\mathbf{1.02}$ & $\mathbf{3.18}$\\
\RODESP & $-3.0\cdot10^{-4}$ & $\mathbf{1.8\cdot10^{-2}}$ & $\mathbf{0.12}$ & $\mathbf{0.12}$ & $\mathbf{0.64}$\\
\RODLSP & $\mathbf{5.8\cdot10^{-4}}$ & $\mathbf{2.1\cdot10^{-2}}$ & $\mathbf{0.13}$ & $\mathbf{0.13}$ & $\mathbf{0.60}$\\
\RODEMP & $\mathbf{5.2\cdot10^{-3}}$ & $\mathbf{1.8\cdot10^{-2}}$ & $\mathbf{0.17}$ & $\mathbf{0.17}$ & $\mathbf{0.86}$\\
\RODLMP & $\mathbf{6.8\cdot10^{-3}}$ & $\mathbf{2.1\cdot10^{-2}}$ & $\mathbf{0.20}$ & $\mathbf{0.20}$ & $\mathbf{0.77}$\\
\bottomrule
\end{tabular}
\end{table}

\Cref{tab:channel} gives the largest real part of the eigenvalues of the whole operator. As in \cref{sec:single-cell}, the lengths are in units of the cell, $h=1$: on a mesh of size $h$ the eigenvalues are $\lambda/h$, so that the table gives $h\max\operatorname{Re}\lambda$. The fitted staircase is stable at all the degrees, as predicted by the energy balance of \cref{sec:wall-energy}, and every other boundary treatment is unstable. With $\PP{2}$ and $\PP{3}$ all of them have eigenvalues with real parts of size $\mathcal{O}(1)$ in units of the cell. The growth rate of the instability is therefore proportional to $1/h$, it increases without bound as the mesh is refined, and it cannot be removed by the time step. With $\PP{1}$ the instability is weaker and depends on the length of the wall: starting from $16$ cells, which contain longer modes along the wall, every boundary treatment other than the fitted staircase is unstable. The positive real parts are, however, of the order of $10^{-3}$--$10^{-2}$, i.e.\ two orders of magnitude smaller than with $\PP{2}$ and $\PP{3}$, so that $\PP{1}$ is a borderline case, which could be stabilized much more easily than the higher degrees. The modified boundary condition that uses the normal of the wall without moving the points is already unstable, which shows that the instability is not caused only by the extrapolation of the polynomial to the mapped points, but also by the mismatch between the normals of the faces of the staircase and of the wall. Contrary to what was observed for linear advection, the face-orthogonal map is no better for the slip wall: with this map all the corrections are unstable as well, already with $\PP{1}$ on the shortest channel, $N=8$, and with real parts of the same size as in \cref{tab:channel} at the higher degrees. The energy balance of the next section makes this precise.

\subsection{Energy balance of the slip wall}
\label{sec:wall-energy}

To understand why the tilted wall is unstable we compute the energy balance of the semi-discrete scheme. The acoustic energy of the discrete solution is $E=\sum_KE_K$, with $E_K=\tfrac12\int_K(p_h^2+|\bu_h|^2)\,d\bx$. On a cell, taking $p_h$, $u^x_h$ and $u^y_h$ as test functions in the DG equations and summing, the volume term becomes $\int_K\nabla\cdot(p_h\bu_h)\,d\bx$. Writing the balance in this form assumes that the integration by parts of the volume operator holds exactly at the discrete level, i.e.\ that the volume operator has the summation-by-parts property~\cite{gassner2013skew}. This is not true for a DG discretization in general: for nonlinear fluxes, or when the volume integrals are not computed exactly, the discrete volume term does not reduce to a boundary term, and summation-by-parts operators, e.g.\ in split form, are needed to recover it. In our setting it holds, because the flux is linear and the Gauss rules with $p+1$ points per direction integrate exactly the products of two polynomials of the cell: the integrands $\partial_x(p_hu^x_h)$, $\partial_y(p_hu^y_h)$ and $p_h\,\bu_h\cdot\tilde\bn$ have degree at most $2p$ in each variable. The same holds on the extended cells of \cref{sec:gap}, whose quadrature rules are exact. The volume term then cancels the trace of the energy flux on the faces, and the balance isolates the contribution of the boundary, which is the object of this section; in a setting without this property, the volume would contribute additional terms, which would have to be controlled separately from the wall. The energy of a cell therefore changes only through its faces,
\begin{equation}
  \frac{dE_K}{dt}=\sum_f\sum_kw_k\,\mathcal B_k,\qquad
  \mathcal B_k=p_k\,\bu_k\cdot\tilde\bn-p_k\,\hat F^p_k-\bu_k\cdot\hat{\mathbf F}^{\bu}_k ,
  \label{eq:energy-point}
\end{equation}
where $p_k$ and $\bu_k$ are the traces of the cell at $\bxt_k$, and each face point contributes the difference between the energy flux of the trace and the energy flux carried by the numerical flux. On an interior face the two contributions add up to $-\tfrac12\bigl([\![p]\!]^2+|[\![\bu]\!]|^2\bigr)$, the dissipation of the Rusanov flux, and on a face with exact data, which are zero in the energy balance, the contribution is $-\tfrac12(p_k^2+|\bu_k|^2)$. Both are non-positive, so that
\begin{equation}
  \frac{dE}{dt}=-\sum_{\text{interior}}\sum_k\frac{w_k}2\Bigl([\![p]\!]^2+\bigl|[\![\bu]\!]\bigr|^2\Bigr)
  -\sum_{\text{exact}}\sum_k\frac{w_k}2\bigl(p_k^2+|\bu_k|^2\bigr)
  +\sum_{\text{wall}}\sum_kw_k\,\mathcal B^{\rm w}_k ,
  \label{eq:energy-mesh}
\end{equation}
and the sign of the energy rate depends only on the wall terms $\mathcal B^{\rm w}_k$.

Inserting the wall flux \eqref{eq:wall-flux} in \eqref{eq:energy-point} gives, at a point of a wall face and for any correction $c_k$,
\begin{equation}
  \mathcal B^{\rm w}_k=c_k\bigl(p_k\,\bar\bn_k\cdot\tilde\bn-\bu_k\cdot\bar\bn_k\bigr)-p_k\,\bu_k\cdot\tilde\bn .
  \label{eq:energy-wall}
\end{equation}
The corrections of \cref{sec:wall-condition} differ from the fitted wall in two ingredients: they use the normal of the real wall instead of that of the face, and they evaluate the normal velocity at the mapped point instead of the face point. We introduce them one at a time.

\emph{Fitted wall.} With the normal of the face and the point on the face, $\bar\bn_k=\tilde\bn$ and $c_k=\bu_k\cdot\tilde\bn$, \eqref{eq:energy-wall} gives
\begin{equation*}
  \mathcal B^{\rm w}_k=-(\bu_k\cdot\tilde\bn)^2\le0:
\end{equation*}
the fitted wall dissipates the normal velocity, and the scheme on the fitted staircase is energy stable.

\emph{Normal of the wall.} We now take the normal of the wall but keep the point on the face, $\bxb_k=\bxt_k$ and $c_k=\bu_k\cdot\bar\bn_k$, so that the normal velocity of the face state vanishes along $\bar\bn_k$ while $p_k$ and $\bu_k\cdot\btaub_k$ are those of the cell. This is the spirit of the boundary treatment of~\cite{krivodonova2006high}, who impose the slip condition on body-fitted meshes with straight faces by reflecting the velocity about the normal of the true boundary at each quadrature point of the face, and it is also the limit of all the corrections when the mapped points approach the face points. Since $\bar\bn_k\ne\tilde\bn$, we decompose the normal of the face and the velocity of the cell in the orthonormal basis $(\bar\bn_k,\btaub_k)$ of the wall,
\begin{equation}
  \tilde\bn=(\bar\bn_k\cdot\tilde\bn)\,\bar\bn_k+(\btaub_k\cdot\tilde\bn)\,\btaub_k,\qquad
  \bu_k\cdot\tilde\bn=(\bu_k\cdot\bar\bn_k)(\bar\bn_k\cdot\tilde\bn)+(\bu_k\cdot\btaub_k)(\btaub_k\cdot\tilde\bn).
  \label{eq:energy-basis}
\end{equation}
Inserting $c_k=\bu_k\cdot\bar\bn_k$ in \eqref{eq:energy-wall} and using \eqref{eq:energy-basis} in the last term, the normal parts cancel and
\begin{equation}
  \mathcal B^{\rm w}_k=-(\bu_k\cdot\bar\bn_k)^2-p_k\,(\bu_k\cdot\btaub_k)\,(\btaub_k\cdot\tilde\bn) .
  \label{eq:energy-normal}
\end{equation}
The first term dissipates the normal velocity, as for the fitted wall. The second one has no sign, since the product of the pressure and of the tangential velocity can take any sign, and \eqref{eq:energy-mesh} gives no bound on it. It vanishes only when $\btaub_k\cdot\tilde\bn=0$, i.e.\ for a wall parallel to the face. This term is the energy flux of the exact solution through the face: by \eqref{eq:energy-basis}, $p\,\bu\cdot\tilde\bn=p(\bu\cdot\bar\bn)(\bar\bn\cdot\tilde\bn)+p(\bu\cdot\btaub)(\btaub\cdot\tilde\bn)$, and on the wall $\bu\cdot\bar\bn=0$, so that what remains is the work of the pressure along the wall that crosses the computational face, which is not a streamline. Any consistent flux must reproduce it, so that it cannot be removed by the choice of the normal in the Riemann solver. In the continuous problem this energy may enter or leave the region between the computational and the real boundary, which in our case is not part of the computational domain. The size of the term is the size of $\btaub_k\cdot\tilde\bn$, the sine of the angle between the face and the wall. On the body-fitted meshes of~\cite{krivodonova2006high} the straight faces are chords of the boundary, the angle is $\mathcal{O}(\kappa h)$, with $\kappa$ the curvature, and the term is a small perturbation of the dissipative fitted wall. On a staircase the angle is $\mathcal{O}(1)$, and close to $90^\circ$ on the steps of the staircase, so that the same boundary treatment carries an indefinite energy flux of the size of the solution itself on every step.

\emph{Polynomial correction.} Finally, we add the correction of the location, with the mapped points $\bxb_k\ne\bxt_k$, and write the correction as the normal velocity of the cell at the face point minus a remainder, $c_k=\bu_k\cdot\bar\bn_k-e_k$, with
\begin{equation}
  e_k=\bu_k\cdot\bar\bn_k-\alpha_k\,\bar\bn_k\cdot\bu_h(\bxb_k),
  \label{eq:energy-e}
\end{equation}
with $\alpha_k=1$ for \SB\ and $\alpha_k=\aSPW(\bxt_k)$ for SP, and in the same way for MP, with $\sum_j\alpha_{kj}\,\bar\bn_j\cdot\bu_h(\bxb_j)$ in the last term. The remainder $e_k$ is the change of the normal velocity between the face and the wall estimated by the correction, and it vanishes when $\bxb_k=\bxt_k$. The same steps as above give
\begin{equation}
  \mathcal B^{\rm w}_k=-(\bu_k\cdot\bar\bn_k)^2\;\underbrace{-\,p_k\,(\bu_k\cdot\btaub_k)\,(\btaub_k\cdot\tilde\bn)}_{T_1}\;\underbrace{+\,e_k\bigl(\bu_k\cdot\bar\bn_k-p_k\,\bar\bn_k\cdot\tilde\bn\bigr)}_{T_2} .
  \label{eq:energy-rod}
\end{equation}
The correction of the location adds the term $T_2$ to the term $T_1$ of the normal, and $T_2$ has no sign either. Already on a wall parallel to the face at constant distance, where $T_1=0$, $T_2$ has an indefinite sign, and it is responsible for the instabilities of the 1D wall in \cref{sec:wall-const}. To control these terms, two properties of the boundary treatment are therefore convenient. The term $T_1$ contains explicitly the product $\btaub_k\cdot\tilde\bn$, so that the normal of the computational boundary should be close to that of the real wall. The term $T_2$ contains the remainder $e_k$ of the polynomial correction, which depends on the distance $|\bxb_k-\bxt_k|$ between the mapped point and the face point, so that the points where the condition is imposed should be close to the real wall. On a staircase neither property holds, and this is consistent with the instabilities of the channel of \cref{sec:wall-tilted}.

\section{A gap method for DG discretizations of hyperbolic systems}
\label{sec:gap}

The energy balance of \cref{sec:wall-energy} points to the geometry of the computational boundary as the origin of the instability, rather than to the particular correction: on a staircase the faces are not aligned with the wall and are far from it on the scale of the cell, and no choice of the coefficients of the correction can change this. For this reason, rather than designing yet another correction on the staircase, in this section we change the paradigm with which the computational boundary is built. We fill the region between the staircase and the wall, which is not part of the computational domain in the classical approach, with extensions of the boundary cells, and we impose the boundary condition on a polygonal boundary inscribed in the wall. We call this construction the gap method, and gap-\SB\ or gap-ROD the corresponding corrections. This paradigm has been developed for continuous Galerkin discretizations of elliptic problems, with the primary goal of achieving optimal accuracy with Neumann boundary conditions~\cite{collins2026gap,gorgi2026isogeometric,li2026gap}. Here we explore its development for discontinuous Galerkin discretizations of hyperbolic systems, to mitigate the staircase instabilities of \cref{sec:acoustics}.

\subsection{Extended cells}
\label{sec:gap-construction}

\begin{figure}
\centering
\includegraphics[width=\linewidth]{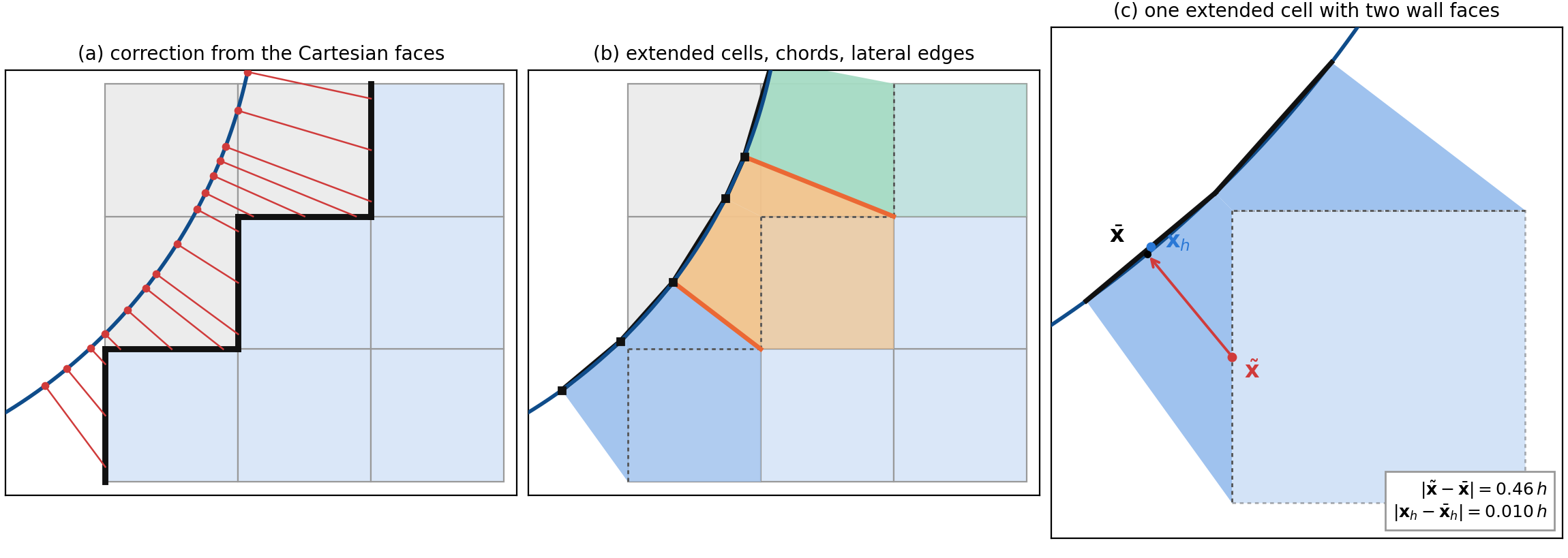}
\caption{The gap method. (a) Classical correction: the boundary condition is imposed on the faces of the staircase (black), and each face point is mapped to the wall (red). (b) Extended cells: each computational face is extended by the quadrilateral between the face and the chord joining the projections of its vertices (black squares); the colours show the extended cells, the orange segments the lateral edges between two extended cells, and the dashed lines the faces of the staircase, which are now inside the extended cells. (c) One extended cell with two computational faces: a point $\tilde\bx$ of a face of the staircase is at distance $0.46\,h$ from its projection $\bar\bx$, while a point of the chord is at distance $0.010\,h$ from its projection.}
\label{fig:gap}
\end{figure}

Let $e=(A,B)$ be a computational face of an active cell $K$, with the vertices ordered counter-clockwise on $K$, and let $P$ denote the closest-point projection on the wall. The \emph{gap quadrilateral} of the face is
\begin{equation}
  K_e=\bigl(A,\ P(A),\ P(B),\ B\bigr),
  \label{eq:gap-quad}
\end{equation}
and the \emph{extended cell} is the union of $K$ with the gap quadrilaterals of all its computational faces (\cref{fig:gap}). A cell at a convex corner of the staircase has two computational faces and two quadrilaterals, which share the segment $(A,P(A))$ of their common vertex. When the vertex $A$ is shared by the computational faces of two different cells, the segment $(A,P(A))$ separates the quadrilaterals of the two cells. The extended cells use the polynomial space of the cell, with the same basis, so that the method has no new unknowns: the polynomial of $K$ is simply extrapolated to its quadrilaterals. Since the projections of the vertices of the staircase are shared by the neighbouring quadrilaterals, the extended cells tile the region between the active cells and the polygonal line made of the segments $P(A)P(B)$. We call these segments the \emph{chords}, since their end points lie on the wall.

\begin{figure}
\centering
\begin{tikzpicture}[scale=3.4,>=stealth]
  \coordinate (A)  at (0,0);      \coordinate (B)  at (1,0);      \coordinate (Ap) at (-1,0);
  \coordinate (PA) at (0.094,-0.407); \coordinate (PB) at (0.942,-0.376); \coordinate (PAp) at (-0.579,-0.684);
  \fill[blue!7] (-1,0) -- (0,0) -- (0,1) -- (-1,1) -- cycle;
  \fill[blue!7] (Ap) -- (PAp) -- (PA) -- (A) -- cycle;
  \fill[blue!22] (0,0) -- (1,0) -- (1,1) -- (0,1) -- cycle;
  \fill[blue!32] (A) -- (PA) -- (PB) -- (B) -- cycle;
  \fill[gray!25] plot[domain=66:134,samples=90] ({0.6+2.25*cos(\x)},{-2.6+2.25*sin(\x)}) -- (-1.25,-1.15) -- (1.3,-1.15) -- cycle;
  \draw[very thick,blue!50!black] plot[domain=66:134,samples=90] ({0.6+2.25*cos(\x)},{-2.6+2.25*sin(\x)});
  \node[blue!50!black] at (1.22,-0.62) {$\bar\Gamma$};
  \draw[dashed,gray!70!black] (Ap) -- (A) -- (B);
  \draw[gray!60] (-1,0) -- (-1,1) -- (0,1);
  \draw[gray!60] (Ap) -- (PAp);
  \draw[red!35,thick] (PAp) -- (PA);
  \draw[thick] (B) -- (1,1) -- (0,1);
  \draw[orange!90!black,line width=1.6pt] (0,0) -- (0,1);
  \draw[orange!90!black,line width=1.6pt] (A) -- (PA);
  \draw[red!80!black,line width=1.6pt] (PA) -- (PB);
  \draw[thick] (PB) -- (B);
  \fill (0.518,-0.3915) circle (0.018);
  \fill[red!80!black] (0.5165,-0.3516) circle (0.018);
  \draw[->,red!80!black] (0.518,-0.3915) -- (0.5165,-0.358);
  \node at (0.44,-0.47) {\small $\bx_k$};
  \node at (0.40,-0.29) {\small $\bxb_k$};
  \draw[->,thick] (0.72,-0.386) -- (0.724,-0.586) node[left] {$\bn_c$};
  \draw[->,thick] (0.047,-0.2035) -- (-0.153,-0.249) node[left] {$\bn_l$};
  \draw[->,thick] (0,0.62) -- (-0.2,0.62) node[above] {$\tilde\bn$};
  \node[circle,fill,inner sep=1pt,label={[label distance=-1pt]above right:$A$}] at (A) {};
  \node[circle,fill,inner sep=1pt,label=above right:$B$] at (B) {};
  \node[circle,fill,inner sep=1pt] at (PA) {}; \node at (0.06,-0.56) {$P(A)$};
  \node[circle,fill,inner sep=1pt] at (PB) {}; \node at (1.02,-0.53) {$P(B)$};
  \node at (0.5,0.5) {$K$};
  \node at (0.72,-0.16) {$K_e$};
  \node[gray!70!black] at (-0.5,0.5) {$K'$};
  \node[gray!70!black] at (-0.55,-0.16) {$K_{e'}$};
  \draw[red!80!black,line width=1.6pt] (1.25,0.95) -- (1.4,0.95) node[right,black] {\small chord: boundary condition};
  \draw[orange!90!black,line width=1.6pt] (1.25,0.78) -- (1.4,0.78) node[right,black] {\small interior fluxes with $\Omega_{K'}$};
  \draw[dashed,gray!70!black] (1.25,0.61) -- (1.4,0.61) node[right,black] {\small staircase: no flux};
  \draw[thick] (1.25,0.44) -- (1.4,0.44) node[right,black] {\small other Cartesian faces};
\end{tikzpicture}
\caption{One extended cell $\Omega_K=K\cup K_e$ (blue), with its gap quadrilateral $K_e=(A,P(A),P(B),B)$, and the neighbouring extended cell $\Omega_{K'}=K'\cup K_{e'}$ (light blue). The computational faces of the staircase (dashed) are inside the extended cells. The chord $P(A)P(B)$ (red) carries the boundary condition, with the correction from each Gauss point $\bx_k$ to its projection $\bxb_k$ on the wall; the lateral Cartesian face and the lateral edge $(A,P(A))$ (orange) carry the Rusanov flux between $\Omega_K$ and $\Omega_{K'}$.}
\label{fig:gap-cell}
\end{figure}
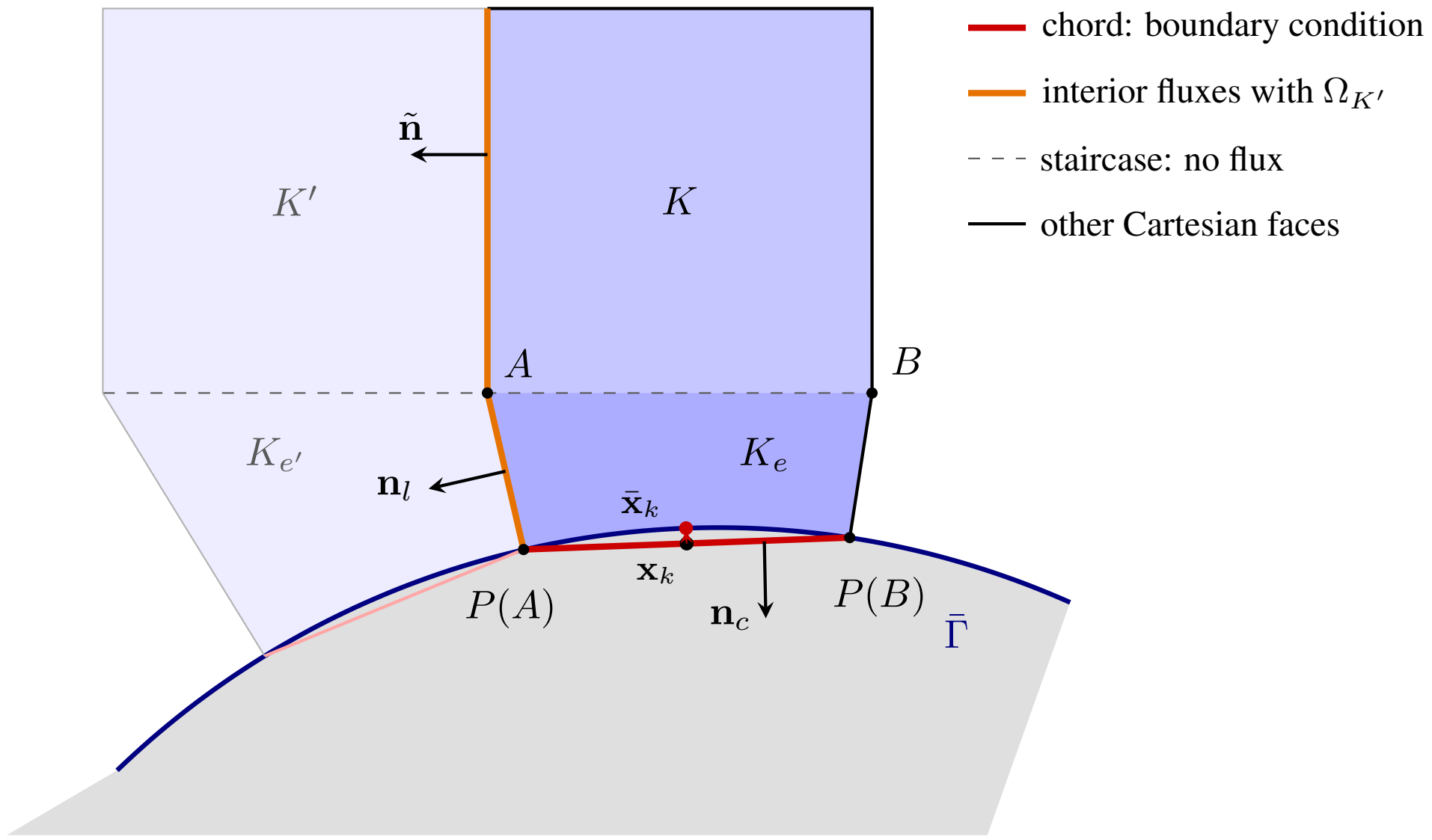

The DG scheme is written on the extended cells. Let $\Omega_K=K\cup\bigcup_eK_e$ be the extended cell of $K$, the union running over the gap quadrilaterals of its computational faces, and consider the system \eqref{eq:system}. The boundary of $\Omega_K$ is made of three kinds of edges (\cref{fig:gap-cell}): the set $\mathcal F_K$ of the Cartesian faces of $K$ that are shared with other active cells or lie on the sides of the box, the set $\mathcal L_K$ of the lateral edges shared with another extended cell, and the set $\mathcal C_K$ of the chords of $K$. The computational faces of the staircase, such as $(A,B)$ in \cref{fig:gap-cell}, are inside $\Omega_K$ and do not appear. For every basis function $\varphi$ of $K$ the scheme reads
\begin{equation}
\begin{aligned}
  \int_{\Omega_K}\varphi\,\partial_t\mathbf U_h\,d\bx
  ={}&\int_{\Omega_K}\bigl(\partial_x\varphi\,\mathbf F(\mathbf U_h)+\partial_y\varphi\,\mathbf G(\mathbf U_h)+\varphi\,\mathbf S\bigr)\,d\bx
  -\sum_{f\in\mathcal F_K}\int_f\varphi\,\hat{\mathbf F}\bigl(\mathbf U_h^K,\mathbf U_h^{K_f},\tilde\bn_f\bigr)\,ds\\
  &-\sum_{l\in\mathcal L_K}\int_l\varphi\,\hat{\mathbf F}\bigl(\mathbf U_h^K,\mathbf U_h^{K_l},\bn_l\bigr)\,ds
  -\sum_{c\in\mathcal C_K}\int_c\varphi\,\hat{\mathbf F}\bigl(\mathbf U_h^K,\mathbf U^{g},\bn_c\bigr)\,ds ,
\end{aligned}
\label{eq:gap-weak}
\end{equation}
where $\mathbf U_h^K$ is the polynomial of $K$, evaluated also on its quadrilaterals $K_e$, $\mathbf U_h^{K_f}$ and $\mathbf U_h^{K_l}$ are the polynomials of the neighbouring cells across the face $f$ and the lateral edge $l$, or the exact data on the sides of the box, and $\tilde\bn_f$, $\bn_l$, $\bn_c$ are the outward normals of $\Omega_K$. The first two sums are the interior fluxes of the extended cell, with the Rusanov flux \eqref{eq:rusanov}. A lateral edge between two quadrilaterals of the same cell, which occurs at a convex corner of the staircase, is inside the extended cell and carries no flux. The last sum is the boundary condition: the ghost state $\mathbf U^g$ is built at each Gauss point $\bx_k$ of the chord with the correction of \cref{sec:wall-condition}, from $\bx_k$ to its projection $\bxb_k=P(\bx_k)$ on the wall, and with the normal $\bn_c$ of the chord in the flux. With $\mathbf U_h^K=\bphi^T\mathbf U_K$, where $\mathbf U_K$ collects the coefficients of the cell as in \cref{sec:wall-operator}, the left-hand side of \eqref{eq:gap-weak} involves the mass matrix of the extended cell,
\begin{equation}
  M_{\rm ext}=M+\sum_{e}\int_{K_e}\bphi\,\bphi^T\,d\bx ,
\end{equation}
which replaces $M$ in the time update and in the initial projection. Since the quadrilaterals are added to a whole Cartesian cell, the extended cell is never smaller than a Cartesian cell: unlike cut-cell methods, the gap method does not create arbitrarily small cells, and it does not suffer from the small cell problem. To integrate exactly on these general quadrilaterals, which are mapped bilinearly on the unit square, the areas are integrated with $(2p+2)^2$ Gauss points, and the chords and lateral edges with $2p+1$ points, so that the divergence theorem holds at the discrete level on the extended cells, as on the Cartesian cells.

\subsection{Energy balance of the slip wall with gap elements}
\label{sec:gap-energy}

Since the divergence theorem holds on each extended cell, the energy balance of \cref{sec:wall-energy} applies to the gap method with the energy of the extended cells,
\begin{equation}
  E=\tfrac12\sum_K\int_{\Omega_K}\bigl(p_h^2+|\bu_h|^2\bigr)\,d\bx ,
\end{equation}
where $\Omega_K$ is the extended cell of $K$. The polynomial that each cell extrapolates to its quadrilaterals is therefore part of the energy, and is controlled by it. The interior faces now include the lateral edges, which carry the Rusanov flux and contribute the dissipation $-\tfrac12([\![p]\!]^2+|[\![\bu]\!]|^2)$ to \eqref{eq:energy-mesh}, and the wall terms are those of \eqref{eq:energy-rod} with the chord in place of the face, $\tilde\bn=\bn_c$ and $\bxt_k=\bx_k$ on the chord. The two terms $T_1$ and $T_2$ are both much smaller than on the staircase, as the following geometric facts show.

Let the wall be smooth enough with curvature at most $\kappa$, and let $\ell$ be the length of a chord, with $\ell\le Ch$. The direction of the chord is that of the tangent of the wall at some point of the arc between its end points, by the mean value theorem, and the tangent turns at the rate $\kappa$ along the arc. Hence, for every point $\bxb_k$ of the arc,
\begin{equation}
  |\btaub_k\cdot\bn_c|\le2\kappa\ell,\qquad |\bx_k-\bxb_k|\le\kappa\ell^2 ,
  \label{eq:chord-geometry}
\end{equation}
the second bound being the distance between a chord and its arc. The same two quantities on the faces of the staircase, for the classical corrections with the geometry-normal map, only satisfy
\begin{equation}
  |\btaub_k\cdot\tilde\bn|\le1,\qquad |\bxt_k-\bxb_k|\le\sqrt2\,h .
  \label{eq:staircase-geometry}
\end{equation}
The first bound cannot be improved: the normal of a face is one of the Cartesian directions, while the wall takes all the directions, so that on the steps where the wall is almost parallel to the axes the face is almost orthogonal to it and $|\btaub_k\cdot\tilde\bn|$ is close to $1$, whatever the mesh size. The second one holds because the neighbour behind a computational face is not active, so that the wall crosses it, and the closest point of the wall is at most as far as the diagonal $\sqrt2\,h$ of this neighbour; with the face-orthogonal map there is no such bound, and the distances can reach several cells, as in \cref{fig:mapping}. Comparing \eqref{eq:chord-geometry} with \eqref{eq:staircase-geometry}, the angle between the boundary and the wall goes from $\mathcal{O}(1)$ to $\mathcal{O}(\kappa h)$, and the distance of the mapped points from $\mathcal{O}(h)$ to $\mathcal{O}(\kappa h^2)$, so that both quantities gain a factor $\kappa h$, which vanishes as the mesh is refined. For instance, in \cref{fig:gap}(c), the distance of a face point from the wall is $0.46\,h$, and that of a chord point $0.010\,h$. The construction gives a Cartesian mesh the boundary faces of the body-fitted meshes of~\cite{krivodonova2006high}, which are chords of the boundary, without adding any degree of freedom: the solution on each extended cell is simply the extrapolation of the polynomial of the original Cartesian cell, while the polynomial corrections remove the remaining distance.

The lateral edges play a second role. On the staircase, two boundary cells that meet at a concave corner share only a vertex, and they exchange no information: each of them sends the energy flux $T_1$ of its face into the region below the staircase, which belongs to no cell. In the gap method this region is part of the extended cells, the energy that crosses the faces of the staircase stays in the computational domain, and it reaches the neighbouring extended cell through the lateral edge, with a conservative and dissipative flux. The lateral edges are not optional: without them each quadrilateral would leave the uncompensated term $\int_l p_h\,\bu_h\cdot\bn_l$ on the edge, of size $\mathcal{O}(1)$ and without sign, as on the steps of the staircase. These two facts lead to the following results.

\begin{proposition}[gap method on a straight wall]\label{prop:gap-straight}
Let the wall be a straight line. Then the gap method with the \SB, SP or MP corrections is energy stable: $dE/dt\le0$.
\end{proposition}

\begin{proof}
The projections of the vertices lie on the line, so that every chord is contained in the wall. Every Gauss point of a chord is then its own projection, $\bxb_k=\bx_k$, and the normal of the chord is that of the wall, $\bn_c=\bar\bn$. Hence $\btaub_k\cdot\bn_c=0$, so that $T_1=0$. With $\bxb_k=\bx_k$ the SP coefficient \eqref{eq:sp-alpha} is $1$, and the MP correction imposes $\bar\bn\cdot\bv_h(\bx_k)=0$ while keeping the tangential velocity at the same points, so that $c_k=\bu_k\cdot\bar\bn$ for all the corrections and $e_k=0$, i.e.\ $T_2=0$. By \eqref{eq:energy-rod} the wall terms reduce to $-(\bu_k\cdot\bar\bn)^2\le0$. The interior faces, including the lateral edges, and the faces with exact data are dissipative, and \eqref{eq:energy-mesh} gives $dE/dt\le0$.
\end{proof}

On the same straight wall the classical corrections are unstable (\cref{tab:channel}), while the gap method dissipates energy. On a curved wall the chords are no longer contained in the wall, and the two terms $T_1$ and $T_2$ no longer vanish. They remain, however, controlled by the geometric bounds \eqref{eq:chord-geometry}. The term $T_1$ is proportional to the angle between the chord and the wall, which is at most $2\kappa\ell$, and the term $T_2$ is proportional to the remainder $e_k$, i.e.\ to the change of the normal velocity over the distance between the chord and the wall, which is at most $\kappa\ell^2$. Both are therefore proportional to the curvature, and they vanish for $\kappa=0$, where \cref{prop:gap-straight} applies. On the staircase the same two quantities, the angle between the face and the wall and the distance of the mapped points relative to the cell size, are $\mathcal{O}(1)$ whatever the mesh, and the channel of \cref{sec:wall-tilted} shows that the resulting energy production grows like $1/h$. The bounds \eqref{eq:chord-geometry} thus suggest that, on a curved wall, the energy that the gap method may produce at the wall is controlled by the local curvature through the product $\kappa h$, which can be reduced by a local refinement of the mesh. In the following section, we test the gap DG method on problems dominated by the slip wall boundary condition.

\subsection{Numerical results}
\label{sec:gap-results}

\subsubsection{Pressure pulse on a straight wall}
\label{sec:gap-pulse}

\begin{figure}
\centering
\includegraphics[width=0.95\linewidth]{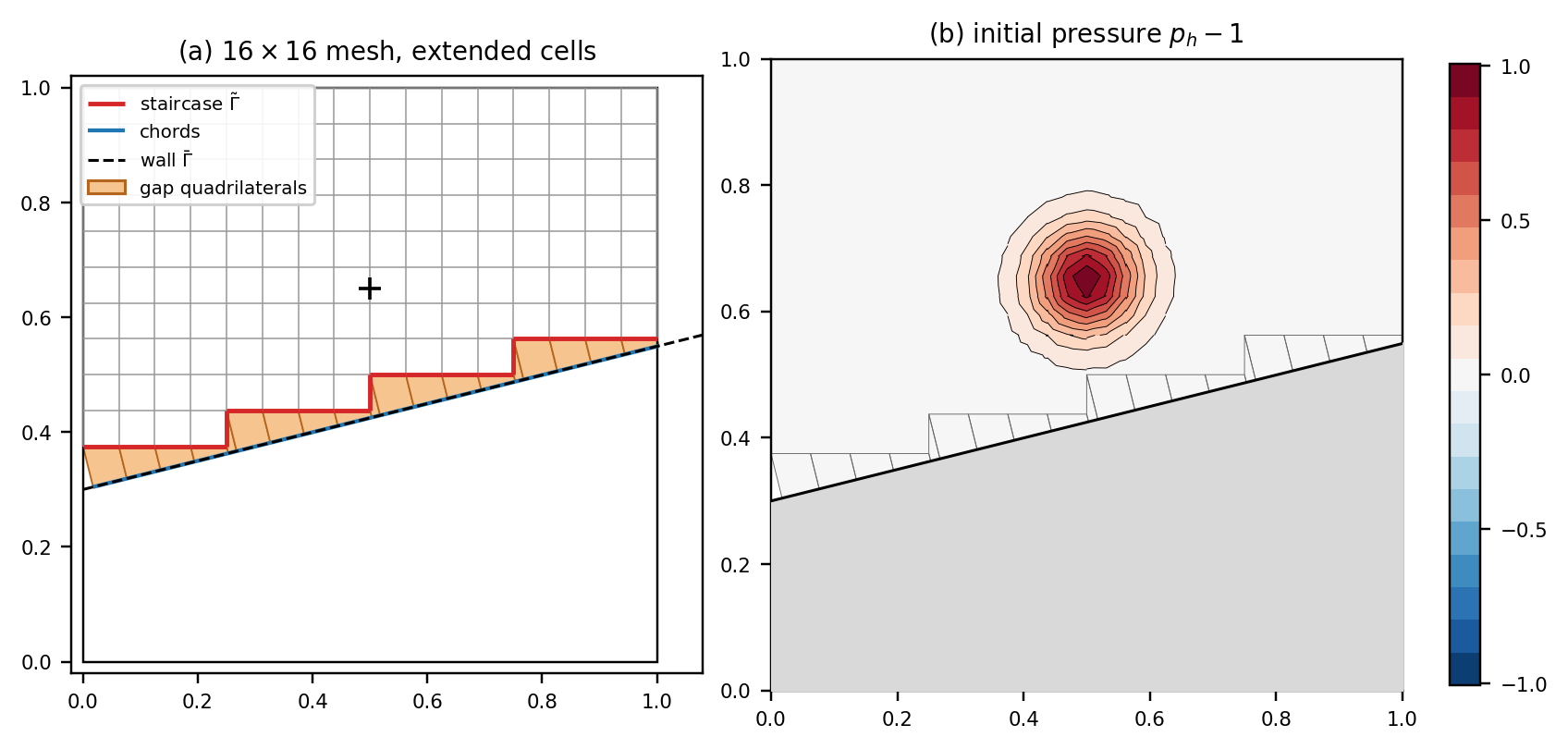}
\caption{Pressure pulse on a straight wall, $16\times16$ mesh: (a) active cells, staircase $\tilde\Gamma$, gap quadrilaterals and chords, which lie on the wall $\bar\Gamma$, and centre of the pulse ($+$); (b) initial pressure \eqref{eq:pulse} on the extended cells, the solid in grey.}
\label{fig:pulse-setup}
\end{figure}

We first verify \cref{prop:gap-straight} in the full solver. The wall is the straight line $y=0.3+x\tan\theta$, with $\theta=14^\circ$ as in the channel of \cref{sec:wall-tilted}, and the fluid is above it, in the box $[0,1]^2$ discretized with $16\times16$ cells, where the cells cut by the wall are discarded (\cref{fig:pulse-setup}(a)). The fluid is initially at rest, $\bu=\mathbf 0$, with the Gaussian pressure pulse
\begin{equation}
  p=1+e^{-|\bx-\bx_0|^2/\sigma^2},\qquad \bx_0=(0.5,\,0.65),\qquad \sigma=0.08,
  \label{eq:pulse}
\end{equation}
at the distance $0.22$ from the wall (\cref{fig:pulse-setup}(b)), and the state at rest is imposed on the sides of the box. To verify \cref{prop:gap-straight}, we use $\PP{2}$ and the time integration of order $3$ with $\text{CFL}=1$, and we compare gap-\RODLSP\ with the \RODLSP\ correction applied on the staircase, with the geometry-normal map. Since the system is linear and the state at rest is stationary, the energy $E$ of the perturbation $(p_h-1,\bu_h)$ satisfies the same balance as in \cref{sec:gap-energy}; we compute it with the mass matrices of the extended cells for the gap method, and of the Cartesian cells for the staircase.

\begin{figure}
\centering
\includegraphics[width=\linewidth]{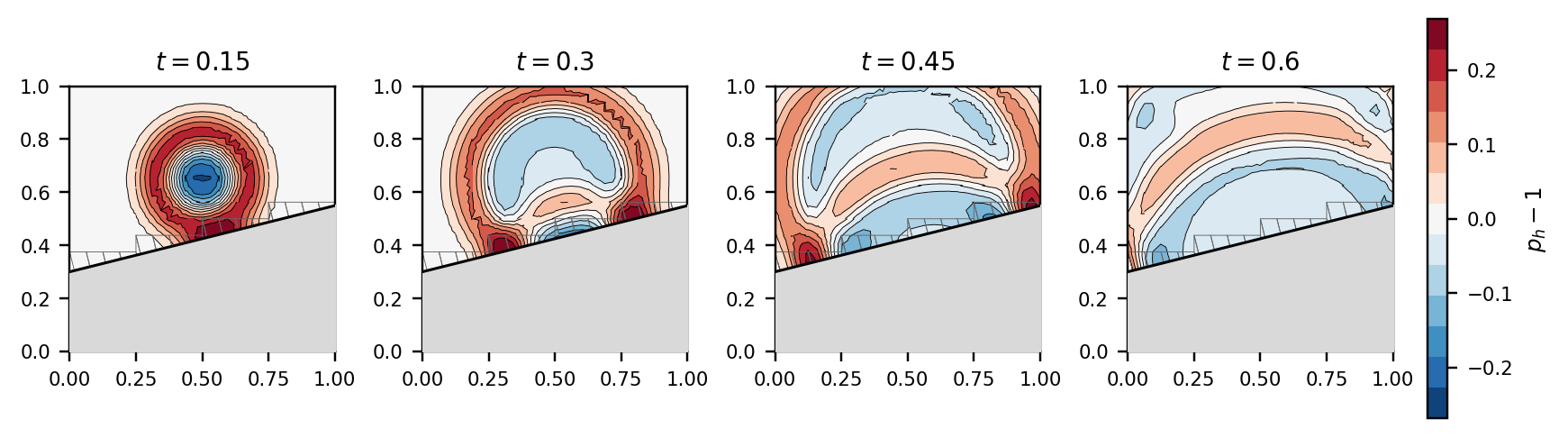}
\caption{Pressure pulse on a straight wall, gap-\RODLSP, $\PP{2}$, $16\times16$ cells: pressure perturbation $p_h-1$ at $t=0.15$, $0.3$, $0.45$ and $0.6$, on the active cells and on the gap quadrilaterals.}
\label{fig:pulse-transient}
\end{figure}

\begin{figure}
\centering
\includegraphics[width=0.95\linewidth]{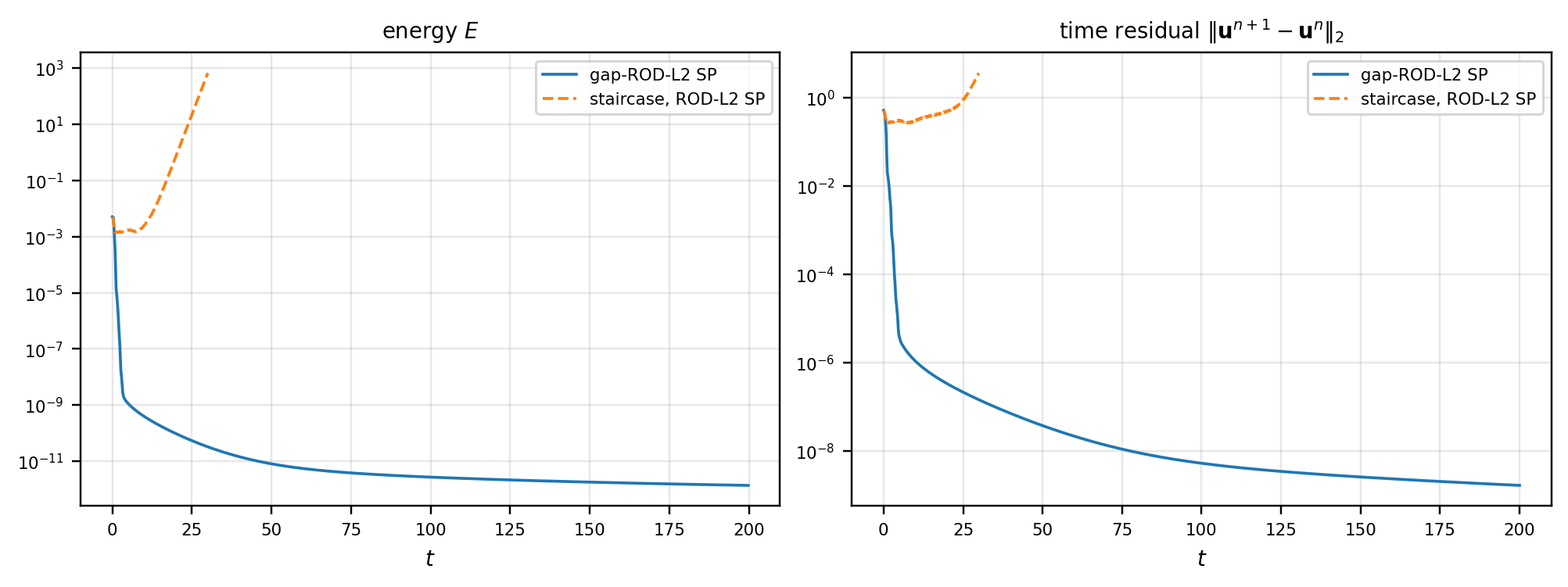}
\caption{Pressure pulse on a straight wall, $\PP{2}$, $16\times16$ cells: energy $E$ of the perturbation of the state at rest (left) and time residual (right), for gap-\RODLSP\ and for \RODLSP\ on the staircase.}
\label{fig:pulse-energy}
\end{figure}

The pulse reaches the wall at $t\approx0.1$, is reflected, and the waves leave the box through its sides (\cref{fig:pulse-transient}), with no visible trace of the staircase. With the gap method the energy never increases and the time residual follows the same trend (\cref{fig:pulse-energy}). With the correction on the staircase, instead, the energy grows exponentially once the pulse has left the box, by five orders of magnitude at $t=30$, as expected from the eigenvalues with positive real part of \cref{tab:channel}.

\subsubsection{Potential flow around a cylinder}
\label{sec:gap-cylinder}

\begin{figure}
\centering
\includegraphics[width=0.95\linewidth]{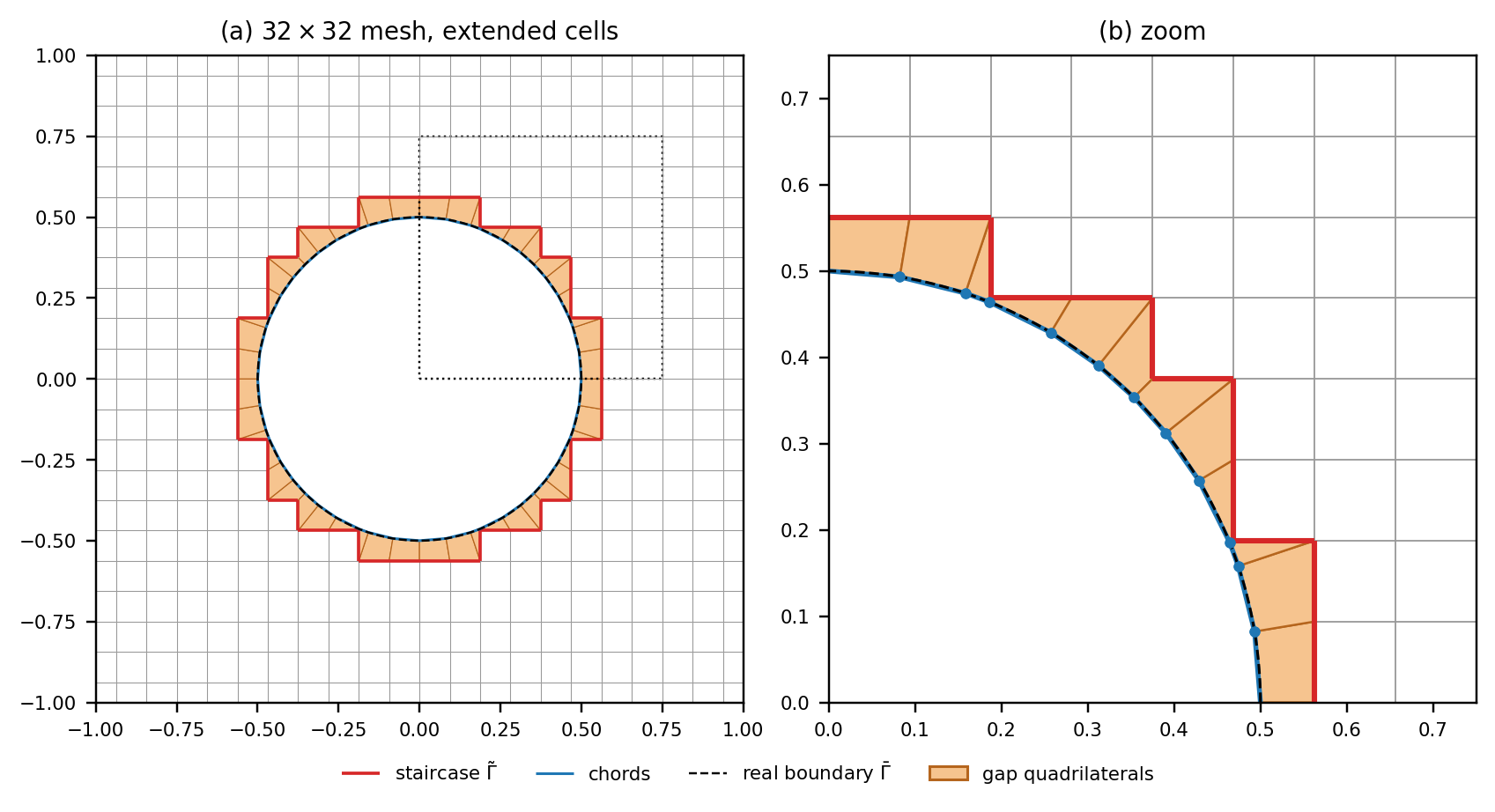}
\caption{Gap method around the cylinder, $32\times32$ mesh. (a) Active Cartesian cells, staircase $\tilde\Gamma$ (computational boundary of the classical corrections), gap quadrilaterals, chords and real boundary $\bar\Gamma$. (b) Zoom on the first quadrant (dotted square in (a)); the dots are the projections $P(A)$ of the vertices of the staircase, which are the ends of the chords.}
\label{fig:gap-mesh}
\end{figure}

We test the gap method in the full solver on the steady potential flow around a cylinder of radius $R=0.5$ centred at the origin, with speed $1$ at infinity,
\begin{equation}
  p=1,\qquad u^x=1-\frac{R^2(x^2-y^2)}{r^4},\qquad u^y=-\frac{2R^2xy}{r^4},\qquad r^2=x^2+y^2,
  \label{eq:potential-flow}
\end{equation}
which satisfies $\nabla p=0$ and $\nabla\cdot\bu=0$, and is therefore a stationary solution of \eqref{eq:acoustics} with $\bu\cdot\bar\bn=0$ on the cylinder. The fluid is outside the cylinder, in the box $[-1.5,1.5]^2$ discretized with $N\times N$ cells, the cells cut by the circle are discarded, the slip wall is imposed on the cylinder and the exact state on the sides of the box. The initial condition is the projection of \eqref{eq:potential-flow}, on the extended cells for the gap method, and the errors are measured in $L^2$ on the active cells at the discrete steady state. Since all the corrections behave in the same way within the gap method, we show the results of gap-\RODLSP.
\Cref{fig:gap-mesh} shows the computational domain of the gap method on the $32\times32$ mesh. The active cells are the same as in the classical approach, and their computational faces form the staircase $\tilde\Gamma$. Each computational face is extended by its gap quadrilateral up to the cylinder, so that the boundary condition is imposed on the chords, a polygon inscribed in the circle with vertices on $\bar\Gamma$. Depending on the position of the cell, a boundary cell is extended by one quadrilateral, when only one of its faces is a computational face, or by two, at the corners of the staircase; the quadrilaterals of neighbouring cells share their lateral edges, which carry the internal fluxes between extended cells.

The purpose of this test is to check that the gap method is stable, reaches the correct steady state, and converges to it with the expected order. This needs some care with the linear acoustic system. Every divergence-free velocity field with constant pressure is a stationary solution of \eqref{eq:acoustics}, and the vorticity $\omega=\partial_xu^y-\partial_yu^x$ satisfies $\partial_t\omega=0$: unlike the acoustic waves, which leave the box through its sides, the vorticity is not transported, and the numerical vorticity produced at the wall, where the terms $T_1$ and $T_2$ of the chords do not vanish on a curved wall (\cref{sec:gap-energy}), stays in the domain. The steady state reached by a time-marching procedure is therefore not determined by the boundary data alone, but retains the vortical content of the initial data and the vorticity produced by the discretization~\cite{jung2022steady}. The potential flow \eqref{eq:potential-flow} is the steady solution selected by the involution $\nabla\times\bu=\mathbf 0$, which is preserved by the continuous dynamics, but not exactly by the discrete one in the presence of an embedded wall. In long computations, the time residual then stagnates instead of decreasing to machine precision, at a level that decreases as the mesh is refined. This is a feature of the acoustic system, in which the vortical modes have zero speed: in preliminary computations, not reported here, in which a transport term by a mean flow $\mathbf U_0$ is added to \eqref{eq:acoustics}, i.e.\ $\partial_tp+\mathbf U_0\cdot\nabla p+\nabla\cdot\bu=0$ and $\partial_t\bu+(\mathbf U_0\cdot\nabla)\bu+\nabla p=\mathbf 0$, with $\mathbf U_0$ the potential flow \eqref{eq:potential-flow} scaled by $0.3$, the vortical modes are advected by the mean flow and leave the domain, and the time residual reaches machine precision much faster. To reach the steady state, in this test we add a penalty on the vorticity to the velocity equation of \eqref{eq:acoustics}, which becomes
\begin{equation}
  \partial_t\bu+\nabla p+\sigma\,\nabla^\perp\omega=\mathbf 0,
  \qquad
  \omega=\partial_xu^y-\partial_yu^x,\qquad
  \nabla^\perp\omega=\begin{pmatrix}\partial_y\omega\\-\partial_x\omega\end{pmatrix},\qquad
  \sigma=\beta\,\frac{h}{(p+1)^2},
  \label{eq:curl-penalty}
\end{equation}
where $h$ is the mesh size and $p$ the polynomial degree, while the pressure equation is unchanged. Since $\nabla^\perp\omega=\nabla\times(\nabla\times\bu)=\nabla(\nabla\cdot\bu)-\Delta\bu$, the penalty is a viscosity that acts only on the rotational part of the velocity: it vanishes on the irrotational solutions, such as \eqref{eq:potential-flow}, so that the scheme remains consistent. In particular, since the exact solution is irrotational, the penalty introduces no consistency error on the steady state, and it only accelerates the decay of the vortical part of the discrete solution, which the scheme would otherwise dissipate very slowly. In the DG discretization the penalty is integrated by parts on each cell and only the volume part is kept, with no contribution on the faces: the equations of $u^x$ and $u^y$ receive on the left-hand side of \eqref{eq:dg-weak} the terms
\begin{equation}
  -\sigma\int_{K}\omega_h\,\partial_y\varphi_m\,d\bx
  \qquad\text{and}\qquad
  \sigma\int_{K}\omega_h\,\partial_x\varphi_m\,d\bx,
  \qquad
  \omega_h=\partial_xu^y_h-\partial_yu^x_h,
  \label{eq:curl-penalty-dg}
\end{equation}
with $\beta=0.36$ for all the degrees and meshes, and the integrals are computed on the extended cells for the gap method. Taking $u^x_h$ and $u^y_h$ as test functions, the two terms add $-\sigma\sum_K\int_{\Omega_K}\omega_h^2\,d\bx$ to the energy balance \eqref{eq:energy-mesh}, so that the penalty dissipates only the vorticity. The penalty is used here only to reach the steady state: it does not modify the boundary terms of the energy balance at the wall (\cref{sec:gap-energy}), on which the stability of the gap method relies.

The computations are run in serial, without any parallelism, and with the explicit schemes the time step is limited by the CFL condition: with $\PP{3}$ on the $32\times32$ mesh, about $9\times10^4$ time steps are needed to reach the steady state, and on the $64\times64$ and $128\times128$ meshes with $\PP{2}$ and $\PP{3}$, where the time step is smaller and the number of cells larger, reaching the steady state with the explicit schemes would take a very long time. For this reason, the time integration is performed with the implicit Euler scheme, with a pseudo-time step that grows as the steady residual $\|\mathbf L(\mathbf U_h)\|_2$ decreases~\cite{mulder1985experiments}, where $\mathbf L$ is the semi-discrete operator, and each step is solved with GMRES~\cite{saad1986gmres} preconditioned by a block symmetric Gauss--Seidel method. To make sure that the implicit scheme does not hide an instability, the computations on the $16\times16$ and $32\times32$ meshes have also been advanced with the explicit schemes of \cref{sec:dg} until the time residual $\|\bu^{n+1}-\bu^n\|_2$ reaches $10^{-13}$, with $\text{CFL}=1$ for $\PP{2}$ and $\PP{3}$ and $\text{CFL}=0.5$ for $\PP{1}$: they reach the same errors as the implicit scheme, in all the digits reported in \cref{tab:gap-conv}.

\begin{table}
\centering
\caption{Potential flow around a cylinder with the slip wall, gap-\RODLSP: $L^2$ errors of $p$ and $u^x$ at the steady state, with the vorticity penalty \eqref{eq:curl-penalty}, and orders of accuracy.}
\label{tab:gap-conv}
\small
\setlength{\tabcolsep}{4pt}
\begin{tabular}{rcccccccccccc}
\toprule
 & \multicolumn{4}{c}{$\PP{1}$} & \multicolumn{4}{c}{$\PP{2}$} & \multicolumn{4}{c}{$\PP{3}$}\\
\cmidrule(lr){2-5}\cmidrule(lr){6-9}\cmidrule(lr){10-13}
$N$ & $p$ & order & $u^x$ & order & $p$ & order & $u^x$ & order & $p$ & order & $u^x$ & order\\
\midrule
16  & 2.54E-2 & --   & 8.62E-2 & --   & 2.48E-3 & --   & 4.77E-3 & --   & 3.01E-4 & --   & 7.18E-4 & --\\
32  & 6.60E-3 & 1.94 & 2.86E-2 & 1.59 & 4.36E-4 & 2.51 & 8.96E-4 & 2.41 & 3.28E-5 & 3.20 & 8.34E-5 & 3.10\\
64  & 1.48E-3 & 2.15 & 7.18E-3 & 2.00 & 5.09E-5 & 3.10 & 9.43E-5 & 3.25 & 1.92E-6 & 4.10 & 4.62E-6 & 4.17\\
128 & 3.65E-4 & 2.03 & 1.48E-3 & 2.28 & 6.42E-6 & 2.99 & 1.11E-5 & 3.09 & 1.29E-7 & 3.89 & 3.25E-7 & 3.83\\
\bottomrule
\end{tabular}
\end{table}

\begin{figure}
\centering
\includegraphics[width=0.74\linewidth]{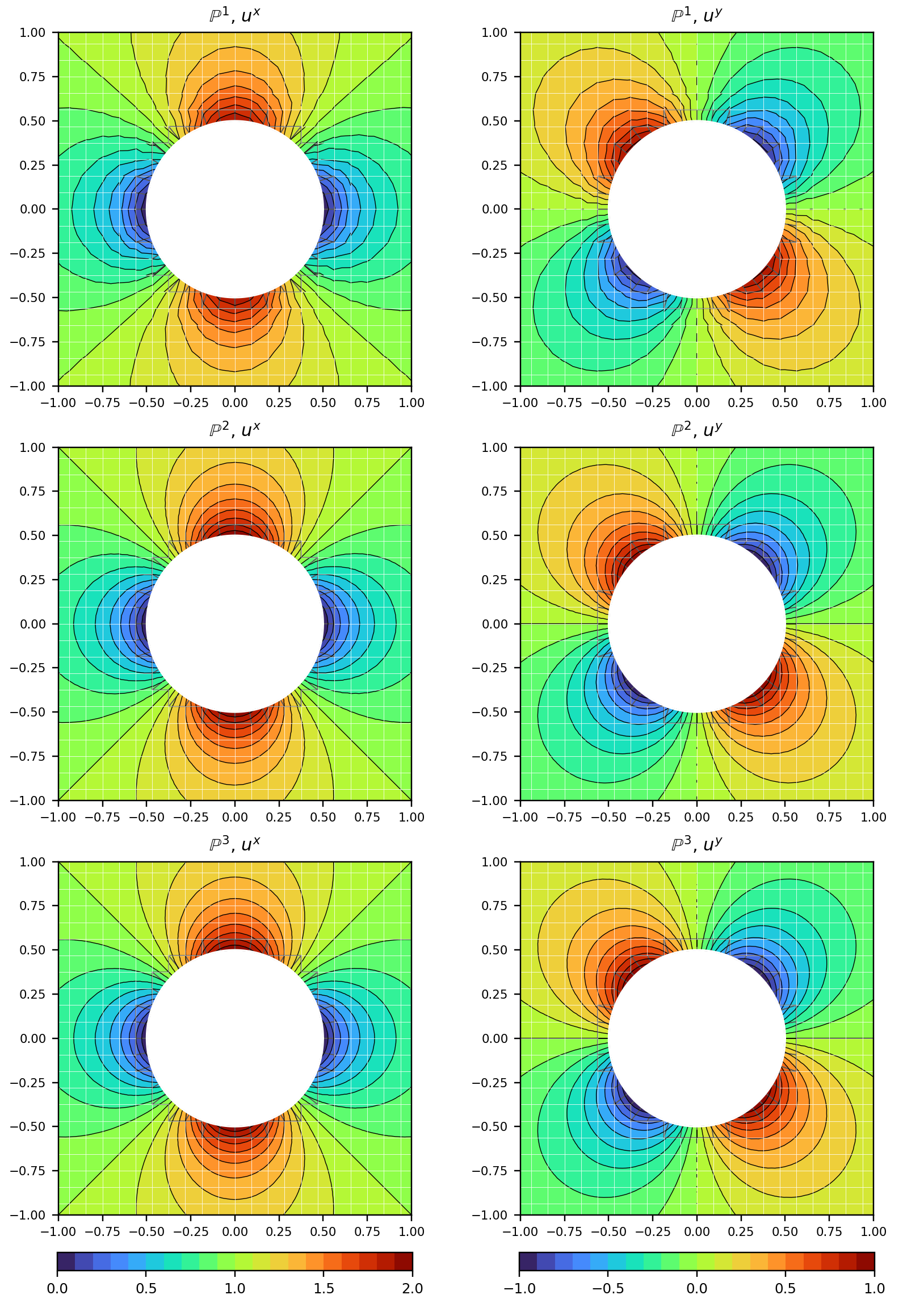}
\caption{Potential flow around a cylinder with the slip wall, gap-\RODLSP, $32\times32$ mesh, steady state: $u^x$ (left) and $u^y$ (right) with $\PP{1}$, $\PP{2}$ and $\PP{3}$, on the active cells and on the gap quadrilaterals, with isolines every $0.1$. White lines: Cartesian cells; grey lines: gap quadrilaterals.}
\label{fig:gap-solution}
\end{figure}

\Cref{tab:gap-conv} reports the convergence analysis on a sequence of nested meshes: with the gap method all the computations reach the steady state, and the errors converge with the expected orders. In the same simulations, the same correction on the staircase diverges. The energy that can be produced at the wall, which is $\mathcal{O}(1)$ on the staircase, is reduced by the chords to terms proportional to $\kappa h$ (\cref{sec:gap-energy}), so that it is controlled by how well the mesh resolves the curvature of the wall, and decreases as the mesh is refined; the small vortical part that it produces is removed by the penalty. \Cref{fig:gap-solution} shows the velocity on the $32\times32$ mesh, where each boundary cell is extended by one or two quadrilaterals up to the cylinder. With such a coarse mesh and $\PP{1}$ the maximum of $u^x$ is underestimated, while with $\PP{2}$ and $\PP{3}$ the solution is much smoother close to the wall with the expected values, with no visible trace of the staircase.

\subsubsection{Start-up of the flow around an inclined ellipse}
\label{sec:gap-ellipse}

\begin{figure}
\centering
\includegraphics[width=0.9\linewidth]{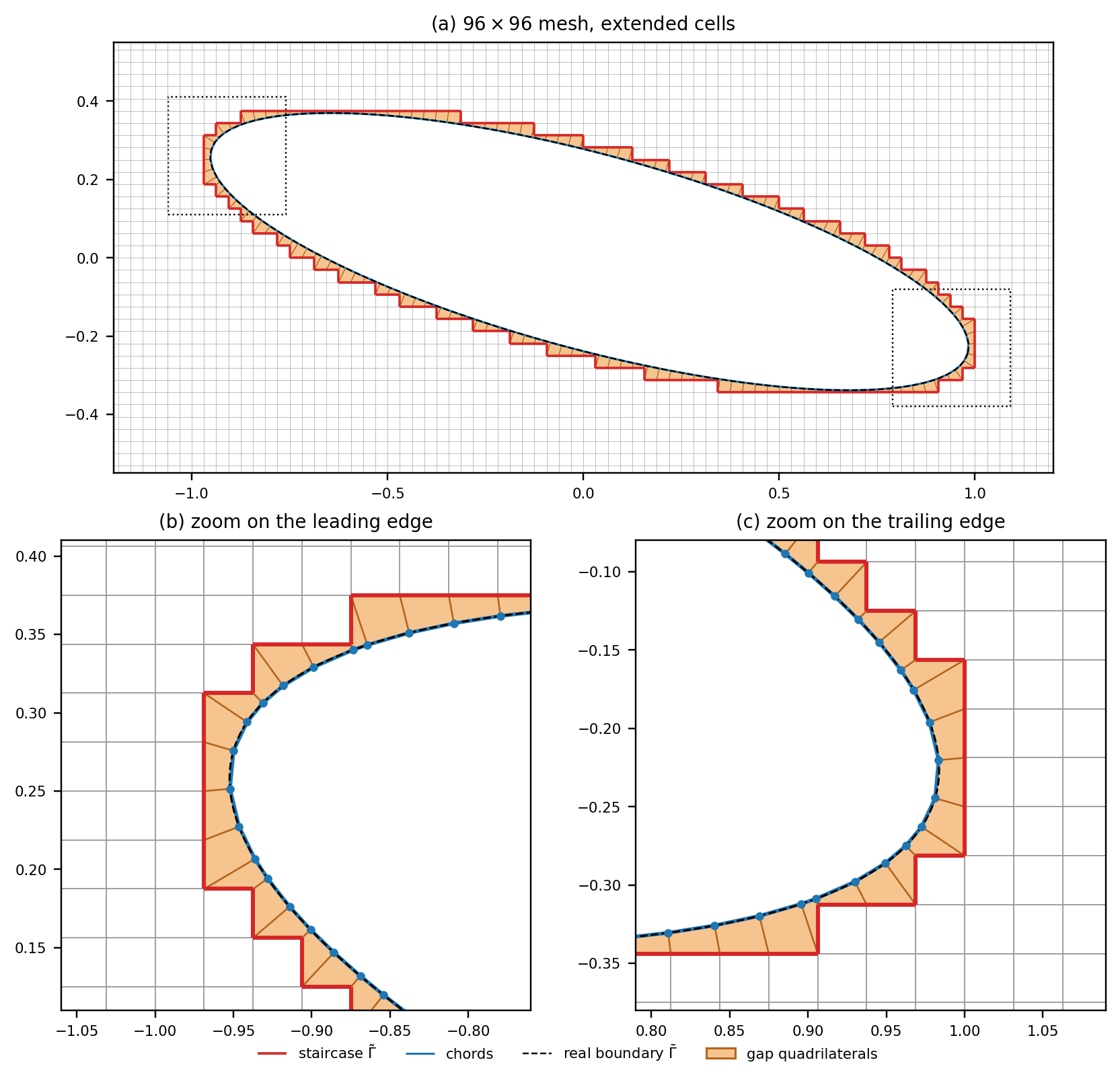}
\caption{Gap method around the inclined ellipse, $96\times96$ mesh. (a) Active Cartesian cells, staircase $\tilde\Gamma$, gap quadrilaterals, chords and real boundary $\bar\Gamma$. (b) Zoom on the leading edge and (c) on the trailing edge (dotted squares in (a)); the dots are the ends of the chords.}
\label{fig:ellipse-mesh}
\end{figure}

In the previous test the initial condition was the exact solution, so that the wall only had to keep the flow in place. We conclude with a computation in which the flow has to be built by the waves that reach the wall. The body is an ellipse with semi-axes $a=1$ and $b=0.25$, centred at the origin and inclined by $15^\circ$, with the leading edge up, so that a flow from the left meets it with an incidence of $15^\circ$, as a thick symmetric profile. Its curvature ranges from $b/a^2=0.25$ on the sides to $a/b^2=16$ at the two ends. The domain is the box $[-1.5,1.5]^2$ with $96\times96$ cells, so that the radius of curvature at the ends, $0.0625$, is twice the mesh size, the cells cut by the ellipse are discarded, the slip wall is imposed with gap-\RODLSP\ and $\PP{3}$, with the vorticity penalty \eqref{eq:curl-penalty}. The transient is computed with the explicit scheme and $\text{CFL}=1$, and the steady state with the implicit scheme of the previous test.
The extended cells around the ellipse are shown in \cref{fig:ellipse-mesh}. At both ends the staircase has a vertical step of four cells, at about half a cell from the end of the ellipse. This is the region where the chords are farthest from the wall, as discussed at the end of this section.

The fluid is initially at rest, $p=1$ and $\bu=\mathbf 0$, and on the sides of the box we impose the exact stationary solution, i.e.\ the potential flow around the ellipse with speed $1$ along $x$ at infinity and no circulation, which is also the limit of the computation for large times. It reads
\begin{equation}
  p=1,\qquad u^x=\operatorname{Re}W,\qquad u^y=-\operatorname{Im}W,\qquad
  W=e^{-i\theta}\,\frac{e^{i\theta}-R^2e^{-i\theta}/\zeta^2}{1-m^2/\zeta^2},
  \label{eq:ellipse-flow}
\end{equation}
where $z=x+iy$, $z'=(z-z_c)e^{-i\theta}$ is the coordinate in the frame of the ellipse, with $z_c$ its centre and $\theta=-15^\circ$ its rotation, $R=(a+b)/2$, $m^2=(a^2-b^2)/4$, and $\zeta$ is the root of $\zeta^2-z'\zeta+m^2=0$ with $|\zeta|\ge m$ (Joukowski map). This stationary state is reached only through the interaction of the waves with the wall.

\begin{figure}
\centering
\includegraphics[width=0.8\linewidth]{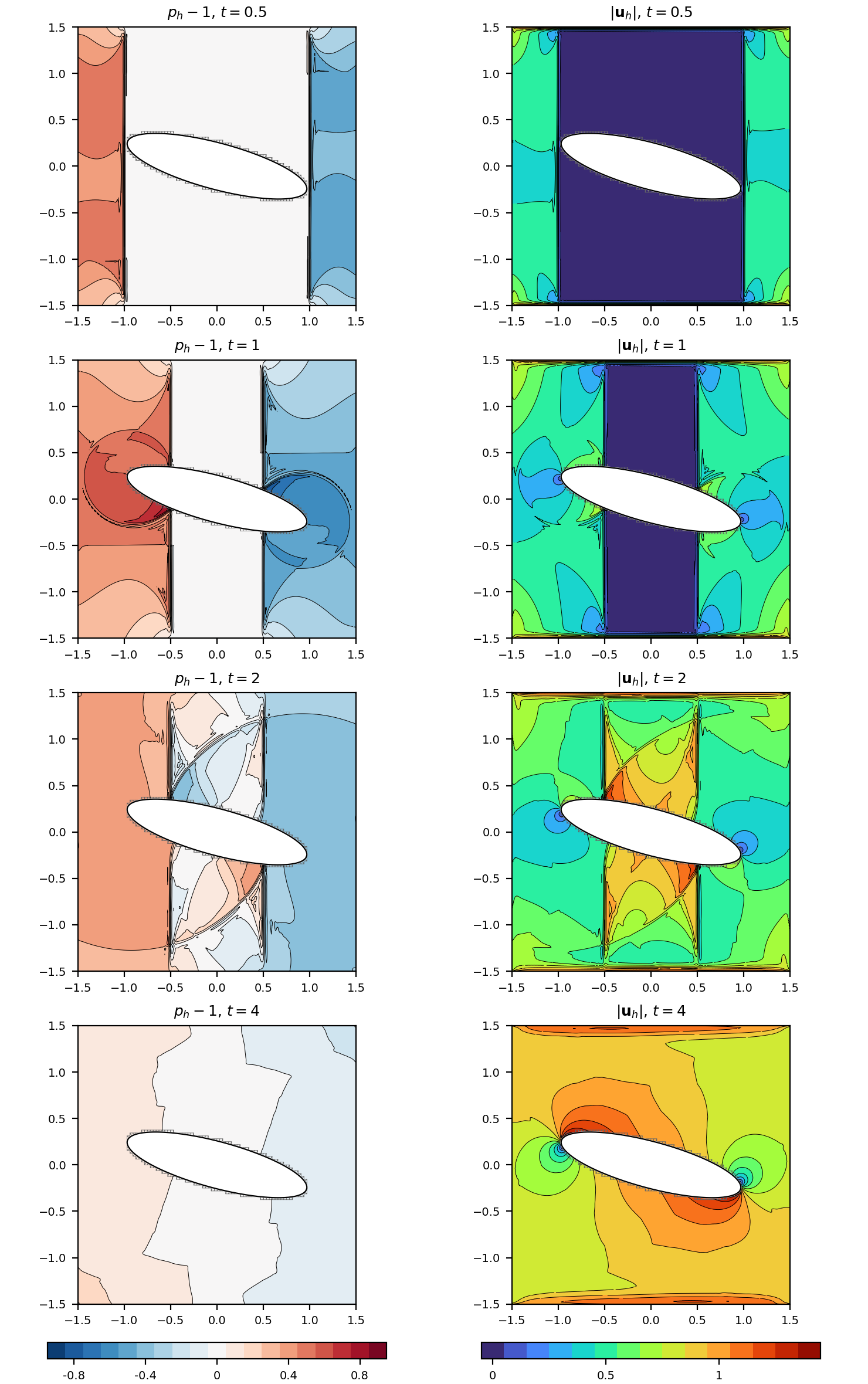}
\caption{Start-up of the flow around the inclined ellipse, gap-\RODLSP, $\PP{3}$, $96\times96$ cells: pressure perturbation $p_h-1$ (left) and speed $|\bu_h|$ (right) of the numerical solution at $t=0.5$, $1$, $2$ and $4$, with isolines at the boundaries of the colour bands; the levels are shifted by half a band, so that the state at rest and the far field $|\bu|=1$ lie inside a band.}
\label{fig:ellipse-startup}
\end{figure}

The transient is shown in \cref{fig:ellipse-startup}. The exact data imposed on the sides of the box at $t=0$ differ from the state at rest, and two plane waves enter the box, a compression from the left, where the flow enters, and an expansion from the right. They reach the two ends of the ellipse at $t\approx0.5$, and at $t=1$ the waves reflected by the ends are visible as circular fronts centred on them, while the plane fronts move along the body. At $t=2$ the plane fronts have crossed the ellipse, and the reflected fronts reach the sides of the box, where they leave the domain through the exact condition. The flow then settles around the ellipse: at $t=4$ the pressure is already close to $1$, and the velocity shows the structure of the potential flow without circulation, with the two stagnation points near the ends, where $|\bu|$ vanishes, and the regions where the flow is accelerated, $|\bu|>1$, on the upper side near the leading edge and on the lower side near the trailing edge.

\begin{figure}
\centering
\includegraphics[width=0.9\linewidth]{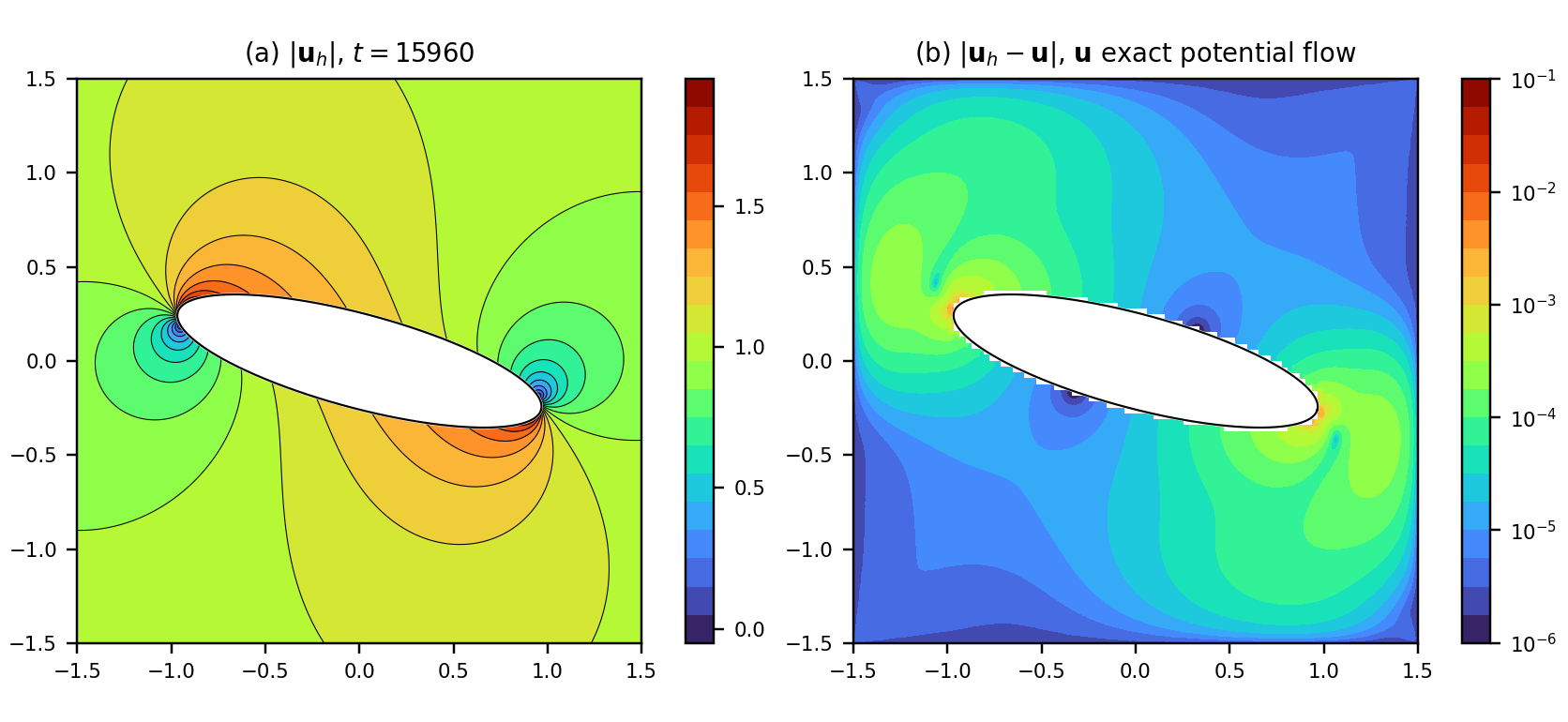}
\caption{Inclined ellipse at the steady state, gap-\RODLSP, $\PP{3}$, $96\times96$ cells: (a) speed of the computed flow, with isolines every $0.1$; (b) pointwise error of the velocity with respect to the exact potential flow \eqref{eq:ellipse-flow}, on the active cells, in logarithmic scale.}
\label{fig:ellipse-final}
\end{figure}

The steady state is computed with the implicit scheme, starting from the same state at rest: in $58$ pseudo-time steps the steady residual decreases by almost twelve orders of magnitude, down to the round-off level, and the $L^2$ errors with respect to \eqref{eq:ellipse-flow} are $3.85\times10^{-5}$ for $p$, $1.77\times10^{-4}$ for $u^x$ and $2.08\times10^{-4}$ for $u^y$. The steady solution is compared with the potential flow \eqref{eq:ellipse-flow} in \cref{fig:ellipse-final}. 
The solution appears smooth and free of spurious waves. The pointwise error of the velocity is shown in logarithmic scale, and it is concentrated in the regions of high curvature at the two ends of the ellipse.

\section{Conclusions}
\label{sec:conclusions}

This paper presented a stability study of high order discontinuous Galerkin methods coupled with polynomial corrections for embedded boundaries on Cartesian meshes, in the context of hyperbolic systems, as a continuation of the 1D analyses of~\cite{ciallella2027stability,ciallella2026minimization}. Combining the spectra of the discrete operators with the energy balance of the boundary flux, we have shown that the stability results obtained in 1D do not carry over to 2D in general. They extend to 2D only when the distance between the surrogate and the true boundary is constant along the boundary, in which case the 2D operator reduces to the 1D one. For the linear advection equation, whose operator is block triangular so that the stability of a mesh is that of its boundary cells, a single straight boundary tilted with respect to the mesh is enough to make all the corrections unstable with the geometry-normal map, depending on the position of the boundary and on the direction of the flow. With the face-orthogonal map, instead, ROD-E remains stable up to $\PP{2}$ and ROD-$L^2$ up to $\PP{3}$, one degree less than in 1D, as also confirmed by the numerical experiments on a disk, likely because with tensor-product bases this map gives a one-dimensional character to the boundary condition.

The main findings show that the impact of the embedded treatment is even larger for the slip wall of the linear acoustic system, whose operator is not block triangular and couples every cell to its neighbours. All the corrections, and even the boundary treatment that only uses the normal of the true wall, are unstable on a tilted straight wall. The energy balance of the surrogate wall flux identifies the origin of this behaviour in the geometry of the surrogate boundary rather than in the corrections: a normal correction term, related to the angle between the faces and the wall, as in the boundary treatment of~\cite{krivodonova2006high}, and a location correction term, related to the distance of the wall from the faces relative to the cell size, are both $\mathcal{O}(1)$, i.e.\ they do not vanish as the mesh is refined, on a surrogate boundary formed by the faces of a Cartesian mesh, where the angle is $\mathcal{O}(1)$ and the distance $\mathcal{O}(h)$.

The size of these terms is therefore set by the geometry of the surrogate boundary, through its angle with the wall and its distance from it, and not by the corrections themselves. This suggests acting on the surrogate boundary itself, by extending the boundary cells towards the wall. To this end, we developed for the first time in a high order DG discretization of hyperbolic systems the gap paradigm, which was originally introduced in the shifted boundary method for continuous Galerkin discretizations of elliptic problems to recover the optimal accuracy with Neumann boundary conditions~\cite{collins2026gap}, and later extended to isogeometric discretizations and solid mechanics~\cite{gorgi2026isogeometric,li2026gap}. The boundary cells are extended to the polygon inscribed in the wall, the condition is imposed on its chords, and the gap cells of neighbouring boundary cells exchange information through the same numerical fluxes used on the interior interfaces, which makes the DG framework a natural setting for this construction. On the chords, the angle with the wall and the distance from it are $\mathcal{O}(\kappa h)$ and $\mathcal{O}(\kappa h^2)$, instead of $\mathcal{O}(1)$ and $\mathcal{O}(h)$ on the staircase, so that both terms gain a factor $\kappa h$. The method is energy stable on straight walls, and on curved walls the terms that can produce energy are reduced to terms proportional to the local curvature. In the numerical experiments, we verify the energy stability on a straight wall, and show the high accuracy of the method on curved boundaries. In all the tests, the corrections applied on the staircase feature unstable modes that grow very fast. The price to pay for this approach is the construction of the gap cells and the computation of the integrals on them, which make the method more expensive than the classical corrections on the surrogate Cartesian faces, which are however unstable for this problem.

Future work will address the remaining limitations of the method. Since the energy production terms are proportional to the product $\kappa h$, they still act where the curvature is poorly resolved by the mesh, and a local refinement of the Cartesian mesh near the embedded boundary, through quadtree meshes, is the natural way to keep $\kappa h$ small where it matters without refining the whole domain. The development of a fully energy stable high order polynomial correction approach for embedded curved boundaries, with an energy estimate that holds independently of the curvature, remains an open question. The extension of the gap method to the nonlinear system of the Euler equations for compressible gas dynamics and to moving boundaries will also be investigated.

\section*{Acknowledgements}
This work has received support under the program ``Investissement d'Avenir'' launched by the French Government and implemented by ANR, with the reference ANR-18-IdEx-0001 as part of its program Emergence. This research was funded, in whole or in part, by the French National Research Agency (ANR) under the project ANR-26-CE46-7233-01. 
The AI assistant Claude (Anthropic) was used during the preparation of this work to speed up the writing and the implementation; the author checked all the content produced in this way, verified the numerical results, and takes full responsibility for the paper.

\FloatBarrier
\bibliographystyle{abbrvnat}
\bibliography{paper_final}

\begin{thebibliography}{57}
\providecommand{\natexlab}[1]{#1}
\providecommand{\url}[1]{\texttt{#1}}
\expandafter\ifx\csname urlstyle\endcsname\relax
  \providecommand{\doi}[1]{doi: #1}\else
  \providecommand{\doi}{doi: \begingroup \urlstyle{rm}\Url}\fi

\bibitem[Adjerid and Moon(2019)]{adjerid2019immersed}
S.~Adjerid and K.~Moon.
\newblock An immersed discontinuous {G}alerkin method for acoustic wave
  propagation in inhomogeneous media.
\newblock \emph{SIAM Journal on Scientific Computing}, 41\penalty0
  (1):\penalty0 A139--A162, 2019.

\bibitem[Antonelli et~al.(2026)Antonelli, Gorgi, Zorrilla, and
  Rossi]{antonelli2026isogeometric}
N.~Antonelli, A.~Gorgi, R.~Zorrilla, and R.~Rossi.
\newblock Isogeometric multipatch coupling with arbitrary refinement and
  parametrization using the {G}ap--shifted boundary method.
\newblock \emph{Computer Methods in Applied Mechanics and Engineering},
  456:\penalty0 118913, 2026.

\bibitem[Atallah et~al.(2022)Atallah, Canuto, and Scovazzi]{atallah2022high}
N.~M. Atallah, C.~Canuto, and G.~Scovazzi.
\newblock The high-order shifted boundary method and its analysis.
\newblock \emph{Computer Methods in Applied Mechanics and Engineering},
  394:\penalty0 114885, 2022.

\bibitem[Bassi and Rebay(1997)]{bassi1997high}
F.~Bassi and S.~Rebay.
\newblock {High-order accurate discontinuous finite element solution of the
  2{D} {E}uler equations}.
\newblock \emph{Journal of Computational Physics}, 138\penalty0 (2):\penalty0
  251--285, 1997.

\bibitem[Berger(2015)]{BERGER2015180}
M.~Berger.
\newblock A note on the stability of cut cells and cell merging.
\newblock \emph{Applied Numerical Mathematics}, 96:\penalty0 180--186, 2015.

\bibitem[Berger(2017)]{berger2017cut}
M.~Berger.
\newblock Cut cells: Meshes and solvers.
\newblock In \emph{Handbook of Numerical Analysis}, volume~18, pages 1--22.
  Elsevier, Amsterdam, 2017.

\bibitem[Berger and Aftosmis(2012)]{berger2012progress}
M.~Berger and M.~Aftosmis.
\newblock Progress towards a {C}artesian cut-cell method for viscous
  compressible flow.
\newblock In \emph{50th AIAA Aerospace Sciences Meeting Including the New
  Horizons Forum and Aerospace Exposition}, page 1301, 2012.

\bibitem[Boscheri and Ciallella(2025)]{BOSCHERI2025114215}
W.~Boscheri and M.~Ciallella.
\newblock {High order treatment of moving curved boundaries:
  Arbitrary-{L}agrangian-{E}ulerian methods with a shifted boundary polynomial
  correction}.
\newblock \emph{Journal of Computational Physics}, 539:\penalty0 114215, 2025.

\bibitem[Burman(2010)]{burman2010ghost}
E.~Burman.
\newblock Ghost penalty.
\newblock \emph{Comptes Rendus. Math{\'e}matique}, 348\penalty0
  (21-22):\penalty0 1217--1220, 2010.

\bibitem[Burman et~al.(2015)Burman, Claus, Hansbo, Larson, and
  Massing]{burman2015cutfem}
E.~Burman, S.~Claus, P.~Hansbo, M.~G. Larson, and A.~Massing.
\newblock {CutFEM}: discretizing geometry and partial differential equations.
\newblock \emph{International Journal for Numerical Methods in Engineering},
  104\penalty0 (7):\penalty0 472--501, 2015.

\bibitem[Burman et~al.(2025)Burman, Hansbo, Larson, and Zahedi]{burman2025cut}
E.~Burman, P.~Hansbo, M.~G. Larson, and S.~Zahedi.
\newblock Cut finite element methods.
\newblock \emph{Acta Numerica}, 34:\penalty0 1--121, 2025.

\bibitem[Ciallella(2026)]{ciallella2026minimization}
M.~Ciallella.
\newblock Minimization-based embedded boundary methods as polynomial
  corrections: a stability study of discontinuous {G}alerkin for hyperbolic
  equations.
\newblock \emph{Journal of Scientific Computing}, 109\penalty0 (1):\penalty0
  28, 2026.

\bibitem[Ciallella(2027)]{ciallella2027stability}
M.~Ciallella.
\newblock Stability analysis of discontinuous {G}alerkin with a high order
  embedded boundary treatment for hyperbolic equations.
\newblock \emph{Applied Mathematics and Computation}, 534:\penalty0 130309,
  2027.

\bibitem[Ciallella and Boscheri(2026)]{ciallellaboscheri2026minimization}
M.~Ciallella and W.~Boscheri.
\newblock Minimization-based polynomial corrections for high-order curved
  boundaries on fixed and moving domains: assessment on finite volume and
  discontinuous galerkin schemes.
\newblock \emph{arXiv preprint arXiv:2608.05768}, 2026.

\bibitem[Ciallella and Visbech(2026)]{ciallella2026generalized}
M.~Ciallella and J.~Visbech.
\newblock Generalized high-order minimization-based polynomial corrections on
  unfitted spectral elements for the {P}oisson problem.
\newblock \emph{arXiv preprint arXiv:2608.06999}, 2026.

\bibitem[Ciallella et~al.(2023)Ciallella, Gaburro, Lorini, and
  Ricchiuto]{ciallella2023shifted}
M.~Ciallella, E.~Gaburro, M.~Lorini, and M.~Ricchiuto.
\newblock Shifted boundary polynomial corrections for compressible flows: high
  order on curved domains using linear meshes.
\newblock \emph{Applied Mathematics and Computation}, 441:\penalty0 127698,
  2023.

\bibitem[Ciallella et~al.(2024)Ciallella, Clain, Gaburro, and
  Ricchiuto]{ciallella2024very}
M.~Ciallella, S.~Clain, E.~Gaburro, and M.~Ricchiuto.
\newblock Very high order treatment of embedded curved boundaries in
  compressible flows: {ADER} discontinuous {G}alerkin with a space-time
  reconstruction for off-site data.
\newblock \emph{Computers \& Mathematics with Applications}, 175:\penalty0
  1--18, 2024.

\bibitem[Cockburn and Shu(2001)]{cockburn2001runge}
B.~Cockburn and C.-W. Shu.
\newblock {Runge--Kutta discontinuous {G}alerkin methods for
  convection-dominated problems}.
\newblock \emph{Journal of Scientific Computing}, 16\penalty0 (3):\penalty0
  173--261, 2001.

\bibitem[Collins et~al.(2023)Collins, Lozinski, and
  Scovazzi]{collins2023penalty}
J.~H. Collins, A.~Lozinski, and G.~Scovazzi.
\newblock A penalty-free shifted boundary method of arbitrary order.
\newblock \emph{Computer Methods in Applied Mechanics and Engineering},
  417:\penalty0 116301, 2023.

\bibitem[Collins et~al.(2026)Collins, Li, Lozinski, and
  Scovazzi]{collins2026gap}
J.~H. Collins, K.~Li, A.~Lozinski, and G.~Scovazzi.
\newblock {G}ap-{SBM}: A new conceptualization of the shifted boundary method
  with optimal convergence for the {N}eumann and {D}irichlet problems.
\newblock \emph{Computer Methods in Applied Mechanics and Engineering},
  453:\penalty0 118793, 2026.

\bibitem[Coppeans et~al.(2025)Coppeans, Fidkowski, and
  Martins]{coppeans2025aerodynamic}
A.~Coppeans, K.~Fidkowski, and J.~R. Martins.
\newblock Aerodynamic shape optimization with curved mesh adaptation.
\newblock In \emph{AIAA SCITECH 2025 Forum}, page 0782, 2025.

\bibitem[Costa et~al.(2018)Costa, Clain, Loub{\`e}re, and
  Machado]{costa2018very}
R.~Costa, S.~Clain, R.~Loub{\`e}re, and G.~J. Machado.
\newblock Very high-order accurate finite volume scheme on curved boundaries
  for the two-dimensional steady-state convection--diffusion equation with
  {D}irichlet condition.
\newblock \emph{Applied Mathematical Modelling}, 54:\penalty0 752--767, 2018.

\bibitem[Costa et~al.(2019)Costa, N{\'o}brega, Clain, Machado, and
  Loub{\`e}re]{costa2019very}
R.~Costa, J.~M. N{\'o}brega, S.~Clain, G.~J. Machado, and R.~Loub{\`e}re.
\newblock Very high-order accurate finite volume scheme for the
  convection-diffusion equation with general boundary conditions on arbitrary
  curved boundaries.
\newblock \emph{International Journal for Numerical Methods in Engineering},
  117\penalty0 (2):\penalty0 188--220, 2019.

\bibitem[Dutt et~al.(2000)Dutt, Greengard, and Rokhlin]{dutt2000spectral}
A.~Dutt, L.~Greengard, and V.~Rokhlin.
\newblock Spectral deferred correction methods for ordinary differential
  equations.
\newblock \emph{BIT Numerical Mathematics}, 40\penalty0 (2):\penalty0 241--266,
  2000.

\bibitem[Fadlun et~al.(2000)Fadlun, Verzicco, Orlandi, and
  Mohd-Yusof]{fadlun2000combined}
E.~A. Fadlun, R.~Verzicco, P.~Orlandi, and J.~Mohd-Yusof.
\newblock Combined immersed-boundary finite-difference methods for
  three-dimensional complex flow simulations.
\newblock \emph{Journal of Computational Physics}, 161\penalty0 (1):\penalty0
  35--60, 2000.

\bibitem[Fedkiw et~al.(1999)Fedkiw, Aslam, Merriman, and Osher]{fedkiw1999non}
R.~P. Fedkiw, T.~Aslam, B.~Merriman, and S.~Osher.
\newblock {A non-oscillatory {E}ulerian approach to interfaces in multimaterial
  flows (the ghost fluid method)}.
\newblock \emph{Journal of Computational Physics}, 152\penalty0 (2):\penalty0
  457--492, 1999.

\bibitem[Gassner(2013)]{gassner2013skew}
G.~J. Gassner.
\newblock A skew-symmetric discontinuous {G}alerkin spectral element
  discretization and its relation to {SBP-SAT} finite difference methods.
\newblock \emph{SIAM Journal on Scientific Computing}, 35\penalty0
  (3):\penalty0 A1233--A1253, 2013.

\bibitem[Giuliani(2022)]{giuliani2022two}
A.~Giuliani.
\newblock {A two-dimensional stabilized discontinuous {G}alerkin method on
  curvilinear embedded boundary grids}.
\newblock \emph{SIAM Journal on Scientific Computing}, 44\penalty0
  (1):\penalty0 A389--A415, 2022.

\bibitem[Gorgi et~al.(2026)Gorgi, Antonelli, Cornejo, Scovazzi, and
  Rossi]{gorgi2026isogeometric}
A.~Gorgi, N.~Antonelli, A.~Cornejo, G.~Scovazzi, and R.~Rossi.
\newblock The isogeometric {G}ap-shifted boundary method.
\newblock \emph{Computer Methods in Applied Mechanics and Engineering},
  462:\penalty0 119293, 2026.

\bibitem[Hansbo and Hansbo(2002)]{hansbo2002unfitted}
A.~Hansbo and P.~Hansbo.
\newblock An unfitted finite element method, based on {N}itsche's method, for
  elliptic interface problems.
\newblock \emph{Computer Methods in Applied Mechanics and Engineering},
  191\penalty0 (47-48):\penalty0 5537--5552, 2002.

\bibitem[Jung and Perrier(2022)]{jung2022steady}
J.~Jung and V.~Perrier.
\newblock Steady low {M}ach number flows: identification of the spurious mode
  and filtering method.
\newblock \emph{Journal of Computational Physics}, 468:\penalty0 111462, 2022.
\newblock \doi{10.1016/j.jcp.2022.111462}.

\bibitem[Krivodonova and Berger(2006)]{krivodonova2006high}
L.~Krivodonova and M.~Berger.
\newblock High-order accurate implementation of solid wall boundary conditions
  in curved geometries.
\newblock \emph{Journal of Computational Physics}, 211\penalty0 (2):\penalty0
  492--512, 2006.

\bibitem[Li et~al.(2026)Li, Mittal, Collins, and Scovazzi]{li2026gap}
K.~Li, T.~Mittal, J.~H. Collins, and G.~Scovazzi.
\newblock The gap-shifted boundary method for solid mechanics problems: linear
  elasticity, hyperelasticity and elastoplasticity.
\newblock \emph{Computer Methods in Applied Mechanics and Engineering},
  462:\penalty0 119282, 2026.

\bibitem[Li et~al.(2016)Li, Shu, and Zhang]{li2016stability}
T.~Li, C.-W. Shu, and M.~Zhang.
\newblock Stability analysis of the inverse {L}ax--{W}endroff boundary
  treatment for high order upwind-biased finite difference schemes.
\newblock \emph{Journal of Computational and Applied Mathematics},
  299:\penalty0 140--158, 2016.

\bibitem[Main and Scovazzi(2018)]{main2018shifted}
A.~Main and G.~Scovazzi.
\newblock The shifted boundary method for embedded domain computations. {P}art
  {I}: {P}oisson and {S}tokes problems.
\newblock \emph{Journal of Computational Physics}, 372:\penalty0 972--995,
  2018.

\bibitem[May and Berger(2017)]{may2017explicit}
S.~May and M.~Berger.
\newblock An explicit implicit scheme for cut cells in embedded boundary
  meshes.
\newblock \emph{Journal of Scientific Computing}, 71\penalty0 (3):\penalty0
  919--943, 2017.

\bibitem[Micalizzi and Torlo(2024)]{micalizzi2024new}
L.~Micalizzi and D.~Torlo.
\newblock A new efficient explicit deferred correction framework: analysis and
  applications to hyperbolic {PDEs} and adaptivity.
\newblock \emph{Communications on Applied Mathematics and Computation},
  6\penalty0 (3):\penalty0 1629--1664, 2024.

\bibitem[Michalak and Ollivier-Gooch(2009)]{michalak2009unstructured}
C.~Michalak and C.~Ollivier-Gooch.
\newblock Unstructured high-order accurate finite-volume solutions of the
  {N}avier--{S}tokes equations.
\newblock In \emph{47th AIAA Aerospace Sciences Meeting including The New
  Horizons Forum and Aerospace Exposition}, page 954, 2009.

\bibitem[Mittal and Iaccarino(2005)]{mittal2005immersed}
R.~Mittal and G.~Iaccarino.
\newblock Immersed boundary methods.
\newblock \emph{Annual Review of Fluid Mechanics}, 37\penalty0 (1):\penalty0
  239--261, 2005.

\bibitem[Monasse et~al.(2012)Monasse, Daru, Mariotti, Piperno, and
  Tenaud]{monasse2012conservative}
L.~Monasse, V.~Daru, C.~Mariotti, S.~Piperno, and C.~Tenaud.
\newblock A conservative coupling algorithm between a compressible flow and a
  rigid body using an embedded boundary method.
\newblock \emph{Journal of Computational Physics}, 231\penalty0 (7):\penalty0
  2977--2994, 2012.

\bibitem[Moxey et~al.(2015)Moxey, Green, Sherwin, and
  Peir{\'o}]{moxey2015isoparametric}
D.~Moxey, M.~Green, S.~Sherwin, and J.~Peir{\'o}.
\newblock An isoparametric approach to high-order curvilinear boundary-layer
  meshing.
\newblock \emph{Computer Methods in Applied Mechanics and Engineering},
  283:\penalty0 636--650, 2015.

\bibitem[Mulder and van Leer(1985)]{mulder1985experiments}
W.~A. Mulder and B.~van Leer.
\newblock Experiments with implicit upwind methods for the {E}uler equations.
\newblock \emph{Journal of Computational Physics}, 59\penalty0 (2):\penalty0
  232--246, 1985.

\bibitem[Ollivier-Gooch and Van~Altena(2002)]{ollivier2002high}
C.~Ollivier-Gooch and M.~Van~Altena.
\newblock A high-order-accurate unstructured mesh finite-volume scheme for the
  advection--diffusion equation.
\newblock \emph{Journal of Computational Physics}, 181\penalty0 (2):\penalty0
  729--752, 2002.

\bibitem[Persson and Peraire(2009)]{persson2009curved}
P.-O. Persson and J.~Peraire.
\newblock {Curved mesh generation and mesh refinement using {L}agrangian solid
  mechanics}.
\newblock In \emph{47th AIAA Aerospace Sciences Meeting including the New
  Horizons Forum and Aerospace Exposition}, page 949, 2009.

\bibitem[Peskin(1977)]{peskin1977numerical}
C.~S. Peskin.
\newblock Numerical analysis of blood flow in the heart.
\newblock \emph{Journal of Computational Physics}, 25\penalty0 (3):\penalty0
  220--252, 1977.

\bibitem[Peskin(2002)]{peskin2002immersed}
C.~S. Peskin.
\newblock The immersed boundary method.
\newblock \emph{Acta Numerica}, 11:\penalty0 479--517, 2002.

\bibitem[Qin and Krivodonova(2013)]{qin2013discontinuous}
R.~Qin and L.~Krivodonova.
\newblock A discontinuous {G}alerkin method for solutions of the {E}uler
  equations on {C}artesian grids with embedded geometries.
\newblock \emph{Journal of Computational Science}, 4\penalty0 (1-2):\penalty0
  24--35, 2013.

\bibitem[Reed and Hill(1973)]{reed1973triangular}
W.~H. Reed and T.~R. Hill.
\newblock Triangular mesh methods for the neutron transport equation.
\newblock Technical report, Los Alamos Scientific Lab., N. Mex.(USA), 1973.

\bibitem[Saad and Schultz(1986)]{saad1986gmres}
Y.~Saad and M.~H. Schultz.
\newblock {GMRES}: A generalized minimal residual algorithm for solving
  nonsymmetric linear systems.
\newblock \emph{SIAM Journal on Scientific and Statistical Computing},
  7\penalty0 (3):\penalty0 856--869, 1986.

\bibitem[Santos et~al.(2024)Santos, Ara{\'u}jo, Barbeiro, Clain, Costa, and
  Machado]{santos2024very}
M.~Santos, A.~Ara{\'u}jo, S.~Barbeiro, S.~Clain, R.~Costa, and G.~J. Machado.
\newblock {Very high-order accurate discontinuous {G}alerkin method for curved
  boundaries with polygonal meshes}.
\newblock \emph{Journal of Scientific Computing}, 100\penalty0 (3):\penalty0
  66, 2024.

\bibitem[Song et~al.(2018)Song, Main, Scovazzi, and Ricchiuto]{song2018shifted}
T.~Song, A.~Main, G.~Scovazzi, and M.~Ricchiuto.
\newblock The shifted boundary method for hyperbolic systems: Embedded domain
  computations of linear waves and shallow water flows.
\newblock \emph{Journal of Computational Physics}, 369:\penalty0 45--79, 2018.

\bibitem[Tan and Shu(2010)]{tan2010inverse}
S.~Tan and C.-W. Shu.
\newblock Inverse {L}ax-{W}endroff procedure for numerical boundary conditions
  of conservation laws.
\newblock \emph{Journal of Computational Physics}, 229\penalty0 (21):\penalty0
  8144--8166, 2010.

\bibitem[Tan et~al.(2012)Tan, Wang, Shu, and Ning]{tan2012efficient}
S.~Tan, C.~Wang, C.-W. Shu, and J.~Ning.
\newblock Efficient implementation of high order inverse {L}ax--{W}endroff
  boundary treatment for conservation laws.
\newblock \emph{Journal of Computational Physics}, 231\penalty0 (6):\penalty0
  2510--2527, 2012.

\bibitem[Visbech et~al.(2025)Visbech, Engsig-Karup, and
  Ricchiuto]{visbech2025spectral}
J.~Visbech, A.~P. Engsig-Karup, and M.~Ricchiuto.
\newblock A spectral element solution of the {P}oisson equation with shifted
  boundary polynomial corrections: influence of the surrogate to true boundary
  mapping and an asymptotically preserving {R}obin formulation.
\newblock \emph{Journal of Scientific Computing}, 102\penalty0 (1):\penalty0
  11, 2025.

\bibitem[Wang et~al.(2013)Wang, Fidkowski, Abgrall, Bassi, Caraeni, Cary,
  Deconinck, Hartmann, Hillewaert, Huynh, et~al.]{wang2013high}
Z.~J. Wang, K.~Fidkowski, R.~Abgrall, F.~Bassi, D.~Caraeni, A.~Cary,
  H.~Deconinck, R.~Hartmann, K.~Hillewaert, H.~T. Huynh, et~al.
\newblock {High-order CFD methods: current status and perspective}.
\newblock \emph{International Journal for Numerical Methods in Fluids},
  72\penalty0 (8):\penalty0 811--845, 2013.

\bibitem[Yang et~al.(2025)Yang, Li, Jiang, Shu, and Zhang]{yang2025inverse}
L.~Yang, S.~Li, Y.~Jiang, C.-W. Shu, and M.~Zhang.
\newblock Inverse lax-wendroff boundary treatment of discontinuous galerkin
  method for 1d conservation laws.
\newblock \emph{Communications on Applied Mathematics and Computation},
  7\penalty0 (2):\penalty0 796--826, 2025.

\bibitem[Zeng et~al.(2025)Zeng, Song, and Scovazzi]{ZENG2025113512}
X.~Zeng, T.~Song, and G.~Scovazzi.
\newblock {A shifted boundary method for the compressible {E}uler equations}.
\newblock \emph{Journal of Computational Physics}, 520:\penalty0 113512, 2025.

\end{thebibliography}

\end{document}